\overfullrule=0pt
\magnification=1200

\hsize=12cm    
\vsize=19cm
\hoffset=0.5cm

\catcode`@=11
\def\@height{height}
\def\@depth{depth}
\def\@width{width}

\newcount\@tempcnta
\newcount\@tempcntb

\newdimen\@tempdima
\newdimen\@tempdimb

\newbox\@tempboxa

\def\@ifnextchar#1#2#3{\let\@tempe #1\def\@tempa{#2}\def\@tempb{#3}\futurelet
    \@tempc\@ifnch}
\def\@ifnch{\ifx \@tempc \@sptoken \let\@tempd\@xifnch
      \else \ifx \@tempc \@tempe\let\@tempd\@tempa\else\let\@tempd\@tempb\fi
      \fi \@tempd}
\def\@ifstar#1#2{\@ifnextchar *{\def\@tempa*{#1}\@tempa}{#2}}

\def\@whilenoop#1{}
\def\@whilenum#1\do #2{\ifnum #1\relax #2\relax\@iwhilenum{#1\relax 
     #2\relax}\fi}
\def\@iwhilenum#1{\ifnum #1\let\@nextwhile=\@iwhilenum 
         \else\let\@nextwhile=\@whilenoop\fi\@nextwhile{#1}}

\def\@whiledim#1\do #2{\ifdim #1\relax#2\@iwhiledim{#1\relax#2}\fi}
\def\@iwhiledim#1{\ifdim #1\let\@nextwhile=\@iwhiledim 
        \else\let\@nextwhile=\@whilenoop\fi\@nextwhile{#1}}

\newdimen\@wholewidth
\newdimen\@halfwidth
\newdimen\unitlength \unitlength =1pt
\newbox\@picbox
\newdimen\@picht

\def\@nnil{\@nil}
\def\@empty{}
\def\@fornoop#1\@@#2#3{}

\def\@for#1:=#2\do#3{\edef\@fortmp{#2}\ifx\@fortmp\@empty \else
    \expandafter\@forloop#2,\@nil,\@nil\@@#1{#3}\fi}

\def\@forloop#1,#2,#3\@@#4#5{\def#4{#1}\ifx #4\@nnil \else
       #5\def#4{#2}\ifx #4\@nnil \else#5\@iforloop #3\@@#4{#5}\fi\fi}

\def\@iforloop#1,#2\@@#3#4{\def#3{#1}\ifx #3\@nnil 
       \let\@nextwhile=\@fornoop \else
      #4\relax\let\@nextwhile=\@iforloop\fi\@nextwhile#2\@@#3{#4}}

\def\@tfor#1:=#2\do#3{\xdef\@fortmp{#2}\ifx\@fortmp\@empty \else
    \@tforloop#2\@nil\@nil\@@#1{#3}\fi}
\def\@tforloop#1#2\@@#3#4{\def#3{#1}\ifx #3\@nnil 
       \let\@nextwhile=\@fornoop \else
      #4\relax\let\@nextwhile=\@tforloop\fi\@nextwhile#2\@@#3{#4}}

\def\@makepicbox(#1,#2){\leavevmode\@ifnextchar 
   [{\@imakepicbox(#1,#2)}{\@imakepicbox(#1,#2)[]}}

\long\def\@imakepicbox(#1,#2)[#3]#4{\vbox to#2\unitlength
   {\let\mb@b\vss \let\mb@l\hss\let\mb@r\hss
    \let\mb@t\vss
    \@tfor\@tempa :=#3\do{\expandafter\let
        \csname mb@\@tempa\endcsname\relax}%
\mb@t\hbox to #1\unitlength{\mb@l #4\mb@r}\mb@b}}

\def\picture(#1,#2){\@ifnextchar({\@picture(#1,#2)}{\@picture(#1,#2)(0,0)}}

\def\@picture(#1,#2)(#3,#4){\@picht #2\unitlength
\setbox\@picbox\hbox to #1\unitlength\bgroup 
\hskip -#3\unitlength \lower #4\unitlength \hbox\bgroup\ignorespaces}

\def\endpicture{\egroup\hss\egroup\ht\@picbox\@picht
\dp\@picbox\z@\leavevmode\box\@picbox}

\long\def\put(#1,#2)#3{\@killglue\raise#2\unitlength\hbox to \z@{\kern
#1\unitlength #3\hss}\ignorespaces}

\long\def\multiput(#1,#2)(#3,#4)#5#6{\@killglue\@multicnt=#5\relax
\@xdim=#1\unitlength
\@ydim=#2\unitlength
\@whilenum \@multicnt > 0\do
{\raise\@ydim\hbox to \z@{\kern
\@xdim #6\hss}\advance\@multicnt \m@ne\advance\@xdim
#3\unitlength\advance\@ydim #4\unitlength}\ignorespaces}

\def\@killglue{\unskip\@whiledim \lastskip >\z@\do{\unskip}}

\def\thinlines{\let\@linefnt\tenln \let\@circlefnt\tencirc
  \@wholewidth\fontdimen8\tenln \@halfwidth .5\@wholewidth}
\def\thicklines{\let\@linefnt\tenlnw \let\@circlefnt\tencircw
  \@wholewidth\fontdimen8\tenlnw \@halfwidth .5\@wholewidth}

\def\linethickness#1{\@wholewidth #1\relax \@halfwidth .5\@wholewidth}

\def\shortstack{\@ifnextchar[{\@shortstack}{\@shortstack[c]}}

\def\@shortstack[#1]{\leavevmode
\vbox\bgroup\baselineskip-1pt\lineskip 3pt\let\mb@l\hss
\let\mb@r\hss \expandafter\let\csname mb@#1\endcsname\relax
\let\\\@stackcr\@ishortstack}

\def\@ishortstack#1{\halign{\mb@l ##\unskip\mb@r\cr #1\crcr}\egroup}

\def\@stackcr{\@ifstar{\@ixstackcr}{\@ixstackcr}}
\def\@ixstackcr{\@ifnextchar[{\@istackcr}{\cr\ignorespaces}}

\def\@istackcr[#1]{\cr\noalign{\vskip #1}\ignorespaces}

\newif\if@negarg

\def\droite(#1,#2)#3{\@xarg #1\relax \@yarg #2\relax
\@linelen=#3\unitlength
\ifnum\@xarg =0 \@vline 
  \else \ifnum\@yarg =0 \@hline \else \@sline\fi
\fi}

\def\@sline{\ifnum\@xarg< 0 \@negargtrue \@xarg -\@xarg \@yyarg -\@yarg
  \else \@negargfalse \@yyarg \@yarg \fi
\ifnum \@yyarg >0 \@tempcnta\@yyarg \else \@tempcnta -\@yyarg \fi
\ifnum\@tempcnta>6 \@badlinearg\@tempcnta0 \fi
\ifnum\@xarg>6 \@badlinearg\@xarg 1 \fi
\setbox\@linechar\hbox{\@linefnt\@getlinechar(\@xarg,\@yyarg)}%
\ifnum \@yarg >0 \let\@upordown\raise \@clnht\z@
   \else\let\@upordown\lower \@clnht \ht\@linechar\fi
\@clnwd=\wd\@linechar
\if@negarg \hskip -\wd\@linechar \def\@tempa{\hskip -2\wd\@linechar}\else
     \let\@tempa\relax \fi
\@whiledim \@clnwd <\@linelen \do
  {\@upordown\@clnht\copy\@linechar
   \@tempa
   \advance\@clnht \ht\@linechar
   \advance\@clnwd \wd\@linechar}%
\advance\@clnht -\ht\@linechar
\advance\@clnwd -\wd\@linechar
\@tempdima\@linelen\advance\@tempdima -\@clnwd
\@tempdimb\@tempdima\advance\@tempdimb -\wd\@linechar
\if@negarg \hskip -\@tempdimb \else \hskip \@tempdimb \fi
\multiply\@tempdima \@m
\@tempcnta \@tempdima \@tempdima \wd\@linechar \divide\@tempcnta \@tempdima
\@tempdima \ht\@linechar \multiply\@tempdima \@tempcnta
\divide\@tempdima \@m
\advance\@clnht \@tempdima
\ifdim \@linelen <\wd\@linechar
   \hskip \wd\@linechar
  \else\@upordown\@clnht\copy\@linechar\fi}

\def\@hline{\ifnum \@xarg <0 \hskip -\@linelen \fi
\vrule \@height \@halfwidth \@depth \@halfwidth \@width \@linelen
\ifnum \@xarg <0 \hskip -\@linelen \fi}

\def\@getlinechar(#1,#2){\@tempcnta#1\relax\multiply\@tempcnta 8
\advance\@tempcnta -9 \ifnum #2>0 \advance\@tempcnta #2\relax\else
\advance\@tempcnta -#2\relax\advance\@tempcnta 64 \fi
\char\@tempcnta}

\def\vector(#1,#2)#3{\@xarg #1\relax \@yarg #2\relax
\@tempcnta \ifnum\@xarg<0 -\@xarg\else\@xarg\fi
\ifnum\@tempcnta<5\relax
\@linelen=#3\unitlength
\ifnum\@xarg =0 \@vvector 
  \else \ifnum\@yarg =0 \@hvector \else \@svector\fi
\fi
\else\@badlinearg\fi}

\def\@hvector{\@hline\hbox to 0pt{\@linefnt 
\ifnum \@xarg <0 \@getlarrow(1,0)\hss\else
    \hss\@getrarrow(1,0)\fi}}

\def\@vvector{\ifnum \@yarg <0 \@downvector \else \@upvector \fi}

\def\@svector{\@sline
\@tempcnta\@yarg \ifnum\@tempcnta <0 \@tempcnta=-\@tempcnta\fi
\ifnum\@tempcnta <5
  \hskip -\wd\@linechar
  \@upordown\@clnht \hbox{\@linefnt  \if@negarg 
  \@getlarrow(\@xarg,\@yyarg) \else \@getrarrow(\@xarg,\@yyarg) \fi}%
\else\@badlinearg\fi}

\def\@getlarrow(#1,#2){\ifnum #2 =\z@ \@tempcnta='33\else
\@tempcnta=#1\relax\multiply\@tempcnta \sixt@@n \advance\@tempcnta
-9 \@tempcntb=#2\relax\multiply\@tempcntb \tw@
\ifnum \@tempcntb >0 \advance\@tempcnta \@tempcntb\relax
\else\advance\@tempcnta -\@tempcntb\advance\@tempcnta 64
\fi\fi\char\@tempcnta}

\def\@getrarrow(#1,#2){\@tempcntb=#2\relax
\ifnum\@tempcntb < 0 \@tempcntb=-\@tempcntb\relax\fi
\ifcase \@tempcntb\relax \@tempcnta='55 \or 
\ifnum #1<3 \@tempcnta=#1\relax\multiply\@tempcnta
24 \advance\@tempcnta -6 \else \ifnum #1=3 \@tempcnta=49
\else\@tempcnta=58 \fi\fi\or 
\ifnum #1<3 \@tempcnta=#1\relax\multiply\@tempcnta
24 \advance\@tempcnta -3 \else \@tempcnta=51\fi\or 
\@tempcnta=#1\relax\multiply\@tempcnta
\sixt@@n \advance\@tempcnta -\tw@ \else
\@tempcnta=#1\relax\multiply\@tempcnta
\sixt@@n \advance\@tempcnta 7 \fi\ifnum #2<0 \advance\@tempcnta 64 \fi
\char\@tempcnta}

\def\@vline{\ifnum \@yarg <0 \@downline \else \@upline\fi}

\def\@upline{\hbox to \z@{\hskip -\@halfwidth \vrule \@width \@wholewidth
   \@height \@linelen \@depth \z@\hss}}

\def\@downline{\hbox to \z@{\hskip -\@halfwidth \vrule \@width \@wholewidth
   \@height \z@ \@depth \@linelen \hss}}

\def\@upvector{\@upline\setbox\@tempboxa\hbox{\@linefnt\char'66}\raise 
     \@linelen \hbox to\z@{\lower \ht\@tempboxa\box\@tempboxa\hss}}

\def\@downvector{\@downline\lower \@linelen
      \hbox to \z@{\@linefnt\char'77\hss}}

\def\dashbox#1(#2,#3){\leavevmode\hbox to \z@{\baselineskip \z@%
\lineskip \z@%
\@dashdim=#2\unitlength%
\@dashcnt=\@dashdim \advance\@dashcnt 200
\@dashdim=#1\unitlength\divide\@dashcnt \@dashdim
\ifodd\@dashcnt\@dashdim=\z@%
\advance\@dashcnt \@ne \divide\@dashcnt \tw@ 
\else \divide\@dashdim \tw@ \divide\@dashcnt \tw@
\advance\@dashcnt \m@ne
\setbox\@dashbox=\hbox{\vrule \@height \@halfwidth \@depth \@halfwidth
\@width \@dashdim}\put(0,0){\copy\@dashbox}%
\put(0,#3){\copy\@dashbox}%
\put(#2,0){\hskip-\@dashdim\copy\@dashbox}%
\put(#2,#3){\hskip-\@dashdim\box\@dashbox}%
\multiply\@dashdim 3 
\fi
\setbox\@dashbox=\hbox{\vrule \@height \@halfwidth \@depth \@halfwidth
\@width #1\unitlength\hskip #1\unitlength}\@tempcnta=0
\put(0,0){\hskip\@dashdim \@whilenum \@tempcnta <\@dashcnt
\do{\copy\@dashbox\advance\@tempcnta \@ne }}\@tempcnta=0
\put(0,#3){\hskip\@dashdim \@whilenum \@tempcnta <\@dashcnt
\do{\copy\@dashbox\advance\@tempcnta \@ne }}%
\@dashdim=#3\unitlength%
\@dashcnt=\@dashdim \advance\@dashcnt 200
\@dashdim=#1\unitlength\divide\@dashcnt \@dashdim
\ifodd\@dashcnt \@dashdim=\z@%
\advance\@dashcnt \@ne \divide\@dashcnt \tw@
\else
\divide\@dashdim \tw@ \divide\@dashcnt \tw@
\advance\@dashcnt \m@ne
\setbox\@dashbox\hbox{\hskip -\@halfwidth
\vrule \@width \@wholewidth 
\@height \@dashdim}\put(0,0){\copy\@dashbox}%
\put(#2,0){\copy\@dashbox}%
\put(0,#3){\lower\@dashdim\copy\@dashbox}%
\put(#2,#3){\lower\@dashdim\copy\@dashbox}%
\multiply\@dashdim 3
\fi
\setbox\@dashbox\hbox{\vrule \@width \@wholewidth 
\@height #1\unitlength}\@tempcnta0
\put(0,0){\hskip -\@halfwidth \vbox{\@whilenum \@tempcnta < \@dashcnt
\do{\vskip #1\unitlength\copy\@dashbox\advance\@tempcnta \@ne }%
\vskip\@dashdim}}\@tempcnta0
\put(#2,0){\hskip -\@halfwidth \vbox{\@whilenum \@tempcnta< \@dashcnt
\relax\do{\vskip #1\unitlength\copy\@dashbox\advance\@tempcnta \@ne }%
\vskip\@dashdim}}}\@makepicbox(#2,#3)}

\newif\if@ovt 
\newif\if@ovb 
\newif\if@ovl 
\newif\if@ovr 
\newdimen\@ovxx
\newdimen\@ovyy
\newdimen\@ovdx
\newdimen\@ovdy
\newdimen\@ovro
\newdimen\@ovri

\def\@getcirc#1{\@tempdima #1\relax \advance\@tempdima 2pt\relax
  \@tempcnta\@tempdima
  \@tempdima 4pt\relax \divide\@tempcnta\@tempdima
  \ifnum \@tempcnta > 10\relax \@tempcnta 10\relax\fi
  \ifnum \@tempcnta >\z@ \advance\@tempcnta\m@ne
    \else \@warning{Oval too small}\fi
  \multiply\@tempcnta 4\relax
  \setbox \@tempboxa \hbox{\@circlefnt
  \char \@tempcnta}\@tempdima \wd \@tempboxa}

\def\@put#1#2#3{\raise #2\hbox to \z@{\hskip #1#3\hss}}

\def\oval(#1,#2){\@ifnextchar[{\@oval(#1,#2)}{\@oval(#1,#2)[]}}

\def\@oval(#1,#2)[#3]{\begingroup\boxmaxdepth \maxdimen
  \@ovttrue \@ovbtrue \@ovltrue \@ovrtrue
  \@tfor\@tempa :=#3\do{\csname @ov\@tempa false\endcsname}\@ovxx
  #1\unitlength \@ovyy #2\unitlength
  \@tempdimb \ifdim \@ovyy >\@ovxx \@ovxx\else \@ovyy \fi
  \advance \@tempdimb -2pt\relax  
  \@getcirc \@tempdimb
  \@ovro \ht\@tempboxa \@ovri \dp\@tempboxa
  \@ovdx\@ovxx \advance\@ovdx -\@tempdima \divide\@ovdx \tw@
  \@ovdy\@ovyy \advance\@ovdy -\@tempdima \divide\@ovdy \tw@
  \@circlefnt \setbox\@tempboxa
  \hbox{\if@ovr \@ovvert32\kern -\@tempdima \fi
  \if@ovl \kern \@ovxx \@ovvert01\kern -\@tempdima \kern -\@ovxx \fi
  \if@ovt \@ovhorz \kern -\@ovxx \fi
  \if@ovb \raise \@ovyy \@ovhorz \fi}\advance\@ovdx\@ovro
  \advance\@ovdy\@ovro \ht\@tempboxa\z@ \dp\@tempboxa\z@
  \@put{-\@ovdx}{-\@ovdy}{\box\@tempboxa}%
  \endgroup}

\def\@ovvert#1#2{\vbox to \@ovyy{%
    \if@ovb \@tempcntb \@tempcnta \advance \@tempcntb by #1\relax
      \kern -\@ovro \hbox{\char \@tempcntb}\nointerlineskip
    \else \kern \@ovri \kern \@ovdy \fi
    \leaders\vrule width \@wholewidth\vfil \nointerlineskip
    \if@ovt \@tempcntb \@tempcnta \advance \@tempcntb by #2\relax
      \hbox{\char \@tempcntb}%
    \else \kern \@ovdy \kern \@ovro \fi}}

\def\@ovhorz{\hbox to \@ovxx{\kern \@ovro
    \if@ovr \else \kern \@ovdx \fi
    \leaders \hrule height \@wholewidth \hfil
    \if@ovl \else \kern \@ovdx \fi
    \kern \@ovri}}

\def\circle{\@ifstar{\@dot}{\@circle}}
\def\@circle#1{\begingroup \boxmaxdepth \maxdimen \@tempdimb #1\unitlength
   \ifdim \@tempdimb >15.5pt\relax \@getcirc\@tempdimb
      \@ovro\ht\@tempboxa 
     \setbox\@tempboxa\hbox{\@circlefnt
      \advance\@tempcnta\tw@ \char \@tempcnta
      \advance\@tempcnta\m@ne \char \@tempcnta \kern -2\@tempdima
      \advance\@tempcnta\tw@
      \raise \@tempdima \hbox{\char\@tempcnta}\raise \@tempdima
        \box\@tempboxa}\ht\@tempboxa\z@ \dp\@tempboxa\z@
      \@put{-\@ovro}{-\@ovro}{\box\@tempboxa}%
   \else  \@circ\@tempdimb{96}\fi\endgroup}

\def\@dot#1{\@tempdimb #1\unitlength \@circ\@tempdimb{112}}

\def\@circ#1#2{\@tempdima #1\relax \advance\@tempdima .5pt\relax
   \@tempcnta\@tempdima \@tempdima 1pt\relax
   \divide\@tempcnta\@tempdima 
   \ifnum\@tempcnta > 15\relax \@tempcnta 15\relax \fi    
   \ifnum \@tempcnta >\z@ \advance\@tempcnta\m@ne\fi
   \advance\@tempcnta #2\relax
   \@circlefnt \char\@tempcnta}

\font\tenln line10
\font\tencirc lcircle10
\font\tenlnw linew10
\font\tencircw lcirclew10

\thinlines   

\newcount\@xarg
\newcount\@yarg
\newcount\@yyarg
\newcount\@multicnt 
\newdimen\@xdim
\newdimen\@ydim
\newbox\@linechar
\newdimen\@linelen
\newdimen\@clnwd
\newdimen\@clnht
\newdimen\@dashdim
\newbox\@dashbox
\newcount\@dashcnt
\catcode`@=12

\def\g{\mathop{\dashv}\nolimits}

\def\d{\mathop{\vdash}\nolimits}


\def\gg{{\bf g}}
\def\h{{\bf h}}
\def\m{{\bf m}}

\def\N{\noindent}
\def\S{\smallskip \par}
\def\M{\medskip \par}
\def\B{\bigskip \par}
\def\BB{\bigskip \bigskip \par}

\def\aa{\alpha}

\def\cc{\gamma}

\def\dd{\delta}
\def\DD{\Delta}
\def\ee{\epsilon}

\def\ss{\sigma}
\def\SS{\Sigma}
\def\oo{\omega}

\def\NN{{\bf N}}
\def\ZZ{{\bf Z}}

\def\sqr#1#2{{\vcenter{\vbox{\hrule height.#2pt
\hbox{\vrule width .#2pt height#1pt \kern#1pt
\vrule width.#2pt}
\hrule height.#2pt}}}}
\def\square{\mathchoice\sqr64\sqr64\sqr{4.2}3\sqr33}

\def \Im{\mathop{\rm Im\,}\nolimits}
\def \Ker{\mathop{\rm Ker\,}\nolimits}

\def \Hom{\mathop{\rm Hom\,}\nolimits}

\def \Id{\mathop{\rm Id}\nolimits}

\def \ad{\mathop{\rm ad\,}\nolimits}

\def \Ind{\mathop{\rm Ind\,}\nolimits}

\def\t{\otimes }

\def \r{\mathop{\rm \rightarrow}\nolimits}


\def\proj{\mathop{\rightarrow \!\!\!\!\!\!\! \rightarrow}}



\def\tv{\vrule height 12pt depth 5pt}
\def\newfunction #1{\expandafter\def\csname #1\endcsname
{\mathop{\rm #1\kern 0pt}\nolimits}}

\def \P{\mathop{\cal P}\nolimits}

\def \Q{\mathop{\cal Q}\nolimits}

\newfunction {Tor}
\newfunction {id}
\newfunction {Tot}
\newfunction {Mor}
\newfunction {Aut}
\newfunction {Mod}
\newfunction {Vect}
\newfunction {Com}
\newfunction {Leib}
\newfunction {Id}
\newfunction {Hom}
\newfunction {Lie}
\newfunction {Pois}
\newfunction {sgn}
\newfunction {As}
\newfunction {alg}
\newfunction {Iso}
\newfunction {Ens}
\newfunction {In}
\newfunction {mod}
\newfunction {Dias}
\newfunction {Di}
\newfunction {alg}

\def\m{\mathop{^{-1}}\nolimits}

\catcode`@=11
\font\@linefnt linew10 at 2.4pt
\catcode`@=12

\def\arbreun{\kern-0.4ex
\hbox{\unitlength=.25pt
\picture(60,40)(0,0)
\put(30,0){\droite(0,1){20}}
\put(30,20){\droite(-1,1){30}}
\put(30,20){\droite(1,1){30}}
\put(20,30){\droite(1,1){20}}
\put(10,40){\droite(1,1){10}}
\endpicture}\kern 0.4ex}

\def\arbredeux{\kern-0.4ex
\hbox{\unitlength=.25pt
\picture(60,40)(0,0)
\put(30,0){\droite(0,1){20}}
\put(30,20){\droite(-1,1){30}}
\put(30,20){\droite(1,1){30}}
\put(20,30){\droite(1,1){20}}
\put(30,40){\droite(-1,1){10}}
\endpicture}\kern 0.4ex}

\def\arbretrois{\kern-0.4ex
\hbox{\unitlength=.25pt
\picture(60,40)(0,0)
\put(30,0){\droite(0,1){20}}
\put(30,20){\droite(-1,1){30}}
\put(30,20){\droite(1,1){30}}
\put(50,40){\droite(-1,1){10}}
\put(10,40){\droite(1,1){10}}
\endpicture}\kern 0.4ex}

\def\arbretroisxyz{\kern-0.4ex
\hbox{\unitlength=.25pt
\picture(90,75)(0,0)
\put(45,0){\droite(0,1){30}}
\put(45,30){\droite(-1,1){45}}
\put(45,30){\droite(1,1){45}}
\put(75,60){\droite(-1,1){15}}
\put(15,60){\droite(1,1){15}}
\put(5,75){$x$}
\put(35,75){$y$}
\put(65,75){$z$}
\endpicture}\kern 0.4ex}

\def\arbrequatre{\kern-0.4ex
\hbox{\unitlength=.25pt
\picture(60,40)(0,0)
\put(30,0){\droite(0,1){20}}
\put(30,20){\droite(-1,1){30}}
\put(30,20){\droite(1,1){30}}
\put(40,30){\droite(-1,1){20}}
\put(30,40){\droite(1,1){10}}
\endpicture}\kern 0.4ex}

\def\arbrecinq{\kern-0.4ex
\hbox{\unitlength=.25pt
\picture(60,40)(0,0)
\put(30,0){\droite(0,1){20}}
\put(30,20){\droite(-1,1){30}}
\put(30,20){\droite(1,1){30}}
\put(40,30){\droite(-1,1){20}}
\put(50,40){\droite(-1,1){10}}
\endpicture}\kern 0.4ex}

\def\arbreA{\kern-0.4ex
\hbox{\unitlength=.25pt
\picture(60,40)(0,0)
\put(30,0){\droite(0,1){20}}
\put(30,20){\droite(-1,1){30}}
\put(30,20){\droite(1,1){30}}
\endpicture}\kern 0.4ex}

\def\arbreAx{\kern-0.4ex
\hbox{\unitlength=.25pt
\picture(60,60)(0,0)
\put(30,0){\droite(0,1){20}}
\put(30,20){\droite(-1,1){30}}
\put(30,20){\droite(1,1){30}}
\put(25, 60){$x$}
\endpicture}\kern 0.4ex}

\def\arbreAz{\kern-0.4ex
\hbox{\unitlength=.25pt
\picture(60,60)(0,0)
\put(30,0){\droite(0,1){20}}
\put(30,20){\droite(-1,1){30}}
\put(30,20){\droite(1,1){30}}
\put(25, 60){$z$}
\endpicture}\kern 0.4ex}

\def\arbreAxgy{\kern-0.4ex
\hbox{\unitlength=.25pt
\picture(100,70)(0,0)
\put(50,0){\droite(0,1){20}}
\put(50,20){\droite(-1,1){50}}
\put(50,20){\droite(1,1){50}}
\put(10, 70){$x\g y$}
\endpicture}\kern 0.4ex}

\def\arbreAxdy{\kern-0.4ex
\hbox{\unitlength=.25pt
\picture(100,70)(0,0)
\put(50,0){\droite(0,1){20}}
\put(50,20){\droite(-1,1){50}}
\put(50,20){\droite(1,1){50}}
\put(10, 70){$x\d y$}
\endpicture}\kern 0.4ex}

\def\arbreAbis{\kern-0.4ex
\hbox{\unitlength=.25pt
\picture(60,40)(0,0)
\put(30,0){\droite(0,1){20}}
\put(30,20){\droite(-1,1){30}}
\put(30,20){\droite(1,1){30}}
\put(-20,60){$y_1$}
\put(40,60){$y_2$}
\endpicture}\kern 0.4ex}

\def\arbreB{\kern-0.4ex
\hbox{\unitlength=.25pt
\picture(60,40)(0,0)
\put(30,0){\droite(0,1){20}}
\put(30,20){\droite(-1,1){30}}
\put(30,20){\droite(1,1){30}}
\put(15,35){\droite(1,1){15}}
\endpicture}\kern 0.4ex}

\def\arbreBxy{\kern-0.4ex
\hbox{\unitlength=.25pt
\picture(60,60)(0,0)
\put(30,0){\droite(0,1){20}}
\put(30,20){\droite(-1,1){30}}
\put(30,20){\droite(1,1){30}}
\put(15,35){\droite(1,1){15}}
\put(10,60){$x$}
\put(35,60){$y$}
\endpicture}\kern 0.4ex}

\def\arbreC{\kern-0.4ex
\hbox{\unitlength=.25pt
\picture(60,40)(0,0)
\put(30,0){\droite(0,1){20}}
\put(30,20){\droite(-1,1){30}}
\put(30,20){\droite(1,1){30}}
\put(45,35){\droite(-1,1){15}}
\endpicture}\kern 0.4ex}

\def\arbreCxy{\kern-0.4ex
\hbox{\unitlength=.25pt
\picture(60,40)(0,0)
\put(30,0){\droite(0,1){20}}
\put(30,20){\droite(-1,1){30}}
\put(30,20){\droite(1,1){30}}
\put(45,35){\droite(-1,1){15}}
\put(5,60){$x$}
\put(35,60){$y$}
\endpicture}\kern 0.4ex}

\def\TRE{\kern-0.4ex
\hbox{\unitlength=1pt
\picture(100,90)(0,0)
\put(50,0){\droite(0,1){40}}
\put(50,40){\droite(-1,1){50}}
\put(50,40){\droite(1,1){50}}
\put(30,60){\droite(1,1){30}}
\put(10,80){\droite(1,1){10}}
\put(50,80){\droite(-1,1){10}}
\put(90,80){\droite(-1,1){10}}

\put(55,35){$\dashv$}
\put(25,50){$\vdash$}
\put(5,70){$\dashv$}
\put(50,70){$\dashv$}
\put(90,70){$\vdash$}
\put(0,95){$x_1$}
\put(20,95){$x_2$}
\put(40,95){$x_3$}
\put(60,95){$x_4$}
\put(80,95){$x_5$}
\put(100,95){$x_6$}
\endpicture}\kern 0.4ex}

\def\TREE{\kern-0.4ex
\hbox{\unitlength=1pt
\picture(100,90)(0,0)
\put(50,0){\droite(0,1){40}}
\put(52,0){\droite(0,1){40}}
\put(50,40){\droite(-1,1){50}}
\put(50,42){\droite(-1,1){18}}
\put(50,40){\droite(1,1){50}}
\put(30,60){\droite(1,1){30}}
\put(30,62){\droite(1,1){18}}
\put(10,80){\droite(1,1){10}}
\put(50,80){\droite(-1,1){10}}
\put(50,82){\droite(-1,1){8}}
\put(90,80){\droite(-1,1){10}}

\put(55,35){$\dashv$}
\put(35,55){$\vdash$}
\put(15,75){$\vdash$}
\put(55,75){$\dashv$}
\put(95,75){$\vdash$}
\put(40,95){$i$}
\endpicture}\kern 0.4ex}

\def\TREEE{\kern-0.4ex
\hbox{\unitlength=1pt
\picture(100,90)(0,0)
\put(50,0){\droite(0,1){40}}
\put(52,0){\droite(0,1){40}}
\put(50,40){\droite(-1,1){50}}
\put(50,42){\droite(-1,1){18}}
\put(50,40){\droite(1,1){50}}
\put(30,60){\droite(1,1){30}}
\put(30,62){\droite(1,1){18}}
\put(10,80){\droite(1,1){10}}
\put(50,80){\droite(-1,1){10}}
\put(50,82){\droite(-1,1){8}}
\put(90,80){\droite(-1,1){10}}

\put(55,35){$\dashv$}
\put(25,50){$\vdash$}
\put(5,70){$\dashv$}
\put(50,70){$\dashv$}
\put(90,70){$\vdash$}
\put(40,95){$x_i$}
\endpicture}\kern 0.4ex}

\def\TREEF{\kern-0.4ex
\hbox{\unitlength=1pt
\picture(100,90)(0,0)
\put(50,0){\droite(0,1){40}}
\put(52,0){\droite(0,1){40}}
\put(50,40){\droite(-1,1){50}}
\put(50,42){\droite(-1,1){18}}
\put(50,40){\droite(1,1){50}}
\put(30,60){\droite(1,1){30}}
\put(30,62){\droite(1,1){18}}
\put(10,80){\droite(1,1){10}}
\put(50,80){\droite(-1,1){10}}
\put(50,82){\droite(-1,1){8}}
\put(90,80){\droite(-1,1){10}}

\put(55,35){$\dashv$}
\put(25,50){$\vdash$}
\put(5,70){$\vdash$}
\put(50,70){$\dashv$}
\put(90,70){$\dashv$}
\endpicture}\kern 0.4ex}

\def\TREEG{\kern-0.4ex
\hbox{\unitlength=1pt
\picture(100,90)(0,0)
\put(50,0){\droite(0,1){40}}
\put(52,0){\droite(0,1){40}}
\put(50,40){\droite(-1,1){50}}
\put(50,40){\droite(1,1){50}}
\put(50,42){\droite(1,1){18}}
\put(10,80){\droite(1,1){10}}
\put(50,80){\droite(1,1){10}}
\put(70,60){\droite(-1,1){30}}
\put(70,62){\droite(-1,1){30}}
\put(50,80){\droite(-1,1){10}}
\put(50,82){\droite(-1,1){8}}
\put(90,80){\droite(-1,1){10}}

\put(55,35){$\vdash$}
\put(75,50){$\dashv$}
\put(5,70){$\vdash$}
\put(45,70){$\dashv$}
\put(90,70){$\dashv$}
\endpicture}\kern 0.4ex}

\def\arbretroisun{\kern-0.4ex
\hbox{\unitlength=.50pt
\picture(80,100)(0,0)
\put(30,0){\droite(0,1){20}}
\put(30,20){\droite(-1,1){40}}
\put(30,20){\droite(1,1){40}}
\put(60,50){\droite(-1,1){10}}
\put(10,40){\droite(1,1){20}}
\put(70,50){... 1}
\put(70,35){... 2}
\put(70,20){... 3}
\endpicture}\kern 0.4ex}

\def\arbretroisdeux{\kern-0.4ex
\hbox{\unitlength=.50pt
\picture(80,100)(0,0)
\put(30,0){\droite(0,1){20}}
\put(30,20){\droite(-1,1){40}}
\put(30,20){\droite(1,1){40}}
\put(50,40){\droite(-1,1){20}}
\put(0,50){\droite(1,1){10}}
\put(70,50){... 1}
\put(70,35){... 2}
\put(70,20){... 3}
\endpicture}\kern 0.4ex}

\def\dessinun{\kern-0.4ex
\hbox{\unitlength=.25pt
\picture(60,40)(0,0)
\put(0,0){\droite(1,1){40}}
\put(20,0){\droite(1,1){40}}
\endpicture}\kern 0.4ex}

\def\leveltree{\kern-0.4ex
\hbox{\unitlength=1pt
\picture(100,90)(0,0)

\put(50,0){\droite(0,1){30}}
\put(50,30){\droite(-1,1){50}}
\put(50,30){\droite(1,1){50}}
\put(10,70){\droite(1,1){10}}
\put(70,50){\droite(-1,1){30}}

\put(200,0){\droite(0,1){30}}
\put(200,30){\droite(-1,1){50}}
\put(200,30){\droite(1,1){50}}
\put(180,50){\droite(1,1){30}}
\put(240,70){\droite(-1,1){10}}

\put(-30, 70){1 -}
\put(-30, 50){2 -}
\put(-30, 30){3 -}
\put(-30, 0){level}
\endpicture}\kern 0.4ex}

\def\arbreEFA{\kern-0.4ex
\hbox{\unitlength=.50pt
\picture(60,40)(0,0)
\put(30,0){\droite(0,1){45}}
\put(30,20){\droite(-1,1){25}}
\put(30,20){\droite(1,1){25}}
\put(15,35){\droite(1,1){10}}
\put(5,30){$\scriptscriptstyle \succ$}
\endpicture}\kern 0.4ex}

\def\arbreEFB{\kern-0.4ex
\hbox{\unitlength=.50pt
\picture(60,40)(0,0)
\put(30,0){\droite(0,1){45}}
\put(30,20){\droite(-1,1){25}}
\put(30,20){\droite(1,1){25}}
\put(45,35){\droite(-1,1){10}}
\put(45,30){$\scriptscriptstyle \prec$}
\endpicture}\kern 0.4ex}

\def\arbreEFC{\kern-0.4ex
\hbox{\unitlength=.50pt
\picture(60,40)(0,0)
\put(30,0){\droite(0,1){20}}
\put(30,20){\droite(-1,1){30}}
\put(30,20){\droite(1,1){30}}
\put(15,35){\droite(1,1){15}}
\put(5,30){$\scriptscriptstyle \prec$}
\put(20,15){$\scriptscriptstyle \prec$}
\endpicture}\kern 0.4ex}

\def\arbreEFD{\kern-0.4ex
\hbox{\unitlength=.50pt
\picture(60,40)(0,0)
\put(30,0){\droite(0,1){20}}
\put(30,20){\droite(-1,1){30}}
\put(30,20){\droite(1,1){30}}
\put(45,35){\droite(-1,1){15}}
\put(45,30){$\scriptscriptstyle \succ$}
\put(30,15){$\scriptscriptstyle \succ$}
\endpicture}\kern 0.4ex}

\def\arbreEFE{\kern-0.4ex
\hbox{\unitlength=.50pt
\picture(60,40)(0,0)
\put(30,0){\droite(0,1){20}}
\put(30,20){\droite(-1,1){30}}
\put(30,20){\droite(1,1){30}}
\put(15,35){\droite(1,1){15}}
\put(5,30){$\scriptscriptstyle \succ$}
\put(20,15){$\scriptscriptstyle \prec$}
\endpicture}\kern 0.4ex}

\def\arbreEFF{\kern-0.4ex
\hbox{\unitlength=.50pt
\picture(60,40)(0,0)
\put(30,0){\droite(0,1){20}}
\put(30,20){\droite(-1,1){30}}
\put(30,20){\droite(1,1){30}}
\put(45,35){\droite(-1,1){15}}
\put(45,30){$\scriptscriptstyle \prec$}
\put(30,15){$\scriptscriptstyle \succ$}
\endpicture}\kern 0.4ex}

\def\arbreEFX{\kern-0.4ex
\hbox{\unitlength=.50pt
\picture(60,40)(0,0)
\put(30,0){\droite(0,1){45}}
\put(30,20){\droite(-1,1){25}}
\put(30,20){\droite(1,1){25}}
\put(28,50){$.$}
\put(11,30){$\scriptscriptstyle \bullet$}
\put(39,30){$\scriptscriptstyle \bullet$}
\put(26,8){$\scriptscriptstyle \bullet$}
\endpicture}\kern 0.4ex}

\def\arbreEFY{\kern-0.4ex
\hbox{\unitlength=.50pt
\picture(60,40)(0,0)
\put(30,0){\droite(0,1){20}}
\put(30,20){\droite(-1,1){25}}
\put(30,20){\droite(1,1){25}}
\put(0,50){$.$}
\put(39,30){$\scriptscriptstyle \bullet$}
\put(26,8){$\scriptscriptstyle \bullet$}
\put(23,25){$\scriptscriptstyle \prec$}
\endpicture}\kern 0.4ex}

\def\arbreEFZ{\kern-0.4ex
\hbox{\unitlength=.50pt
\picture(60,40)(0,0)
\put(30,0){\droite(0,1){20}}
\put(30,20){\droite(-1,1){25}}
\put(30,20){\droite(1,1){25}}
\put(50,50){$.$ }
\put(11,30){$\scriptscriptstyle \bullet$}
\put(26,8){$\scriptscriptstyle \bullet$}
\put(26,25){$\scriptscriptstyle \succ$}
\endpicture}\kern 0.4ex}

\def\trivalent{\kern-0.4ex
\hbox{\unitlength=.50pt
\picture(60,40)(0,0)
\put(30,0){\droite(0,1){45}}
\put(30,20){\droite(-1,1){25}}
\put(30,20){\droite(1,1){25}}
\endpicture}\kern 0.4ex}

\def\arbretroisbis{\kern-0.4ex
\hbox{\unitlength=1pt
\picture(60,40)(0,0)
\put(30,0){\droite(0,1){20}}
\put(30,20){\droite(-1,1){30}}
\put(30,20){\droite(1,1){30}}
\put(50,40){\droite(-1,1){10}}
\put(10,40){\droite(1,1){10}}
\put(70,40){... 1}
\put(70,30){... 2}
\put(70,20){... 3}
\endpicture}\kern 0.4ex}

\def\arbretroisun{\kern-0.4ex
\hbox{\unitlength=.50pt
\picture(80,100)(0,0)
\put(30,0){\droite(0,1){20}}
\put(30,20){\droite(-1,1){40}}
\put(30,20){\droite(1,1){40}}
\put(60,50){\droite(-1,1){10}}
\put(10,40){\droite(1,1){20}}
\put(70,50){... 1}
\put(70,35){... 2}
\put(70,20){... 3}
\endpicture}\kern 0.4ex}

\def\arbretroisdeux{\kern-0.4ex
\hbox{\unitlength=.50pt
\picture(80,100)(0,0)
\put(30,0){\droite(0,1){20}}
\put(30,20){\droite(-1,1){40}}
\put(30,20){\droite(1,1){40}}
\put(50,40){\droite(-1,1){20}}
\put(0,50){\droite(1,1){10}}
\put(70,50){... 1}
\put(70,35){... 2}
\put(70,20){... 3}
\endpicture}\kern 0.4ex}

\def\SUBTREEun{\kern-0.4ex
\hbox{\unitlength=1pt
\picture(90,60)(0,0)
\put(40,0){\droite(0,1){20}}
\put(40,20){\droite(-1,1){40}}
\put(30,50){\droite(-1,1){10}}
\put(70,50){\droite(-1,1){10}}
\put(20,40){\droite(1,1){20}}
\put(40,20){\droite(1,1){40}}
\endpicture}\kern 0.4ex}

\def\SUBTREEdeux{\kern-0.4ex
\hbox{\unitlength=1pt
\picture(90,60)(0,0)
\put(45,0){\droite(0,1){15}}
\put(47,0){\droite(0,1){15}}
\put(45,15){\droite(-1,1){45}}
\put(47,15){\droite(-1,1){35}}
\put(40,50){\droite(-1,1){10}}
\put(42,50){\droite(-1,1){10}}
\put(75,45){\droite(-1,1){15}}
\put(77,45){\droite(-1,1){15}}
\put(10,50){\droite(1,1){10}}
\put(12,50){\droite(1,1){10}}
\put(25,35){\droite(1,1){25}}
\put(27,35){\droite(1,1){25}}
\put(70,50){\droite(1,1){10}}
\put(72,50){\droite(1,1){10}}
\put(45,15){\droite(1,1){45}}
\put(47,15){\droite(1,1){30}}
\endpicture}\kern 0.4ex}

\def\SUBTREEtrois{\kern-0.4ex
\hbox{\unitlength=1pt
\picture(90,60)(0,0)
\put(40,0){\droite(0,1){20}}
\put(40,20){\droite(-1,1){40}}
\put(65,45){\droite(-1,1){15}}
\put(15,45){\droite(1,1){15}}
\put(40,20){\droite(1,1){40}}
\endpicture}\kern 0.4ex}

\def\Mysmalldiagram{
\matrix{
&&{\bf Dend}&&&&{\bf Dias}&\cr
&&&&&&&\cr
&\nearrow&&\searrow&&\nearrow&&\searrow\cr
&&&&&&&\cr
{\bf Zinb}&&&&{\bf As}&&&&{\bf Leib}\cr
&&&&&&&\cr
&\searrow&&\nearrow&&\searrow&&\nearrow\cr
&&&&&&&\cr
&&{\bf Com}&&&&{\bf Lie}&\cr
}}

\centerline {\bf DIALGEBRAS}
\BB

\centerline {\bf Jean-Louis LODAY}

\vglue 1cm
 There is a notion of ``non-commutative Lie algebra" called {\it Leibniz algebra}, which
 is characterized by the following property. The bracketing $[-,z]$ is a
derivation for the bracket operation, that is, it satisfies the Leibniz identity: 
$$
[[x,y],z] = [[x,z],y] + [x, [y,z]].
$$
 cf. [L1]. When it happens that the bracket is  skew-symmetric, we get a Lie algebra 
since the Leibniz identity becomes equivalent to the Jacobi identity.

Any associative algebra gives rise to a Lie algebra by $[x,y] = xy - yx$. The purpose of
this article is to introduce and study a new notion of algebra which gives, by a
similar procedure, a Leibniz algebra. The idea is to start with two distinct operations
for the product  $xy$ and the product $yx$, so that the bracket is not necessarily
skew-symmetric any more. Explicitly, we define an {\it associative dialgebra} (or simply
dialgebra for short) as a vector space
$D$ equipped with two associative operations $\g $ and $ \d $, called respectively left
and right product, satisfying 3 more axioms:
$$\left\{ \eqalign{
x \dashv (y\dashv z) & = x \dashv (y \vdash z) , \cr
  (x\vdash y)\dashv z & = x\vdash (y\dashv
z),\hfill\cr
 (x \dashv y)\vdash z & =(x \vdash y)\vdash z.\cr } \right.
$$
It is immediate to check that $[x,y]:= x\g y - y\d x$ defines a Leibniz bracket. Hence any dialgebra
gives rise to a Leibniz algebra.

A typical example of dialgebra is constructed as follows. Let $(A,d)$ be a differential associative
algebra, and put 
$$
x\g y:= x\ dy \quad \hbox{and }\quad  x\d y:= dx\ y .
$$
One easily checks that $(A, \g , \d )$ is a dialgebra. For instance there is a natural dialgebra
structure on the de Rham complex of a manifold.

Observe that, since the relations defining a
dialgebra do not involve sums, there is a well-defined notion of {\it dimonoid}. 

In this
article we construct and study a (co)homology theory for dialgebras. Since an associative
algebra is a particular case of dialgebra, we get a new (co)homology theory for associative
algebras as well. The surprizing fact, in the construction of the chain complex, is the
appearance of the combinatorics of {\it planar binary trees}. The principal result about this
homology theory $HY$ is its vanishing on free dialgebras. In order to state some of the 
properties of the theory $HY$, we introduce  another type of algebras with two operations:
the {\it dendriform  algebras} (sometimes called dendriform dialgebras). This notion dichotomizes the notion of
associative algebra in the following sense: there are two operations $\prec$ and $\succ$, such that the product
$*$ made of the sum of them
$$
x*y:=x\prec y + x \succ y,
$$ 
is associative. The axioms relating these two products are
$$
\eqalign{
\hbox{(i)}\quad &(a\prec b)\prec c=a\prec(b\prec c)+a\prec(b\succ c),\cr
\hbox{(ii)}\quad& (a\succ b)\prec c=a\succ (b\prec c),\cr
\hbox{(iii)}\quad& (a\prec b)\succ c+(a\succ b)\succ c=a\succ (b\succ c).\cr}
$$
The free dendriform  algebra can be constructed by means of the planar binary trees, whence the
terminology.

The results intertwining associative dialgebras and dendriform  algebras are best expressed in
the framework of {\it algebraic operads}. The notion of dialgebra defines an algebraic operad
$Dias$, which is binary and quadratic. By the theory of Ginzburg and Kapranov (cf. [GK]), there is a
well-defined ``dual operad"
$Dias^!$. We show that this is precisely the operad $Dend$ of the  ``dendriform  algebras",
in other words a dual associative dialgebra is nothing but a dendriform  algebra. The vanishing of $HY$ of a free
dialgebra implies that those two operads are of a special kind: they are ``Koszul operads". As a
consequence the cohomology of a dialgebra is a graded dendriform  algebra and, a fortiori, a graded
associative algebra. The explicit description of the free dendriform  algebra in terms of trees permits
us to describe the notion of {\it  strong homotopy associative dialgebra}.

The categories of algebras over these operads assemble into a commutative diagram of functors which
reflects the Koszul duality.
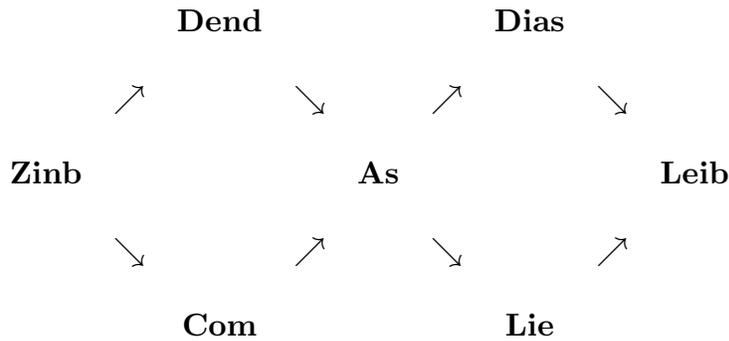
\B

\qquad \qquad $\Mysmalldiagram$

\BB
In this diagram {\bf Zinb} denotes the categories of Zinbiel algebras, which are Koszul dual to the Leibniz algebras.
\M

This paper is part of a long-standing project whose ultimate aim is to study periodicity phenomenons
in algebraic $K$-theory. This project is described in [L4]. The next step would consist in computing
the dialgebra homology of the augmentation ideal of $K[GL(A)]$, for an associative algebra
$A$.
 \B

 Here is the content of this article. In the
first section we introduce the notion of associative dimonoid, or dimonoid for short,  and develop the calculus in a
dimonoid. In particular we describe the free dimonoid on a given set. In the second section we introduce the notion
of dialgebra and give several examples. We explicitly describe the free dialgebra over a vector space. In
the third section we construct the chain complex of a dialgebra $D$, which gives rise to
homology and cohomology groups denoted $HY(D)$. The main tool is made of the planar binary
trees and operations on them. We prove that $HY$ of a free dialgebra vanishes (hence the
operad  associated to dialgebras is a Koszul operad). We also introduce a variation of the
chain complex by replacing the trees by increasing trees, or, equivalently, by permutations.
This variation appears naturally in the computation of the Leibniz homology of dialgebras of matrices
(cf. [F1]). (Co)homology of dialgebras with non-trivial coefficients is treated by Alessandra
Frabetti in [F4].

Section 4 is devoted to the relationship between Leibniz algebras and dialgebras. The
functor which assigns to any dialgebra  $(D,  \vdash , \dashv)$ the Leibniz
algebra $(D, [x,y]:= x\dashv y - y\vdash x)$ has a left adjoint which is  the {\it universal
enveloping dialgebra} of a Leibniz algebra. Then we compare the diverse types of free algebras and we propose a
definition for a {\it Poisson dialgebra}.  The Hopf-type properties of the universal enveloping
dialgebra are studied by Fran\c cois Goichot in [Go].

In the fifth section we introduce the notion of {\it dendriform  algebra}, which is closely
connected to the notion of associative dialgebra. For instance the tensor product of a dialgebra and
of a dendriform  algebra is naturally equipped with a structure of Lie algebra. The main result of
this section is to make explicit the free dendriform  algebra. It turns out that it is best
expressed in terms of planar binary trees. The dendriform  algebra structure on the vector
space generated by the planar binary trees is the core of this section. It uses the
grafting operation and the nesting operation on trees and it induces a graded associative algebra
structure on the same vector space. In a sense associative algebras are closely connected with the
integers (including addition and multiplication). Similarly dendriform  algebras are closely connected
with planar binary trees and a calculus on them. This arithmetic aspect of the theory will be treated
elsewhere.

In section 6 we
construct (co)homology groups for dendriform  algebras. They vanish on free dendriform  algebras.

In section 7 we relate
dendriform  algebras with Zinbiel algebras (i.e. dual-Leibniz algebras) and associative algebras. It is
based on the relationship between binary trees and permutations as described in Appendix A.

The aim of the eighth section is to interpret the preceding results in the context of
algebraic operads. The basics on algebraic operads and Koszul duality are recalled in
Appendix B. We show that the operads associated to dialgebras and to dendriform  algebras
are dual in the operad sense. Then we show that the (co)homology groups $HY$ for
dialgebras (resp. $H^{Dend}$ for dendriform  algebras) constructed in section 3 (resp. 4) are the
ones predicted by the operad theory. Hence, by the vanishing of $HY$ of a free dialgebra,
both operads  $Dias$ and   $Dend$ are Koszul. It implies, among several
consequences, the vanishing of the homology of a free dendriform  algebra. Some of the theorems
in sections 2 to 6 can be proved either directly or by appealing to the operad
theory. In general we write down the most elementary one.

The last section describes the notion of {\it strong homotopy associative dialgebras}. For any Koszul operad the
notion of algebra up to homotopy is theoretically  well-defined from the bar construction over the
dual operad. Since, in our case, we know explicitly the structure of a free dendriform  algebra, we can
make the notion of dialgebra up to homotopy completely explicit.
\M

Part of the results of this article has been announced in a ``Note aux
Comptes Rendus" [L2]. I thank  Ale Frabetti, Benoit Fresse, Victor Gnedbaye, Fran\c cois Goichot, Phil Hanlon, 
Muriel Livernet, Teimuraz Pirashvili and Maria Ronco
for fruitful conversations on this subject.
\vfill
\eject

\N {\bf Contents}.
\B 

\noindent 1. Dimonoids

\noindent 2. Associative dialgebras

\noindent 3. (Co)homology of associative dialgebras

\noindent 4. Leibniz algebras, associative dialgebras and homology

\noindent 5. Dendriform  algebras

\noindent 6. (Co)homology of dendriform  algebras

\noindent 7. Zinbiel algebras, dendriform  algebras and homology

\noindent 8. Koszul duality for the dialgebra operad

\noindent 9. Strong homotopy associative dialgebras

\noindent Appendix A. Planar binary trees and permutations
 
\noindent Appendix B. Algebraic operads

\vfill\eject



\noindent {\bf 1. DIMONOIDS}
\M

\noindent {\bf 1.1. Definition.} 
An {\it associative dimonoid}, or dimonoid for short,  is a set $X$ equipped with two maps
called respectively {\it left product} and {\it right product}:
$$\leqalignno{
& \dashv\ : X\times X\r X,&{(\rm left)}\cr
& \vdash\ : X\times X\r X,&{(\rm right)}\cr}
$$
satisfying the following axioms

$$
\cases {
x\dashv(y\dashv z)\buildrel{1}\over =(x\dashv y)\dashv z\buildrel{2}\over =x\dashv
(y\vdash z) ,\cr 
(x\vdash y)\dashv z \buildrel{3}\over= x\vdash (y\dashv z) ,\cr
(x\dashv y)\vdash z\buildrel{4}\over=
x\vdash
(y\vdash z)\buildrel{5}\over=(x\vdash y)\vdash z ,\cr}
$$
for all $x,y$ and $z\in X$.

In the notation $x\dashv y$, $y\vdash x$, the element $x$ is said to be on the {\it
pointer} side and the element $y$ is said to be on the {\it bar} side.

The numbers $1$ to $5$ of the relations are for future reference. 

Observe that relations 1 and 5 are the ``associativity'' of the products  $\dashv$ and
$\vdash$  respectively.
Relation 3 will be referred to as ``inside associativity'', since the products point inside.
Relations 2 and 4 can be replaced by the relations 12 and 45:
$$
x\dashv(y\dashv z)\buildrel{12}\over=x\dashv(y\vdash z)\quad \hbox{and}\quad
(x\dashv y)\vdash z\buildrel{45}\over=(x\vdash y)\vdash z, 
$$
which can be summarized as ``on the bar side, does not matter which product''.
All these relations are referred to as ``diassociativity''.

A {\it morphism} of dimonoids is a map $f:X\r Y$ between two dimonoids $X$ and $Y$ such
that $f(x\dashv x')=f(x)\dashv f(x')$ and $f(x\vdash x')=f(x)\vdash f(x')$ for any
$x,x'\in X$.

Observe that one can define a di-object in any monoidal category. One does not need the
monoidal category to be symmetric since in each relation the variables stay in the same
order.
\B

\noindent {\bf 1.2. Bar-unit.} 
An element $e\in X$ is said to be a {\it bar-unit} of the dimonoid $X$ if
$$
x\dashv e=x=e\vdash x, \quad \hbox{for\ any}\ x\in X.
$$
So it is only assumed that $e$ acts trivially from the bar side. There is no reason for
a bar-unit to be unique. The set of bar-units is called the  {\it halo}.

A morphism of dimonoids is said to be {\it unital} if the image of a bar-unit is a bar-unit.
\M

\noindent {\bf1.3. Examples.} 

a) Let $M$ be a monoid (without unit), that is a set $M$ with an associative product
$(m,m')\mapsto mm'$. Putting $m\dashv m'=mm'=m\vdash m'$ gives a dimonoid structure on
$M$. Indeed each relation 1 to 5 is the associativity property. A unit of the monoid is a bar-unit
of the associated dimonoid.

Conversely, if in a dimonoid $D$ there is a unit, that is an element  $1\in D$ which
satisfies either
$1\dashv x=x $ or $x=x\vdash 1$ for all $x\in D$, then, by axiom 3 or 5, one has $\dashv\
=\ \vdash$ and $D$ is simply the dimonoid associated to a unital monoid.
\medskip

b) Let $X$ be a set and define 
$$
x\dashv y=x=y\vdash x,\quad\hbox{for any}\quad x,y\in X.
$$
Then, obviously, $X$ is a (not so interesting) dimonoid and it coincides with its halo.
\M

c) Let $M$ be a monoid. Put $D=M\times M$ and define the products by

$$
 \biggl\{ \matrix {&(m,n)\dashv (m',n'):=(m,nm'n')\cr
&(m,n)\vdash (m',n'):=(mnm',n').\cr}
$$
With these definitions $D=(D,\dashv,\vdash)$ is a dimonoid. Let us check relation 3 for
instance:
$$
\eqalign{
((m,n)\vdash (m',n'))\dashv (m'',n'')=(mnm',n')\dashv(m'',n'')=(mnm',n'm''n'')\cr
(m,n)\vdash ((m',n')\dashv (m'',n''))=(m,n)\vdash (m',n'm''n'')=(mnm',n'm''n'').\cr}
$$
Let $1\in M$ be a unit for $M$. Then $e=(1,1)$ is a bar-unit for $D$, but one has 
$e\dashv x\ne x$ and $x\vdash e\ne x$ in $D$ in general. For any invertible element $m$ the
element $(m,m^{-1})\in D$ is a bar-unit.
\medskip

d) Let $G$ be a group and $X$ a $G$-set. The following formulas define a dimonoid
structure on $X\times G$ (cf. 7.5): 
$$
 \biggl\{ \matrix {
(x,g)\g (y,h):=(x,gh),\cr 
(x,g)\d (y,h):=(g\cdot y,gh).\cr
}
$$ 

\noindent {\bf 1.4. Opposite dimonoid.} 
Let $D$ be a dimonoid. Define new operations $\dashv '$ and $\vdash '$ on $D$ by
$$
\eqalign{
x\dashv 'y:=y\vdash x,\cr
x\vdash 'y:=y\dashv x.\cr}
$$
It is immediate to check that $(D,\dashv ',\vdash ')$ is a new dimonoid which we call the
{\it opposite} dimonoid that we denote by $D^{op}$.

Observe that if we put
$$
 \Bigl\{ \matrix {
x\dashv''y:=y\dashv x,\cr
x\vdash''y:=y\vdash x,\cr}
$$
then $(D,\vdash'',\dashv'')$ is {\it not} a dimonoid.
\B

\noindent {\bf 1.5. Monomials in a dimonoid.} 
Let $x_1,\dots ,x_n$ be elements in the dimonoid $D$. A {\it monomial} in $D$ is a
parenthesizing together with product signs, for instance
$$
((x_1\dashv x_2)\vdash (x_3\dashv x_4))\dashv (x_5\vdash x_6),
$$
giving rise to an element in $D$. Such a monomial is completely determined by a 
binary tree, where each vertex is labelled by $\dashv$ or $\vdash$:
\M

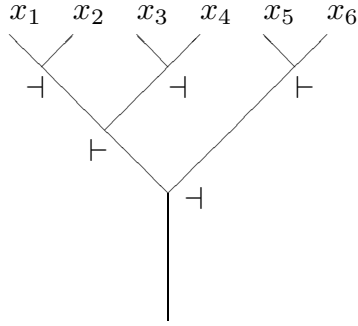
$$\TRE$$

\noindent {\bf 1.6. The middle of a monomial.} 
Given a monomial as above we define the {\it middle} of the monomial as being the entry
$x_i$ determined by the following algorithm. Starting at the root of the tree one goes up
by choosing the route indicated by the pointer. The middle of the monomial is the abutment
of the path. In this example $x_3$ is the middle.
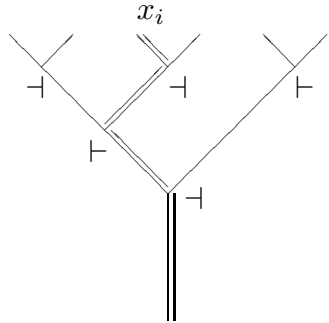
\M

$$\TREEE$$

\noindent {\bf  1.7. Theorem  (Dimonoid calculus).} 
 {\it Let $x_i$, $i\in {\bf Z}$, be elements in a dimonoid $D$.

a) Any parenthesizing of
$$
x_{-n}\vdash x_{-n+1}\vdash \dots \vdash x_{-1}\vdash x_0\dashv x_1\dashv \dots \dashv x_{m-1}
\dashv x_m
$$
gives the same element in $D$, which we denote by
$$
x_{-n}\dots x_{-1}\check x_0x_1\dots x_m.
$$
\medskip

b) Let $m=x_1\dots x_k$ be a monomial in $D$. Let $x_i$ be its middle entry. Then
$m=x_1\dots \check x_i\dots x_k$.
\medskip

c) One has the following formulas in $D$:
$$
\eqalign{
&(x_1\dots \check x_i\dots x_k)\dashv (x_{k+1}\dots \check x_j\dots x_{\ell})=x_1\dots \check
x_i\dots x_kx_{k+1}\dots x_j\dots x_{\ell}\cr
&(x_1\dots \check x_i\dots x_k)\vdash (x_{k+1}\dots \check x_j\dots x_{\ell})=x_1\dots 
x_i\dots x_kx_{k+1}\dots \check x_j\dots x_{\ell}.\cr}
$$}
For instance, in the above example, one has
$$
((x_1\dashv x_2)\vdash (x_3\dashv x_4))\dashv (x_5
\vdash x_6)=x_1x_2\check x_3x_4 x_5x_6.
$$
{\it Proof}. By axiom 1 (associativity of $\dashv$) any parenthesizing of
$x_1\dashv\dots\dashv x_m$ gives the same element. So, in such a monomial we can ignore
the parentheses (and analogously for $\vdash$ thanks to axiom 5).

Consider a generic monomial with first entry $x_{-n}$, last entry $x_m$ and middle entry
$x_0$ (where $-n\le 0\le m$). By axioms $1-3-5$ it is clear that the element
$$
(x_{-n}\vdash \dots\vdash x_{-1})\vdash x_0\dashv(x_1\dashv\dots\dashv x_m)\leqno(*)
$$
is well-defined. We denote it by $x_{-n}\dots \check x_0\dots x_m$.

Consider the labelled tree of our generic monomial. Let $v$ be a vertex which is on the
route from the root to the middle entry $x_0$. Thanks to axioms 12 and 45 all the vertices
on the bar side of $v$ can be labelled with the same label as $v$. In our example
$$\TREEF$$

Then by axiom 3 we can modify the tree so that all labels $\vdash$ come first:

$$\TREEG$$

This new tree corresponds to a monomial of the form $(*)$ and therefore we have proved that
our starting monomial has value $x_{-n}\dots \check x_0\dots x_m$. So parts a) and b) are proved.

By a) and b) it follows that in order to compute 
$$
(x_1\dots \check x_i\dots x_k)\dashv
(x_{k+1}\dots \check x_j\dots x_{\ell})\quad {\rm and} \quad (x_1\dots \check x_i\dots x_k)\vdash
 (x_{k+1}\dots \check
x_j\dots x_{\ell})
$$
 it suffices to determine which entry is the middle of these monomials. By
the algorithm described in 1.6, the middle entry is $x_i$ in the first case and $x_j$ in
the second case. 
\hfill$\square$
\M

\N {\bf 1.8. Corollary.} {\it The free dimonoid on the set $X$ is the disjoint union
$${\cal D}(X) =
\bigsqcup_{n\geq 1}(\underbrace {X^n\cup \cdots \cup X^n}_{n \rm \; copies}).
$$
Denoting by
$x_1\dots
\check x_i\dots x_n$ an element in the $i$-th summand, the products are given by 
$$\eqalign{
(x_1\dots \check x_i\dots x_k)\dashv
(x_{k+1}\dots \check x_j\dots x_{\ell})&= x_1\dots \check x_i\ \dots x_{\ell}\ , \cr
\quad (x_1\dots \check x_i\dots x_k)\vdash
 (x_{k+1}\dots \check
x_j\dots x_{\ell})&= x_1\dots\  \check x_j\dots x_{\ell}\ .\cr}
$$}
\hfill $\square$
\vfill\eject



\noindent {\bf  2. ASSOCIATIVE DIALGEBRAS}
\medskip

In the sequel $K$ denotes a field referred to as the {\it
ground field}. Later on it will be supposed to be of characteristic zero. The tensor product over $K$ is
denoted by
$\otimes_K$ or, more often, by $\otimes$.

After introducing the notion of dialgebra, we give some examples, including free
dialgebras, which we describe explicitly, and define modules and representations over a
dialgebra. 
\M

\noindent {\bf 2.1. Definition.} 
An {\it associative dialgebra}, or {\it dialgebra} for short, over $K$ is a $ 
K $-module $D$ equipped with two $K$-linear maps
$$
\eqalign{
&\dashv\ : D\otimes D\longrightarrow D,\cr
&\vdash\ : D\otimes D\longrightarrow D,\cr}
$$
satisfying the {\it di-associativity} axioms
$$\left\{ {
\matrix { 
(1)\qquad (x\dashv y)\dashv z = x\dashv (y\vdash z),\cr
(2)\qquad (x\dashv y)\dashv z = x\dashv(y\dashv z),  \cr
(3)\qquad (x\vdash y)\dashv z = x\vdash (y\dashv z), \cr 
(4)\qquad (x\dashv y)\vdash z= x\vdash (y\vdash z) , \cr
(5)\qquad (x\vdash y)\vdash z= x\vdash (y\vdash z) . \cr
}}\right.
$$
The maps $\dashv$ and $\vdash$ are called respectively the {\it left product} and the {\it  right
product.}

Here is an equivalent formulation of these axioms:
the products $\g$ and $\d$ are associative and satisfy:
$$\left\{ \eqalign{
(12)\qquad x \dashv (y\dashv z) & = x \dashv (y \vdash z) , \cr
(3)\qquad (x\vdash y)\dashv z & = x\vdash (y\dashv
z),\hfill\cr
 (45)\qquad (x \dashv y)\vdash z & =(x \vdash y)\vdash z.\cr } \right.
$$
Observe that the analogue of formula (3), but with the product symbols pointing outward, is not valid
in general: $ (x\dashv y)\vdash z \ne x\dashv (y\vdash z)$.

A {\it morphism} of dialgebras from $D$ to $D'$ is a $K$-linear map  $f:D\r D'$ such
that $$
f(x\dashv y)=f(x)\dashv f(y)
 \quad\hbox{and}\quad f(x\vdash y)=f(x)\vdash f(y)\quad \hbox{ for all}\quad x,y\in D.$$
We denote by {\bf Dias} the category of dialgebras.
\medskip

A {\it bar-unit} in $D$ is an element $e\in D$ such that
$$
x\dashv e=x=e\vdash x\quad \hbox{for\ all}\quad x\in D.
$$
A bar-unit need not be unique. The subset of bar-units of $D$ is called its {\it halo}.

A {\it unital dialgebra} is a dialgebra with a specified bar-unit $e$. This choice gives rise to a
preferred $K$-linear map $K\hookrightarrow D, \lambda \mapsto \lambda e$.

A morphism of dialgebras is said to be {\it unital} if the image of any bar-unit is a
bar-unit.

Observe that if a dialgebra has a unit $\epsilon$, that is an element which satisfies
$\epsilon \g x = x$ for any $x$, then $\g = \d$ by axiom 12, and $D$ is an associative algebra
with unit.

An {\it ideal} $I$ in a dialgebra $D$ is a submodule of $D$ such that $x\dashv y$ and
$x\vdash y$ are in $I$ whenever one of the variables is in $I$. Clearly the quotient $D/I$
is a dialgebra. Conversely, the kernel of a dialgebra morphism is an ideal.
\M

\noindent {\bf  2.2. Examples.}
\S 

{\it a) Associative algebra.} If $A$ is an associative algebra over $K$, then the formulas $a\dashv
b=ab=a\vdash b$ define a structure of dialgebra on $A$. If 1 is a unit of the associative
algebra, then $e=1$ is a unit of the dialgebra and the halo is just $\{ 1 \} $.
\S

{\it b) Differential associative algebra.} Let $(A, d)$ be a differential associative algebra. So, by
hypothesis,
$d(ab) = da\, b+a\, db$ (here we work in the non-graded setting) and $d^2=0$. Define left and right products
on $A$ by the formulas
$$
x\g y:= x\ dy \quad {\rm and}\quad x\d y:= dx\ y.
$$
It is immediate to check that $A$ equipped with these two products is a dialgebra. A similar 
construction holds in the graded (or more accurately super) algebra framework.
\S

{\it c) Dimonoid algebra.} Let $X$ be a dimonoid, and denote by $K[X]$ the free $K$-module on $X$.
Then obviously $K[X]$ is a dialgebra.
\S

{\it d) Bimodule map.} Let $A$ be an associative algebra and let $M$ be an $A$-bimodule. Let $f:M\r A$ be an
$A$-bimodule map. Then one can put a dialgebra structure on $M$ as follows:
$$
\eqalign{
&m\dashv m':= mf(m'),\cr
&m\vdash m':= f(m)m',\cr}
$$
The verification is left to the reader. One can systematize this procedure by considering the tensor
category of linear maps as follows (cf. [LP2], [Ku] for details). The category of linear maps over
$K$ is made of the $K$-linear maps $f:V\to W$ as objects. It can be equipped with a tensor product
by 
$$(V{\buildrel f \over \to }W)\t (V'{\buildrel f' \over \to }W') =V\t W' \ \oplus \ W\t V'\  
{\buildrel f\t 1 + 1\t f' \over \longrightarrow }\  W\t W'.
$$
An associative algebra in this tensor category defines a dialgebra structure on the source object.

The particular case of the projection
$M\oplus A
\to A$ shows that there is a dialgebra structure on $M\oplus A$ (cf. P. Higgins [Hi]).
\S

{\it e) Tensor product, matrices. } If $D$ and $D'$ are two dialgebras, then the tensor product $D \otimes
D'$ is also a dialgebra by $(a\otimes a') \star (b\otimes b') = (a \star b)\otimes (a' \star b')$ for
$\star = \dashv$ and $\vdash$. For instance the module of $n\times n-$matrices ${\cal
M}_n(D)={\cal M}_n(K)\otimes D$ is a dialgebra. The left and right products are
given by
$$
(\alpha \dashv \beta    )_{ij}=\sum_k\alpha _{ik}\dashv \beta    _{kj}\quad \hbox{and}\quad (\alpha \vdash \beta    )_{ij}=
\sum_k\alpha _{ik}\vdash \beta    _{kj}.$$
\S

{\it f)  Opposite dialgebra}. As for dimonoids, the opposite dialgebra of $D$ is the
dialgebra $D^{op}$ with the same underlying $K$-module and with products given by
$$
x\dashv' y=y\vdash x,\ x\vdash'y=y\dashv x.
$$

 {\it g)} Let $A$ be an associative algebra over $K$. Put
$D=A\otimes A$ and define $$
\eqalign{
&a\otimes b\dashv a'\otimes b':=a\otimes ba'b',\cr
&a\otimes b\vdash a'\otimes b':=aba'\otimes b'.\cr}
$$
Extending these formulas by linearity on $A\otimes A$ gives well-defined product
maps $\dashv$ and $\vdash$ on $D$ which satisfy the diassociativity axioms. If $1\in A$
is a unit of the associative algebra, then $1\otimes 1$ is a bar-unit for the dialgebra. More
generally, for any invertible element $x$ in $A$, the element $x\otimes x\m $ is a bar-unit. If
$I$ is a left ideal and $J$ is a right ideal, then the same formulas define a diassociative
algebra structure on $I\otimes_KJ.$ 
\S

{\it h)} Let $A$ be an associative algebra and $n$ be a positive integer. On the module of
$n$-vectors $D=A^n$ one puts:
$$
\eqalign{
(x\dashv y)_i&=x_i\Bigl(\sum^n_{j=1}y_j\Bigr)\quad \hbox{for}\quad
1\le   i\le   n\quad\hbox{and}\cr 
(x\vdash y)_i&=\Bigl(\sum^n_{j=1}x_j\Bigr)y_i\quad \hbox{for}\quad
1\le   i\le   n.\cr} $$
One easily checks that $D$ is a dialgebra. For $n=1$, this is example (a).
In fact this construction can be extended to any dialgebra $A$.
\B

\noindent {\bf  2.3. Module, bimodule, extension.}
A {\it left module over a dialgebra} $D$ is a $K$-module $M$ equipped with two linear maps
$$
\leqalignno{
&\dashv \ :D\otimes M\r M,&\hbox{(right structure)}\cr
 &\vdash \ : D\otimes M\r M,&\hbox{(left structure)}\cr}
$$
satisfying the axioms (1)-(5) whenever they make sense. There is, of course, a similar definition for
right modules.

A {\it bimodule} over a dialgebra $D$, also called a {\it representation}, is a $K$-module $M$
equipped with four linear maps
$$
\leqalignno{
&\dashv,\vdash \ : M\otimes D\r M,&\hbox{(right structures)}\cr
 &\dashv,\vdash \ : D\otimes M\r M,&\hbox{(left structures)}\cr}
$$
satisfying the axioms (1) to (5), whenever one of the entries $x,y$ or $z$ is in $M$ and
the two others are in $D$.

Obviously a bimodule over $D$ is, a fortiori, a left module and also a right module over $D$ ; and
$D$ is a bimodule over itself.
\S

 Let 
$$
0\to M\to \overline D \to  D \to 0
$$
be an abelian extension of dialgebras, that is an exact sequence of dialgebras such that any product of two
elements in $M$ is trivial. Then, it is immediate to check that $M$ is a representation of
$D$ in the above sense. \M

\noindent {\bf  2.4. Free associative dialgebra.} 
Let $V$ be a $K$-module. By definition the {\it free dialgebra} on $V$ is the
dialgebra $Dias (V)$ equipped with a $K$-linear map $i:V\r Dias(V)$ such that for any
$K$-module map $f:V\r D$, where $D$ is a dialgebra over $K$, there is a unique
factorization
$$
f:V\buildrel{i}\over\longrightarrow Dias (V)\buildrel{\phi}\over\longrightarrow
 D,
$$
where $\phi$ is a dialgebra morphism.

Equivalently the functor $Dias:(K{\rm -}\Mod )\r {\bf Dias}$ is left adjoint to the forgetful functor. The
following proposition proves the existence of the free dialgebra $ Dias (V)$ and gives an explicit
description of it in terms of the tensor module
$$
T(V):= K\oplus V\oplus V^{\otimes 2}\oplus\cdots \oplus V^{\otimes n}\oplus\cdots .
$$
\M

 \noindent {\bf 2.5. Theorem.} 
{\it  The free dialgebra on $V$ is the $K$-module 
$$Dias (V)=T(V)\otimes V\otimes
T(V)$$
 equipped with the two products induced by:
$$
\displaylines{
(v_{-n}\cdots v_{-1}\otimes v_0\otimes v_1\cdots v_m)\dashv (w_{-p}\cdots w_{-1}\otimes w_0\otimes 
w_1\cdots w_q)\hfill\cr
\hfill =v_{-n}\cdots v_{-1}\otimes v_0\otimes 
v_1\cdots v_mw_{-p}\cdots w_q,\cr
(v_{-n}\cdots v_{-1}\otimes v_0\otimes v_1\cdots v_m)\vdash (w_{-p}\cdots w_{-1}\otimes w_0\otimes 
w_1\cdots w_q)\hfill\cr
\hfill =v_{-n}\cdots v_mw_{-p}\cdots w_{-1}\otimes 
w_0\otimes w_1\cdots w_q,\cr}
$$
where $v_i$, $w_j\in V$.}
\medskip

With our notation (cf. 1.7) any additive generator of $Dias (V)$ can be written 
$$
v_{-n}\cdots v_{-1}\otimes v_0\otimes v_1\cdots v_m=v_{-n}\cdots v_{-1}\check v_0v_1\cdots v_m.
$$
{\it Proof}. It is immediate to check that $Dias (V)=(T(V)\otimes V\otimes
T(V),\dashv,\vdash$) is a dialgebra (cf. 1.7). The map $i:V\r Dias (V)$ is the
composite $V\simeq 1\cdot K\otimes V\otimes 1\cdot K\hookrightarrow T(V)
\otimes V\otimes T(V)$.  Starting with $f:V\r D$ the map
$\phi: Dias (V) \r D$ is given by
$$
\phi(v_{-n}\cdots v_{-1}\check v_0v_1\cdots v_m)=f(v_{-n})\cdots f(v_{-1})f(v_0)^{\check{}}f(v_1)\cdots 
f(v_m).
$$
It is obviously a dialgebra morphism. Moreover, by theorem 1.7, it is uniquely determined since it
should coincide with $f$ on $V\cong 1\cdot K\otimes V\otimes 1\cdot K$ and it should be a
morphism of dialgebras. Hence the inclusion $V\to Dias(V)$ is universal. \hfill $\square$
\medskip
\N {\it Remark.} A free dialgebra is a particular case of example 2.2.d, with $M=
T(V)\otimes V\otimes T(V)$ , $A=T(V)$ (the associative tensor algebra) and $f:M\to A$ being the
concatenation.
\medskip

Let $V$ be finite dimensional over $K$ generated by $x_1,\cdots ,x_n$. 
Let us describe the degree $n$
part of $T(V)\otimes V\otimes T(V)$ which is generated by all the monomials containing
$x_i$ once and only once, $1\le   i\le  n$. We denote it by $Dias (n)$.
These
monomials are the elements 
$$
(\sigma  ,i)(x_1,\cdots ,x_n):=x_{\sigma  ^{-1}(1)}\cdots \check x_{\sigma  ^{-1}(i)}\cdots
x_{\sigma  ^{-1}(n)},
\quad  \sigma  \in S_n,\ 1\le   i\le   n, 
$$
where $S_n$ is the symmetric group.
Therefore, as a left $S_n$-module, the multilinear part of this space is isomorphic to $n$ copies
of the regular representation of $S_n$:
$$
 Dias(n)\cong nK[S_n].
$$
The element $\sigma  $ in the $i$-th copy corresponds to the operation $(\sigma  ,i)$ described above (cf.
Corollary 1.8).
\medskip

{\it Examples}: 
$$
\eqalign{
&\bullet\ n=1,\  \hbox {one generator:}\ \check x_1.\cr
&\bullet\ n=2,\ \hbox{four generators:}\ \check x_1x_2,\check x_2x_1,x_1\check
x_2,x_2\check x_1.\cr
&\bullet\ n=3,\ \hbox{eighteen generators:}\ \check x_ix_jx_k,x_i\check
x_jx_k,x_ix_j\check x_k\cr} $$
for all permutations $i,j,k$ of 1,2,3.
\M

\noindent {\bf 2.6. Associative algebra associated to a dialgebra.} For any dialgebra $D$ let
$D_{As}$ be the quotient of
$D$ by the ideal generated by the elements $x\g y - x\d y$, for all  $x,y\in D$. It is clear that $\g
= \d $ in
$D_{As}$, hence $D_{As}$ is an associative algebra (non-unital in general). The quotient map
$\mu: D\proj D_{As}$ is  universal among the  maps from $D$ to associative algebras. In other words the
associativization functor $(-)_{As}: {\bf Dias} \to {\bf As}$ is left adjoint to $inc:  {\bf As} \to 
{\bf Dias}$.

Axioms 12 and 45 imply that the element $x\d y\g z$ in $D$ depends only on the values of $x$ and $z$ in
$D_{As}$. Hence $D$ is a $D_{As}$-bimodule and the projection map $\mu$ is a $D_{As}$-bimodule map. On the
other hand the dialgebra structure of $D$ is completely determined by $\mu$ and the $D_{As}$-bimodule
structure on the space $D$ since 
$$
x\g y = x\, \mu(y)\quad {\rm and }\quad x\d y = \mu (x)\, y \, , 
$$
cf. example 2.2.d. It is useful to write the element  $x\d y\g z$ as  $x \check y z$. Under this notation
the dialgebra calculus rules are
$$\displaylines{
x \check y z \g s \check t u = x \check y z s t u , \cr 
x \check y z \d s \check t u = x  y z   s \check t u . \cr }
$$

\vfill\eject



\noindent {\bf  3. (CO)HOMOLOGY OF ASSOCIATIVE DIALGEBRAS}
\medskip
In this section we introduce a chain complex which permits us to define homology groups
$HY_*(D)$ and cohomology groups $HY^*(D)$ of a dialgebra $D$. The main ingredient is the set of
{\it planar binary trees}. The main result of this section is the vanishing of the dialgebra
homology of a free dialgebra.

An extension of this theory to a theory with coefficients is to be found in [F4].
\M

\noindent {\bf  3.1. Planar binary trees}. 
A planar tree is {\it binary} if any vertex is trivalent. We denote by $Y_n$ the set of
planar binary trees with $n+1$ leaves. Since we only use planar
binary trees in this section we abbreviate it into tree (or $n$-tree to specify that it has $n+1$
leaves, or, equivalently, $n$ interior vertices).
\bigskip

\N $Y_0 =\{\  |\  \}\ ,\  Y_1= \{\  \arbreA \  \}\ ,\quad Y_2=  \{\   \arbreB ,\arbreC \  \},$
\quad 
$Y_3=  \{\  \arbreun ,\arbredeux ,\arbretrois ,\arbrequatre ,\arbrecinq \ 
\}.$
\medskip
\N We will use the permutation-like notation of trees (cf. Appendix A):
\S
$\quad [0] \qquad ,\qquad [1]\qquad ,\qquad [12],\; [21] \qquad ,\qquad [123],\; [213],\;
[131],\; [312],\; [321]$ 
\M
The number of elements in $Y_n$ is the Catalan number $c_n={(2n)!\over n!(n+1)!}$.
For any $y\in Y_n$ we label the $n+1$ leaves by $\{0,1,\cdots  ,n\}$ from left to right.
 
\M

\N {\bf 3.2. Face and degeneracy maps.} For any $i$, $0\le   i\le   n$,
there is a map, called a {\it face map},  $d _i:Y_n\r Y_{n-1}$ which assigns to the tree $y$ the
tree
$d_iy$ obtained from $y$ by deleting the $i$-th leaf. For instance:

$$ d_0[213] = [12],\  d_1[213] = [12],\  d_2[213] = [12],\  d_3[213] = [21].$$

 For any $i$, $0\le   i\le   n$,
there is a map, called a {\it degeneracy map},  $s _i:Y_n\r Y_{n+1}$ which assigns to the tree $y$
the $(n+1)$-tree
$s_iy$ obtained by bifurcating the $i$-th leaf, that is replace it by $\arbreA$ . For
instance 
$$s_0[0] = [1] , s_0[1] = [12] \ , s_1[1] = [21].$$
The face and degeneracy maps satisfy all the classical simplicial relations, {\it except} for the relation
$s_is_i = s_{i+1}s_i$. Indeed, this relation is not fulfilled on trees, because
$$
s_0s_0([0]) = [12] \ {\rm and }\  s_1s_0([0]) = [21]\ .
$$ 
So $Y_.$ is not a simplicial set, but only an {\it almost simplicial set}, (cf. [F3]).
\S

For any $i$, $1\le   i\le   n-1$, there is a map 
$$\circ_i:Y_n\ \r\ \{\dashv,\vdash\}$$
 defined as
follows. The image of $y\in Y_n$ is  $\circ^y_i=\ \dashv$ (resp. $\vdash$) if the $i$-th
leaf points from the vertex to the left (resp. to the right). For instance:
$$
\circ^{[131]}_1 =\  \vdash \quad {\rm and } \quad  \circ^{[131]}_2 =\  \dashv\quad .
$$
More generally one has
$$\circ_i^{[j_1, \cdots ,j_n]} =\ \g\  {\rm if}\  j_i>j_{i+1}\  {\rm and}\ 
\circ_i^{[j_1, \cdots ,j_n]}  =\  \d\  {\rm if}\ 
j_i<j_{i+1},\  {\rm for}\  1\leq i\leq n-1.$$
Here is the table in low dimension:
\medskip

\def\tv{\vrule height 12pt depth 5pt}
$$
\vbox{\offinterlineskip\halign{
\hfill\quad#\quad\hfill\tv&
\hfill\quad#\quad\hfill\tv&
\hfill\quad#\quad\hfill\tv\cr
$y$ &$d_1\qquad \circ_1$ & $d_2\qquad \circ_2$\cr
\noalign{\hrule}
$[12] $ &$[1] \quad \d$ &\cr
$[21] $ &$[1] \quad \g$ &\cr
$[123] $ &$[12] \quad \d$ &$[12] \quad \d$\cr
$[213] $ &$[12]\quad \g$ &$[12] \quad \d$\cr
$[131] $ &$[21] \quad \d$ &$[12] \quad \g$\cr
$[312] $ &$[21] \quad \g$ &$[21]\quad \d$\cr
$[321] $ &$[21] \quad \g$ &$[21] \quad \g$\cr
}}
$$
\M

\noindent {\bf  3.3. The chain complex of a dialgebra.} 
Let $D$ be a dialgebra over $K$. Define the module of $n$-chains by
$$
CY_n(D):=K[Y_n]\otimes D^{\otimes n},
$$
in particular $CY_1(D)\cong D$, $CY_2(D)\cong D^{\otimes 2}\oplus D^{\otimes 2},$
more generally $CY_n(D)$ is isomorphic to the direct sum of $c_n$ copies of $D^{\otimes
n}$ (indexed by $Y_n$).

Define a linear map $d:CY_n(D)\r CY_{n-1}(D)$ by the following formula:
$$
d(y;a_1,\cdots ,a_n):= -\sum^{n-1}_{i=1}(-1)^i(d
_i(y);a_1,\cdots  ,a_{i-1},a_i\circ^y_i a_{i+1},\cdots  ,a_n),
$$
where $y\in Y_n$ and $a_i\in D$.
This formula has a meaning since $\circ^y_i=\ \dashv$ or $\vdash$ and $D$ is a
dialgebra. It is convenient to define
$$
d_i(y;a_1,\cdots  ,a_n):=(d_i(y);a_1,\cdots  ,a_{i-1},a_i\circ^y_ia_{i+1},\cdots  ,a_n)
$$
so that $d=-\sum^{n-1}_{i=1}(-1)^id_i$.
\M

\noindent {\bf  3.4. Lemma.} {\it The face maps $d_i: CY_n(D) \to CY_{n-1}(D)$ satisfy the
simplicial relations
$$d_i d_j = d_{j-1} d_i,\  \hbox {for  any}\  1\leq i< j\leq n-1.$$}

\N {\it Proof}. We first
prove this identity in the lowest dimension, that is
$$
d_1d_2=d_1d_1: CY_3(D) \to CY_2(D)\leqno(*)
$$
The computation of $d_id_j(y;a,b,c)$ splits into 5 cases corresponding to the
five trees with 4 leaves (cf. 3.1).
\medskip
$\bullet$ {\it Case} $[123]$ : 
$$
\eqalign{
&d_1d_2([123] ;a,b,c)=d_1([12] ;a,b\vdash c)=([1] ;a\vdash(b\vdash c)),\cr
&d_1d_1([123] ;a,b,c)=d_1([12] ;(a\vdash b,c)=([1] ;(a\vdash b)\vdash
c).\cr} $$
So relation $(*)$ follows from axiom $5$.
\medskip

$\bullet$ {\it Case} $[213]$ :
$$
\eqalign{
&d_1d_2([213] ;a,b,c)=d_1([12] ;a,b\vdash c)=([1] ;a\vdash(b\vdash
c)),\cr &d_1d_1([213] ;a,b,c)=d_1([12] ;(a\dashv b,c)=([1] ;(a\dashv
b)\vdash c).\cr} $$
So relation $(*)$ follows from axiom $4$.
\medskip

$\bullet$ {\it Case} $[131]$ :
$$
\eqalign{
&d_1d_2([131] ;a,b,c)=d_1([12] ;a,b\dashv c)=([1] ;a\vdash(b\dashv
c)),\cr &d_1d_1([131] ;a,b,c)=d_1([21] ;(a\vdash b,c)=([1] ;(a\vdash
b)\dashv c).\cr} $$
So relation $(*)$ follows from axiom $3$.
\medskip

$\bullet$ {\it Case} $[312]$ :
$$
\eqalign{
&d_1d_2([312] ;a,b,c)=d_1([21] ;a,b\vdash c)=([1] ;a\dashv(b\vdash
c)),\cr &d_1d_1([312] ;a,b,c)=d_1([21] ;(a\dashv b,c)=([1] ;(a\dashv
b)\dashv c).\cr} $$
So relation $(*)$ follows from axiom 2.
\medskip

$\bullet$ {\it Case}  $[321]$ :
$$
\eqalign{
&d_1d_2([321] ;a,b,c)=d_1([21] ;a,b\dashv c)=([1] ;a\dashv(b\dashv
c)),\cr &d_1d_1([321];a,b,c)=d_1([21] ;(a\dashv b,c)=([1] ;(a\dashv
b)\dashv c).\cr} $$
So relation $(*)$ follows from axiom $1$.
\medskip

The proof of the general case $d_id_j=d_{j-1}d_i$ for $i<j$ splits into two different cases.

First, if $j=i+1$, then the proof is exactly as in low dimension and so follows from the
axioms of a dialgebra. Second, if $j>i+1$, then both operations $d_id_j$ and $d_{j-1}d_i$
amount to perform the same modification: removing the leaves number $j$ and $i$ of the
tree $y$, and replace $(a_1,\cdots  ,a_n)$ by
$$
(a_1,\cdots  ,a_i\textstyle\circ^y_ia_{i+1},\cdots  ,a_j\textstyle\circ^y_ja_{j+1},\cdots  ,a_n).
$$
The point is that the leaf number $j$ of $y$ is the leaf number $j-1$ of $d_i(y)$.

So we have proved that $d_id_j=d_{j-1}d_i$ for $i<j$.\hfill  $\square$
\M 

\noindent {\bf  3.5. Proposition.} 
{\it One has $d\circ d=0$ and so $(CY_*(D),d)$ is a chain-complex.}
\medskip

\N {\it Proof.}  This is an immediate consequence of the previous lemma, like for a
pre-simplicial module. \hfill $\square$
\S

Observe that in the chain complex 
$$
\cdots  \r K[Y_n]\otimes D^{\otimes n}\r\cdots  \r 
K[Y_3]\otimes D^{\otimes 3}\r 
K[Y_2]\otimes D^{\otimes 2}\ {\buildrel {(\dashv,\vdash)}\over \longrightarrow}\  D\leqno
{CY_*(D):}
$$
 the module of $n$-chains is the direct sum of $c_n$ copies of $D^{\otimes n}$ (indexed by
the set of trees $Y_n$), the first differential is induced by the two products
$\dashv,\vdash$, and the first relation $d^2=0$ coincides precisely with the 5  
axioms of a dialgebra.
\M

\noindent {\bf  3.6. Homology and cohomology of a dialgebra.} 
 By definition the {\it homology of the dialgebra} $D$ is the homology of the chain-complex
$CY_*(D)$:
$$
HY_n(D):=H_n(CY_*(D),d),\ n\ge 1.
$$
For $n=1$ it is immediate that $HY_1(D)$ is the quotient of $D$ by the submodule generated by all
the elements $x\dashv y$ and $x\vdash y$,
$$
HY_1(D)=D/\{x\dashv y,x\vdash y\ |\ x,y\in D\},
$$
 which we denote, sometimes, by $D/D^2$.

 By definition the {\it cohomology of the dialgebra} $D$ is
$$
HY^n(D):= H^n(\Hom (CY_*(D), K)),\ n\ge 1.
$$
\noindent {\bf  3.7. The chain bicomplex of a dialgebra.} 
The chain complex of a dialgebra is in fact the total chain complex associated to a
bicomplex. Indeed, let $Y_{p,q}$ be the subset of $Y_n$ made of the trees which are obtained
by grafting a $p$-tree with a $q$-tree (cf. Appendix A), where $p+q+1=n$. For instance

$$Y_{0,2} = \{[321], [312]\},\quad Y_{1,1} = \{[131]\},\quad Y_{2,0} = \{[213], [123]\}.$$

Let $CY_{p,q}:=K[Y_{p,q}]\otimes D^{\otimes n}$. Since for any $y\in Y_{p,q}$ the
element  $d_i(y)$ is either in  $Y_{p-1,q}$ or in  $Y_{p,q-1}$, the face map $d_i$ takes value
either in $CY_{p-1,q}$ or in $CY_{p,q-1}$. So the chain bicomplex $CY_{**}(D)$ is
well-defined and its associated total complex is $CY_*(D)$. Remark that, with our choice of
notation, one has $CY_n = \bigoplus _{p+q+1=n} CY_{p,q}$. This bicomplex gives rise to two
 spectral sequences abutting to $HY_*(D)$.
\medskip

\noindent {\bf  3.8. Theorem.}  {\it Let $V$ be a vector space over $K$ and
$Dias(V)=TV\otimes V\otimes TV$ be the free dialgebra over $V$ (cf.2.5). 
Then, one has
}
$$
\eqalign{
&HY_1(Dias(V))\cong V,\cr
&HY_n(Dias(V))=0,\  for \ n>1.\cr}
$$

\N {\it Proof}. The first statement is obvious since $V\cong K\otimes V\otimes K$
is the quotient of $CY_1 = TV\otimes V\otimes TV$ by the submodule generated by all the
products of elements of $V$.

To show that $HY_n = 0$ for $n>1$ we construct a homotopy 
$$
h=h_n: K[Y_n]\otimes
D^{\otimes n} \r  K[Y_{n+1}]\otimes
D^{\otimes n+1}
$$
 such that, for $n>1$,
$$ 
dh_n + h_{n-1}d = \id _n.
$$
In order to write down $h$ explicitly we use the degeneracy maps introduced in 3.2 and the
following construction. Given an $n$-tree $y$, we denote by $p_n(y)$ the $(n+1)$-tree obtained from
$y$ by adding a new leaf at the left of the last one and parallel to it: 

$$\arbreC\quad  \buildrel {p_n}\over {\mapsto } \quad \arbrequatre$$

There are five different formulas for $h_n(x)$ depending on the form of $x\in K[Y_n]\t D^{\t n}$.

\medskip

\noindent $\bullet$ Case (a): $x = (y; \check \omega _1, ..., \check \omega _{n-1},
 \check \omega _n u)$. One puts 
$$h_n(y; \check \omega _1, ..., \check \omega _{n-1}, \check \omega _n u):=
(-1)^n (s _n(y);  \check \omega _1, ..., \check \omega _{n-1}, \check \omega
_n , \check u).$$
First, one has $d_ih_n(x) + h_{n-1}d_i(x) = 0$ for
$1\le   i\le   n-1$ because the modifications on $x$ performed by $h$ and $d$ are disjoint.
Second, one has 
$$(-1)^nd_nh_n(x)=d_n(s _n(y);  \check \omega _1, ..., \check \omega _{n-1}, \check \omega
_n , \check u) = (y; \check \omega _1, ..., \check \omega _{n-1},
 \check \omega _n u) = x$$ since $d_n s _n(y) = y$ and $\circ_n^{s _n(y)} =\  \dashv
.$ So we have proved relation (*) in this case.
\medskip

\noindent $\bullet$ Case (b): $x = (y; \check \omega _1, ..., \check \omega _{n-1},\omega _n
 v \check u)$. One puts 
$$h_n(y; \check \omega _1, ..., \check \omega _{n-1},\omega _n v \check u):= 
(-1)^n(p  _n(y);\check \omega _1, ..., \check \omega _{n-1},\omega _n\check v , \check
u).$$
First, one has $d_ih_n(x) + h_{n-1}d_i(x) = 0$ for
$1\le   i\le   n-1$. Second,
one has 
$$(-1)^nd_nh_n(x)=d_n(p  _n(y);\check \omega _1, ..., \check \omega _{n-1},\omega _n\check v , \check
u) = (y; \check \omega _1, ..., \check \omega _{n-1},\omega _n
 v \check u) = x$$ since $d_n p  _n(y) = y$ and $\circ_n^{p  _n(y)} =\  \vdash .$
So we have proved relation (*) in this case.

\medskip

\noindent $\bullet$ Case (c): $x = (y; \check \omega _1, ..., \check \omega _{n-1},
 \check u)$ and the last two leaves of $y$ have the shape $\arbreA$.  One puts 
$$h_n(y; \check \omega _1, ..., \check \omega _{n-1},\check u):= 
0.$$
First, one has $d_ih_n(x) = 0$ for
$1\le   i\le   n$ and $h_{n-1}d_i(x)=0$ for $1\le   i\le   n-2$. Second, one has 
$$
(-1)^{n-1}h_{n-1}d_{n-1}(x)=(-1)^{n-1}h_{n-1}(d_{n-1}(y);\check \omega _1, ..., \check
\omega _{n-1} u) = x$$ since (by using case (b)) $s _nd_{n-1}(y) = y$ for such $y$.
So we have proved relation (*) in this case.

\medskip

\noindent $\bullet$ Case (d): $x = (y; \check \omega _1, ..., \check \omega _{n-1}v,
 \check u)$ and the last two leaves of $y$ have the shape $\dessinun$. One puts
  $$h_n(y; \check \omega _1, ..., \check \omega _{n-1}v, \check u):= (-1)^n
(s _n(y) - p  _{n-1}(y);\check \omega _1, ..., \check \omega _{n-1},\check v ,
\check u).$$
  Let us write
$h_n = (-1)^n\underline s _n + (-1)^{n-1}\underline p  _{n-1}$.
First, one has $d_ih_n(x) + h_{n-1}d_i(x) = 0$ for
$1\le   i\le   n-2$. Second, one has $d_{n-1}h_n(x) + h_{n-1}d_{n-1}(x) = x$ since 
$d_{n-1}\underline s _n(x) = h_{n-1}d_{n-1}(x)$ and $d_{n-1}\underline p  _{n-1}(x) = x$.
Third, one has $d_nh_n(x)=0$ since $d_n\underline p  _{n-1}(x) = d_n \underline s _n(x)$.
So we have proved relation (*) in this case.

\medskip

\noindent $\bullet$ Case (e): $x = (y; \check \omega _1, ..., \omega _{n-1}\check v,
 \check u)$ and the last two leaves of $y$ have the shape $\dessinun$. One puts
  $$
h_n(y; \check \omega _1, ..., \omega _{n-1}\check v, \check u):= (-1)^n (s _n(y) -
s _{n-1}(y);\check \omega _1, ...,\check \omega_{n-1}\ ,\check v , \check u).
$$
  Let us write
$h_n = (-1)^n\underline s _n + (-1)^{n-1}\underline s _{n-1}$.
First, one has $d_ih_n(x) + h_{n-1}d_i(x) = 0$ for
$1\le   i\le   n-2$. Observe that in many cases $d_ih_n = 0 = h_{n-1}d_i$. Second, one has
$d_{n-1}h_n(x) + h_{n-1}d_{n-1}(x) = x$ since  $ d_{n-1}\underline
s _n(x) = h_{n-1}d_{n-1}(x)$ and $d_{n-1}\underline s _{n-1}(x) = x$. Third, one has
$d_nh_n(x)=0$ since $d_n\underline s _{n-1}(x) = d_n \underline s _n(x)$. So we have
proved relation (*) in this case.\hfill $\square$
\M

\noindent {\bf  3.9. Theorem.}  {\it For any dialgebra $D$ the graded module $HY_*(D)$ is a
graded dual-codialgebra and the graded module $HY^*(D)$ is a graded dendriform  algebra (see section
5). As a consequence  $HY^*(D)$ is a graded associative algebra.}
\S
\N {\it Proof.} Though one could prove these statements directly, they are consequences of
general facts about Koszul operads (see Appendix B). We will show in section 6 that the operad of 
dendriform algebras is dual to the operad of associative dialgebras. Moreover, by theorem
3.8 these operads are Koszul, hence the statement follows from general properties of Koszul
operads (cf. Appendix B5d). The last statement is a consequence of the preceding one and of
Lemma 7.3.
\hfill
$\square$
\M

\N {\bf 3.10. Simplicial properties of the chain-modules.} We have seen in Lemma 3.4 that the
face maps
$d_i: CY_n(D) \to CY_{n-1}(D)$ satisfy the standard simplicial relations. Suppose that D is
equipped with a bar-unit $e$ and let us define $s_j: CY_n(D) \to CY_{n+1}(D) $ by
$$
s_j(y;a_1,\cdots, a_n):= (s_j(y);a_1,\cdots, a_j, e, a_{j+1},\cdots, a_n),\ 0\leq j\leq n,
$$
where $s_j(y)$ is described in 3.2.
From the properties of the bar-unit, it is immediate to check that
$$\eqalign{
d_is_j &= \cases {s_{j-1}d_i& for $i<j$,\cr
{\rm id}&for $i=j,\ i=j+1$,\cr
s_jd_{i-1}& for $i>j+1$,\cr}
\cr
s_is_j&=s_{j+1}s_i \quad  {\rm for}\   i< j.\cr}
$$
So the family $(CY_n(D); d_i, s_j)_{n\geq 0}$ is an {\it almost simplicial module}, that is the face
 and degeneracy operators satisfy all the standard relations of a simplicial module, {\it except}
for the relation $s_is_i = s_{i+1}s_i$  (cf. 3.2).

A variation of the Eilenberg-Zilber theorem is still valid for pseudo-simplicial modules (cf.
Inassaridze [I]) and a fortiori for almost simplicial modules. It is used in the proof of the following result which
is due to Alessandra Frabetti.
\M

 \N {\bf 3.11. Theorem} [F4]. {\it If
$D$ is a dialgebra equipped with a bar-unit, then
$$
HY_n(D) = 0,\  {\it for\ any}\ n\ge 0.
$$}
\hfill $\square$
\S

\N {\it Comment.} This result is similar to the vanishing of the bar-homology for a unital associative
algebra.
\medskip

\N {\bf 3.12. Generalization of the homology of a dialgebra to the symmetric group.} There is a
generalization of the complex $CY_*$ consisting in replacing the set of planar binary trees $Y_n$ by
the symmetric group $S_n$, or, equivalently, by the set $\tilde Y_n$ of {\it binary increasing trees} (cf.
Appendix A). 

The formulas for the maps $d_i$ are the same as in 3.2 and 3.3 once $S_n$ has been identified with  $\tilde
Y_n$ (observe that deleting a leaf in an increasing tree still gives an increasing tree). Hence we get a new complex
$$
\cdots  \r K[S_n]\otimes D^{\otimes n}\r\cdots  \r 
K[S_3]\otimes D^{\otimes 3}\r 
K[S_2]\otimes D^{\otimes 2}\ {\buildrel {(\dashv,\vdash)}\over \longrightarrow}\  D\leqno
{CS_*(D):}
$$
for any dialgebra $D$, and new homology groups $HS_*(D)$.

The boundary map is still the alternate sum of face maps. When the dialgebra is bar-unital, then
there also exist degeneracy maps. All these maps satisfy the simplicial relations {\it except}
the relations involving only the degeneracy maps. Such an object is called a {\it
pseudo-simplicial} module (cf. [I]).

 Forgetting the levels gives a map
$\Psi: S_n =  {\tilde Y_n}
\to Y_n$ (cf. Appendix A) which induces a chain map
$CS_*(D) \to CY_*(D)$ and hence a morphism
$$
HS_*(D) \to HY_*(D).
$$
This new theory $HS_*$, or, more accurately, its variant with non trivial coefficients, crops up
naturally when one wants to compute the Leibniz homology of the Leibniz algebra of matrices $\gg
l(D)$ over a dialgebra $D$ (cf. Frabetti [F2]).
\medskip

\N {\bf 3.13. Remark.} We will show in section 6 that $CY_*(D)$ is the chain complex of the
dialgebra $D$ predicted by the operad theory. One could wonder if there is another notion of
algebra for which $CS_*(D)$ would be the predicted chain complex. If this would be the case, then
the Poincar\'e series of the dual operad would be the inverse of the series
$$
\sum_{n\geq 1} (-1)^n \# S_n x^n = \sum_{n\geq 1} (-1)^n n! x^n.
$$
However this inverse is not of the form $\sum_{n\geq 1} (-1)^n a_n x^n$ with $a_n\in
{\NN }$, hence this operad, even if it existed, could not be a Koszul operad.
\M



\vfill\eject



\noindent {\bf 4.  LEIBNIZ ALGEBRAS, ASSOCIATIVE DIALGEBRAS AND HOMOLOGY}
\B
A Leibniz algebra is a non-commutative version of a Lie algebra. In this section we show that, when
we replace Lie algebras by Leibniz algebras, then the role of associative algebras is played by  the
associative dialgebras. In particular we show that any Leibniz algebra has a {\it universal
enveloping associative dialgebra}.
\M

\noindent {\bf 4.1. Leibniz algebras} [L1], [LP].
Recall that a {\it Leibniz algebra} over $K$ is a $K$-module $\gg$ equipped
with a binary operation (called a bracket): 
$$[-,-]:\gg\otimes \gg\r \gg ,$$
which satisfies the {\it Leibniz identity}:
$$
[x,[y,z]]=[[x,y],z]-[[x,z],y],
$$
for all $x,y,z$ in $\gg$. This is in fact a {\it right} Leibniz algebra. For the opposite
 structure, that is $[x,y]' = [y,x]$, the left Leibniz 
identity is 
$$
[[x,y]',z]' = [x,[y,z]']' -  [y,[x,z]']'.
$$
If the bracket happens to be anticommutative, then $\gg$ is a
Lie algebra. Quotienting the Leibniz algebra $\gg $ by the ideal generated by the elements $[x,x]$
for all $x\in \gg$ gives a Lie algebra that we denote by $\gg _{Lie}$.

To any Leibniz algebra $\gg$ is associated a chain-complex
$$\ldots\longrightarrow \gg ^{\otimes n}\buildrel {d}\over
\longrightarrow \gg ^{\otimes
n-1}\buildrel {d}\over
\longrightarrow\ldots\buildrel {d}\over
\longrightarrow \gg ^{\otimes 2}\buildrel {d}\over
\longrightarrow\gg \leqno CL_*(\gg ):$$
where
$$
d(x_1,\ldots ,x_n)=\sum_{1\leq i<j\leq n}(-1)^j(x_1,\ldots
,x_{i-1},[x_i,x_j],x_{i+1},\ldots ,\hat
x_j,\ldots ,x_n).
$$
The homology groups of this complex are denoted by $HL_n(\gg)$, for $n\geq 1$.
\medskip

\noindent {\bf 4.2. Proposition.}
{\it  Let $D$ be a dialgebra. Then the bracket 
$$[x,y]:=x\dashv y-y\vdash x$$
makes $D$
into a Leibniz algebra, denoted by $D_{Leib}$.}
\medskip

{\it Proof}. It is a straightforward checking in which axioms (1),
(2), (4), (5) are used once and axiom (3) twice since
$$
\eqalign{
[x,[y,z]] \ &= \ x\g (y\g z) - x\g (z\d y) - (y\g z)\d x + (z\d y)\d x , \cr
-[[x,y],z]] &= -(x\g y)\g z + (y\d x)\g z + z\d (x\g y) - z\d (y\d x) ,\cr
[[x,z],y]] \ &= \ (x\g z)\g y - (z\d x)\g y - y\d (x\g z) + y\d (z\d x) . \cr 
} $$
\hfill $\square$

This construction defines a functor 
$$ {\bf Dias} \buildrel{-}\over \longrightarrow {\bf Leib}
$$
from the category {\bf Dias} of dialgebras to the category   {\bf Leib} of 
Leibniz algebras.
\medskip

\N {\bf 4.3. Example.} For any dialgebra $D$ over $K$ the $n\times n$-matrices with entries in  $D$ 
form a new dialgebra ${\cal M}_n(D)$. Its associated Leibniz algebra is denoted
$\gg l_n(D)$. It is not a Lie algebra in general. The homology of $\gg l(D)$, when
$D$ has a bar-unit, has been computed by Frabetti [F2].
\M

\noindent {\bf 4.4. Proposition.} {\it
The following diagram of categories and functors is commutative}
$$
\matrix{
\hbox{{\bf Dias}}&{\buildrel - \over \longrightarrow} &\hbox{{\bf Leib}}\cr
\uparrow &&\uparrow \cr
\hbox{{\bf As}}&{\buildrel - \over \longrightarrow} &\hbox{{\bf Lie}.}\cr
}
$$
\hfill $\square$
\M

\noindent {\bf 4.5. Remark.}
In the definition of a dialgebra, axiom (3) could be relaxed slightly, though proposition
4.2 remains valid. It could be replaced by the weaker axiom:
$$
((y\vdash x)\dashv z)+(z\vdash (x\dashv y))=(y\vdash (x\dashv z))+((z\vdash x)\dashv
y). $$
Observe that in this formula the variables do not stay in the same order in  the
monomials. Hence the associated operad would not be a non-$\SS$-operad anymore.
\medskip

\noindent {\bf 4.6. Universal enveloping associative dialgebra of a Leibniz algebra.}
The functor $-:\hbox{{\bf As}}\r \hbox{{\bf Lie}}$ has a left adjoint which is
the universal enveloping algebra of a Lie algebra:
$$
{\overline U}(\gg)={\overline T}(\gg)/\{[x,y]-x\otimes y+y\otimes x\mid x,y\in \gg\}.
$$
(${\overline U}(\gg)$ is the augmentation ideal of the classical enveloping unital algebra
$U(\gg)$). Similarly, define the {\it universal enveloping dialgebra} of a Leibniz algebra $\gg$
as the following quotient of the free dialgebra on $\gg$:
$$
Ud(\gg):=T(\gg)\otimes \gg\otimes T(\gg)/\{[x,y]-x\dashv y+y\vdash x\mid
x,y\in \gg\}. $$

Under our previous notation, the elements generating the ideal are denoted by
$[x,y]\check{} -\check xy+y\check x$.
\medskip

 \noindent {\bf 4.7. Proposition.} {\it The functor $Ud:\hbox{\bf 
Leib}\r  \hbox{{\bf Dias}}$  is left adjoint to the functor} 
$-:\hbox{\bf Dias}\r\hbox{\bf Leib}$.
\medskip

\N {\it Proof.} Let $f:\gg\r D_{Leib}$ be a morphism of Leibniz algebras. There is a unique
extension of $f$ as a morphism of dialgebras from $T(\gg)\otimes \gg\otimes T(\gg)$ to
$D$. Since the image of $[x,y]-x\dashv y+y\vdash x$ under this morphism is $0$, it defines
a morphism from $Ud(\gg)$ to $D$.

On the other hand the restriction of the morphism of dialgebras $g:Ud(\gg)\r D$ to $
\gg=K\otimes \gg\otimes K$ yields a morphism of Leibniz algebras $\gg\r
D_{Leib}$.

It is now immediate to check that these two constructions give rise to isomorphisms
$$
\Hom_{\bf{ Dias}}(Ud(\gg),D)\cong \Hom_{\bf {Leib}}(\gg,D_{Leib}).
$$
\hfill $\square$
\S
It is well-known that the universal enveloping algebra of a Lie algebra is not only an
associative algebra but a Hopf algebra. Similarly the universal enveloping dialgebra of a
Leibniz algebra possesses co-operations. They are studied by Goichot in [Go].
 \M

\N {\bf 4.8. Lemma.} {\it For any Leibniz algebra $\gg $, one has $Ud(\gg)_{As} = U(\gg_{Lie})$.}
\S
\N {\it Proof.} Since the functor $(-)_{As}: {\bf Dias} \to {\bf As}$ is left adjoint to $inc:
{\bf As} \to  {\bf Dias}$ and since  $(-)_{Lie}: {\bf Leib} \to {\bf Lie}$ is left adjoint to
$inc: {\bf Lie }
\to  {\bf Leib}$ both composites $U \circ (-)_{Lie} $ and $(-)_{As} \circ Ud$ are left adjoint
to the composite $ {\bf Leib} \to {\bf Lie} \to {\bf As}$, and so are equal.\hfill $\square$
\M

\N {\bf 4.9. Proposition.} {\it The universal enveloping dialgebra $Ud(\gg)$ is isomorphic to
$U(\gg _{Lie})\t \gg$, equipped with the dialgebra structure issued from a $U(\gg_{Lie})$-bimodule
structure and the bimodule map $U(\gg_{Lie})\t \gg \to U(\gg_{Lie})$ (cf. example 2.2.d).}
\S

\N {\it Proof.} Let us define a $U(\gg _{Lie})$-bimodule structure on $U(\gg _{Lie})\otimes \gg$.
The left module structure is given by multiplication in the left factor. The right module
structure is induced by 
$$
(\oo \t x)\cdot \bar y:= \oo \t [x,y] + \oo \bar y \t x ,
$$
where $\oo \in U(\gg _{Lie}), x\in \gg , \bar y \in \gg _{Lie}$ and $y\in \gg$ is a lifting of
$\bar y \in \gg _{Lie}$. It is a well-defined element because the bracket $[x,y]$ in the Leibniz
algebra $\gg$ depends only on the class of $y$ in  $\gg _{Lie}$. Let us check that this formula
provides a representation of $\gg _{Lie}$.

Let $\bar y$ and $\bar z$ be elements in $\gg _{Lie}$ and $y, z$ be liftings in $\gg$. On one hand
one gets 
$$
((\oo\t x \cdot \bar y)\cdot \bar z =
\oo \t [[x,y],z] + \oo \bar z \t [x,y] + \oo \bar y \t
[x,z] + \oo \bar y \bar z \t x, 
$$
and 
$$\
((\oo\t x \cdot \bar z)\cdot \bar y =
\hfill  \oo \t [[x,z],y] + \oo \bar y \t [x,z] + \oo \bar z \t
[x,y] + \oo \bar z \bar y \t x. 
$$
Hence one has 
$$\eqalign {
((\oo\t x \cdot \bar z)\cdot \bar y -((\oo\t x \cdot \bar z)\cdot \bar y 
&=\oo\t ([x,y],z] - [x,z],y]) - \oo (\bar y \bar z - \bar z \bar y )\t x \cr
&= \oo \t [x, [y,z]] - \oo \overline { [y,z] } \t x\cr
&= (\oo\t x) \cdot [y,z]\, .
 \cr}
$$
The right and left module structures are immediately seen to be compatible, hence $U(\gg _{Lie})\t
\gg$ is a $U(\gg _{Lie})$-bimodule.

The linear map $U(\gg _{Lie})\t \gg \to U(\gg _{Lie}), \oo\t x \mapsto \oo \bar x$ is a $U(\gg
_{Lie})$-bimodule because $\oo \bar x \bar y = \oo \overline {[x,y]} + \oo \bar y \bar x$.

So, it follows that $U(\gg _{Lie})\t \gg$ is equipped with a dialgebra structure (cf. 22.d). The
(nonunital) associative algebra associated to this dialgebra is the augmentation ideal of $U(\gg
_{Lie})$, cf. 2.6.

There is a well-defined dialgebra map
$$
Ud(\gg) \to  U(\gg _{Lie})\t \gg
$$
which sends $\oo\t x \t 1$ to $\bar {\oo } \t x$ (for $\oo \in T(\gg)$ and $\bar {\oo}$ its image
in
$U(\gg _{Lie})$). Indeed, any element in $Ud(\gg)$ can be written as a linear combination of
elements of the form $\oo\t x \t 1$ since 
$$
\oo\t x \t y = \oo \t [x,y]\t 1 + \oo y\t x\t 1.
$$
Since by lemma 4.8 one has $Ud(\gg)_{As} = U(\gg_{Lie})$, it follows that the element $\oo \check
x \oo '$ in $Ud(\gg)$ depends only on the class of $\oo \in T(\gg)$ (resp.  $\oo '\in
T(\gg))$ in  $ U(\gg_{Lie})$. So, on can define a dialgebra map 
$$
 U(\gg_{Lie})\t \gg \to Ud(\gg )
$$
by sending $\bar {\oo }\t x$ to $\oo \check x$, where $\oo$ is a lifting of $\bar {\oo}$.

It is immediate to check that both composites are the identity, whence the isomorphism.
 \hfill $\square$
\M

\noindent {\bf 4.10. Free algebras and free dialgebras.}
Let $V$ be a $K$-module and let $\overline{T}(V)=V\oplus V^{\otimes 2}\oplus\cdots \oplus
V^{\otimes n}\oplus\cdots \quad $ be the tensor module. We denote by $\cc $ the endomorphism of
$\overline{T}(V)$ defined inductively by $\cc (v)=v$ for $v\in V$ and $\cc (\oo \otimes v)=\oo \otimes
v-v\otimes \oo$ for $\oo \in V^{\otimes n}$ and $v\in V$. It is well-known that $\Im\cc $ is
isomorphic to the free Lie algebra $Lie(V)$ over $V$. Recall that the free
associative algebra over $V$ is $\overline{T}(V)$ equipped with the concatenation product, and
the free Leibniz algebra over $V$ is $\overline{T}(V)$ equipped with the unique Leibniz bracket
which satisfies $[\oo ,v]=\oo \otimes v$ for $\oo \in V^{\otimes n}$ and $v\in V$ (cf. [LP]). In
the sequence 
$$
\overline{T}(V)\buildrel{\cc }\over\rightarrow Lie(V)\hookrightarrow \overline{T}(V),
$$
the first map is a map of Leibniz algebras, the second one is a map of Lie algebras (for
the Lie structure of $\overline{T}(V)$ coming from its associative algebra structure). From
Proposition 4.4 it follows that there is a commutative diagram
$$
\matrix{
Leib(V)=\overline{T}(V)&\buildrel{\cc ^{\check{}}_1}\over\longrightarrow&
\overline{T}(V)\otimes V\otimes \overline{T}(V)&=Dias(V)\cr
&&&\cr
 \gamma \downarrow\quad &&\qquad\quad\downarrow {fusion}\cr
&&&\cr
Lie(V)&\hookrightarrow &\overline{T}(V)&=As(V)\cr }
$$
where the maps {\it fusion} and $\cc ^{\check{}}_1$  are described as follows. 

The fusion map
consists in forgetting the symbol ${\check{}}\ $, that is $\oo \check v\oo '\mapsto \oo v\oo '$. 

The image of
$v_1\cdots v_n$ by $\cc ^{\check{}}_1$ is $\cc (\check v_1v_2\cdots v_n)$, which means writing
$\cc (v_1\cdots v_n)$ and putting the symbol ${\check{}}$ on the variable $v_1$ of each monomial. For
instance
$\cc ^{\check{}}_1(v_1)=\check v_1,\cc ^{\check{}}_1(v_1\otimes
v_2)=\check v_1v_2-v_2\check v_1,\cc ^{\check{}}_1(v_1\otimes v_2\otimes
v_3)=\check v_1v_2v_3-v_2\check v_1v_3 -v_3\check v_1v_2+v_3v_2\check v_1$. 
\S

One observes that this map is very similar to the map $\aa$ used by Lodder in [Lo] to describe the
loop suspension of a wedge product of topological spaces.
\M

\noindent {\bf 4.11. Proposition.} {\it Let $V$ be a K-module and $Leib(V) = \overline{T} (V)$ be
the free Leibniz algebra on $V$. Then one has an isomorphism}
$$ Ud(Leib(V)) \cong Dias(V) = T(V)\otimes V \otimes T(V). $$
\medskip

\N {\it Proof}. The functor $Leib$ is left adjoint to the forgetful functor from Leibniz
algebras to modules. Similarly the functor $Ud$ is left adjoint to the functor from
dialgebras to Leibniz algebras, therefore the composite is left adjoint to
 the forgetful functor from dialgebras to modules, so it is the free diassociative algebra functor.
\hfill $\square$
\medskip

\N {\bf 4.12. Theorem.} {\it
For any dialgebra $D$ there is a natural transformation
$$
 HL_*(D_{Leib})\longrightarrow HS_*(D)
$$
induced by the chain complex map 
$$\displaylines{
\ee_n: CL_n(D_{Leib})=D^{\otimes n} \longrightarrow K[S_n]\otimes D^{\otimes n}= CS_n(D),\cr
\ee_n (x_1, \cdots, x_n) = \sum _{\ss \in S_n} {\rm sgn}(\ss) \sigma \t
\sigma ^{-1}(x_1, \cdots ,x_n).\cr}
$$}
 \S
\N {\it Proof.} The boundary map $d_L$ of the Leibniz complex $CL_*$ is described in 4.1. For any
element $( {\underline x}, y):= (x_1, \ldots , x_n, y) \in D^{\t n+1}$ the map $d_L$ is defined
inductively by the formula:
$$
d_L({\underline x}, y)= (d_L( {\underline x}), y) + (-1)^n \ad (y) ({\underline x}),\leqno
(4.12.1)
$$
where $\ad (y) ({\underline x}) = \ad (y)(x_1, \ldots , x_n):= \sum_{i=1}^n (x_1, \ldots ,
x_i \g y - y \d x_i, \ldots ,  x_n)$.

The boundary map $d_S$ of the symmetric complex $CS_*$ is described in 3.12 and 3.2.

One extends the operator $\ad (y)$  to $K[S_n] \t  D^{\t n}$ by putting $\ad (y) (\ss
\t {\underline x}) = \ss \t \ad (y) ({\underline x})\in K[S_n] \t  D^{\t n}$, so that it
obviously commutes with $\ee_n$:
$$\ee_n \ad (y) = \ad (y) \ee_n.\leqno (4.12.2)$$
It turns out that this new operator is homotopic to 0. The homotopy $h(y): K[S_n] \t  D^{\t n}
\to K[S_{n+1}] \t  D^{\t n+1}$ is given by 
$$h(y)(\ss\t {\underline x}):= \sum_{i=0}^n (-1)^i s_i(\ss) \t (x_1, \ldots ,
x_i, y,  x_{i+1}, \ldots ,  x_n),$$
where $ s_i(\ss) $ is the $i$-th degeneracy of $\ss$ (bifurcate the $i$-th leaf of the
corresponding increasing tree, cf. 3.2).
The checking of
$$d_S h(y) + h(y) d_S = \ad (y)\leqno (4.12.3)
$$
is tedious but straightforward.

The comparison of the homotopy operator $h(y)$ with the symmetrization operator $\ee_n$ gives:
$$\ee_{n+1}({\underline x}, y) = (-1)^n h(y) \ee_n ({\underline x}).\leqno (4.12.4)$$
\S 
The proof of the commutation relation
$$
d_S \circ \ee_n = \ee_{n-1} \circ d_L \leqno (*)_n
$$
is done by induction on $n$ as follows.

For $n=1$, one has $\ee_1 =\Id$.

 For $n=2$, the map $\ee_2: D^{\t 2} \to  K[S_2]\t  D^{\t 2}$
is given by 
$$\ee_2(x,y) = [12]\t (x,y) - [21]\t (y,x).$$
 Since $d_L(x,y) = [x,y]$ and $d_S([12]\t (x,y)) = x\g
y
$, 
 $d_S([21]\t (x,y))= x\d y $, it follows from $[x,y] = x\g y - y \d x$ that  $d_S \circ \ee_2 =
\ee_1
\circ d_L $.
\M

By induction we suppose that $(*)_n$ holds and we will prove $(*)_{n+1}$.

One gets:
$$
\matrix{
d_S \ee_{n+1}({\underline x}, y)&=(-1)^n d_S h(y) \ee_n({\underline x})\hfill & by\ 
(4.12.4)\hfill\cr
     &=(-1)^n(\ad (y) -h(y)d_S) \ee_n({\underline x})\hfill& by\  (4.12.3)\hfill\cr
     &=(-1)^n \ad (y) \ee_n({\underline x}) + (-1)^{n-1}   h(y) \ee_{n-1} d_L({\underline x})
 \hfill&by\  induction\hfill\cr
     &=(-1)^n \ad (y) \ee_n({\underline x}) +\ee_n(d_L({\underline x}),y)\hfill & by\ 
(4.12.4)\hfill\cr
     &=(-1)^n  \ee_n({\underline x}) \ad (y) +\ee_n(d_L({\underline x}),y)\hfill & by\ 
(4.12.2)\hfill\cr
     &= \ee_n d_L ({\underline x},y) \hfill
& by\ 
(4.12.1).\hfill\cr}
$$
\hfill $\square$
\M
\N {\bf Remark.} This proof is mimicked on the proof for the Lie case as done in [L0],
Proposition 1.3.5. This Proposition has been extended to homology of dialgebras with
coefficients by Alessandra Frabetti in her thesis (unpublished).
\M

\N {\bf 4.13. Proposition.} {\it The composite 
$$
 HL_*(D_{Leib})\rightarrow HS_*(D)\buildrel {\psi_*} \over {\longrightarrow } HY_*(D),
$$
where $\psi_n:S_n \proj Y_n$ is the surjective map described in Appendix A, 
is the map induced, through the operad theory, by the morphism of operads $ Dias^! \to Leib^!
$.}
\S

\N {\it Proof.} We show in section 8 that the complexes $CL_*$ and $CY_*$ are the complexes
predicted by the operad theory. So, by Appendix B5e, it suffices to check that the natural map
$ K[Y_n]
\t K[S_n] = Dias^!(n)\to Leib^!(n)= K[S_n]$ on the dual operads (see section 8) is given by $y\t \oo \mapsto
\sum_{\{\ss\in S_n
\mid
\psi(\ss) = y\}} \ss \oo\ $, and this is precisely Theorem 7.5. \hfill $\square$
\M

\N {\bf 4.14. Comparison of $HL_*(\gg)$ and $HY_*(Ud(\gg))$ for a Leibniz algebra $\gg$.}
For a Lie algebra $\h$ the composite map
$$
H_*^{Lie}(\h) \to  H_*^{Lie}({\bar U}(\h)_{Lie})\rightarrow H_*^{As}({\bar U}(\h)),
$$
 is known to be an isomorphism (cf. for instance [L0]).

However, for a Leibniz algebra
$\gg$, the composite map
$$
HL_*(\gg) \to  HL_*(Ud(\gg)_{Leib})\rightarrow HY_*(Ud(\gg))
$$
is no longer an isomorphism, contrarily to what was mistakenly announced in [L2]. The main point in the proof of the Lie case, is that the graded space associated to the
filtration of $U(\gg)$ depends only on the vector space $\gg$, and not on the
Lie structure. In the Leibniz case, the graded space associated to $Ud(\gg) =
U(\gg_{Lie})\t \gg$  (cf. Proposition 4.9) depends on $\gg_{Lie}$, that is, on the Leibniz
structure of $\gg$. Hence the map  $HL_*(\gg) \to HY_*(Ud(\gg))$ is part of a spectral sequence
involving the derived functors of ``Lie-zation".
\M

\noindent {\bf 4.15. Poisson dialgebra.} By definition a
{\it Poisson dialgebra} $P$ is a vector space $P$ equipped with a dialgebra structure
$\dashv$ and $\vdash$, and a Leibniz structure $[-,-]$ which are compatible in the
sense that they satisfy the following 4 relations:
$$[x, y \dashv z] = y \vdash [x,z] + [x,y] \dashv z = [x, y \vdash z],$$
$$[x \dashv y, z] = x \dashv [y,z] + [x,z] \dashv y,$$
$$[x \vdash y, z] = x \vdash [y,z] + [x,z] \vdash y.$$
This definition generalizes the ``non-commutative Poisson algebra" as defined in [Kub],
[KS] and [Ak].

\vfill\eject



\noindent {\bf 5. DENDRIFORM ALGEBRAS}
\medskip

In this chapter we construct a new type of algebra with two binary operations, which
dichotomizes the notion of associative algebra. It is closely related to associative dialgebras. In
fact, we show in section 8 that its operad is  Koszul dual to the operad of associative dialgebras.
The terminology is due to the structure of the free dendriform algebra, which is best
described in terms of planar binary trees.
\medskip

\noindent {\bf 5.1. Definition.}
A {\it dendriform algebra} $E$  over $K$ is a
$K$-vector space $E$ equipped with two binary operations
$$
\eqalign{
&\prec\ : E\otimes E\r E,\cr
&\succ\ : E\otimes E\r E,\cr}
$$
which satisfy the following axioms:
$$
\eqalign{
\hbox{(i)}\quad &(a\prec b)\prec c=a\prec(b\prec c)+a\prec(b\succ c),\cr
\hbox{(ii)}\quad& (a\succ b)\prec c=a\succ (b\prec c),\cr
\hbox{(iii)}\quad& (a\prec b)\succ c+(a\succ b)\succ c=a\succ (b\succ c),\cr}
$$
for any elements $a, b$ and $c$ in $E$.
\M
It is sometimes preferable  to call this object   {\it dendriform dialgebra} to insist on the fact that it is
defined by two operations, but we do not use this terminology here.  It is important to observe that, like for
associative algebras and associative dialgebras, the monomials
 involved in the relations keep the variables in the same order. As a consequence the associated operad is a
non-$\Sigma$-operad.

Observe that there is no ``monoid" version of dendriform algebras, since relations (i) and (iii)
involve sums. In other words, the associated operad does not come from a set-operad.

By introducing the operation 
$$x*y:= x\prec y + x \succ y,$$
 these relations take the following  more concise form:
$$
\eqalign{
\hbox{(i)}\quad &(a\prec b)\prec c = a\prec(b* c)  ,\cr
\hbox{(ii)}\quad& (a\succ b)\prec c = a\succ (b\prec c)  ,\cr
\hbox{(iii)}\quad& (a*b)\succ c = a\succ (b\succ c)  .\cr}
$$

\noindent  {\bf 5.2. Lemma.} {\it For any  dendriform  algebra $E$  the product defined by 
$$ 
x* y:=x\prec y+x\succ y.
$$ 
 is associative.}
\S
\noindent  {\it Proof.} Adding up the three equalities (i), (ii) and (iii) we get $(x*y)*z$ on the left hand side
and 
$x*(y*z)$ on the right hand side, whence the statement. 
 \hfill $\square$
\M

It follows from this lemma that a dendriform algebra is in fact an associative algebra, whose product has
some special property.  The category of dendriform algebras is denoted by ${\bf Dend}$. 
\medskip

\noindent {\bf 5.3. Proposition.}
{\it Let $D$ be a dialgebra and $E$ a
dendriform algebra. Then, on the tensor product
$D\otimes E$, the bracket
$$
\eqalign{
[x\otimes a,y\otimes b]:=&(x\dashv y)\otimes (a\prec b)-(y\vdash x)\otimes (b\succ a)\cr
&-(y\dashv x)\otimes (b\prec a)+(x\vdash y)\otimes (a\succ b),\cr}
$$
where $x,y\in D,a,b\in E$, defines a structure of Lie algebra.}

\medskip

\N {\it Proof}. The bracket is antisymmetric by definition. Hence, it suffices to show that the
Jacobi identity is fulfilled.

The Jacobi identity for $x\otimes a$, $y\otimes b$, $z\otimes c$ gives a total of 48
terms, in fact $8\times 3!$ terms. There are 8 terms for which $x,y,z$ (and also $a,b,c$)
stay in the same order. The other sets of 8 terms are permutations of this set which reads:
$$
\eqalign{
&x\dashv (y\dashv z)\otimes a\prec(b\prec c)\ -\ (x\dashv y)\dashv z\otimes (a\prec
b)\prec c,\cr
&x\vdash (y\dashv z)\otimes a\succ(b\prec c)\ -\ (x\vdash y)\dashv z\otimes (a\succ
b)\prec c,\cr
&x\dashv (y\vdash z)\otimes a\prec(b\succ c)\ -\ (x\dashv y)\vdash z\otimes (a\prec
b)\succ c,\cr
&x\vdash (y\vdash z)\otimes a\succ(b\succ c)\ -\ (x\vdash y)\vdash z\otimes (a\succ
b)\succ c.\cr}
 $$
The terms 1 and 3 in column 1 together with the term 1 in column 2 cancel due to axioms
(1), (2) and (i). Similarly the terms 41, 32 and 42 cancel due to axioms (4), (5) and
(iii). Finally the terms 21 and 22 cancel due to axioms (3) and (ii).\ $\square$
\M

This result may also be seen as a consequence of Koszul duality (cf. Appendix B).
\medskip

\noindent {\bf 5.4. Examples of dendriform algebras.}

\N {\it (a) Shuffle algebra.} Let $V$ be a vector space and let $T'(V)$ be the reduced tensor module over $V$
equipped with the shuffle product (which is associative and commutative). The shuffle of two generating
elements
$v_1 \cdots v_p$ and $v_{p+1}\cdots v_{p+q}$ can be split into two parts depending on the fact that the first
element is $v_1$ or $v_{p+1}$. The first part gives the left product and the second part gives the right
product. One can show that the shuffle algebra is then a dendriform algebra (this fact had been
previously remarked by Gian-Carlo Rota).

\N {\it (b) Matrices over dendriform algebras.} Since  in the
axioms of a dendriform algebra the variables $a,b,c$
stay in this order in all the monomials, the tensor product of two dendriform algebras is
naturally a dendriform algebra. Similarly,  let ${\cal M}_n(E)$
be the module of
$n\times n$-matrices with entries in the dendriform algebra $E$. Then the formulas
$$
(\alpha \prec \beta )_{ij}=\sum_k\alpha _{ik}\prec \beta _{kj}\quad \hbox{and}\quad (\alpha \succ
\beta ) _{ij}=\sum_k\alpha _{ik}\succ
\beta _{kj} $$
make ${\cal M}_n(E)$ into a 
dendriform algebra.

\N {\it (c) Free dendriform algebra.}
Let $V$ be a $K$-module and denote by $ Dend(V)$ the free dendriform algebra over
$V$. It is a dendriform algebra which satisfies the classical universal property. We will prove its
existence and give an explicit description in 5.7. As a first step we describe the free
dendriform algebra on one generator by using the sets of planar binary trees $Y_n$, and some
operations on them.
\M

\N {\bf 5.5. Grafting operation on trees.} By definition the  {\it grafting} of  the
trees $y\in Y_p$ and $z\in Y_q$  is the tree $y\vee z \in Y_{p+q+1}$ obtained from
$y$ and $z$ by joining their roots together and adding a new root. Observe that the number of
internal vertices of  $y\vee z$ is the sum of the numbers of internal vertices of $y$ and of $z$
plus 1.

Given a tree $y$ (different from $\vert$ ) there is a unique decomposition $y=y_1\vee y_2$. For
instance one has:
$$\arbreA =\  \vert \vee \vert ,\quad  \arbreB = \arbreA \vee \vert ,\quad   \arbreC =\ 
\vert \vee \arbreA .$$
which, with our notation, reads 
$$[1] = [0]\vee [0],\quad [12]=[1]\vee [0], \quad [21] =
[0]\vee [1].$$ 

The grafting operation is easy to write down in terms of the permutation-like notation. Indeed,
for $y\in Y_p$ and $z\in Y_q$, one has
$$
[y] \vee [z] = [y\  p+q+1\ z].
$$
For instance one has:
$$[1]\vee [1] = [1\  3\  1],\qquad  [1\  3\  1]\vee
 [2\  1] = [1\  3\  1\  6\  2\  1].$$

We put $K[Y_{\infty}]:= \oplus_{n\geq 0} K[Y_n] $
and 
 $\overline {K[Y_{\infty}]}:= \oplus_{n\geq 1} K[Y_n] $.
We introduce recursively the following operations on $K[Y_{\infty}]$:
$$\leqalignno {
y\prec z &:= y_1 \vee (y_2 * z),& (5.5.1)\cr
y\succ z &:= (y * z_1)\vee z_2,& (5.5.2)\cr
y * z &:= y\prec z + y\succ z , &(5.5.3)\cr
x\prec \vert &:= x =: \vert \succ x \hbox { and}\  x\succ \vert := 0 =: \vert \prec x, \hbox { for } x\neq \vert. &
(5.5.4)\cr}
$$
for $y\in Y_p$ and $z\in Y_q$. Observe that $\vert$ is a unit for $*$.

Since the decomposition $y = y_1 \vee y_2$ is unique, it is clear that these formulas are
well-defined by recursion. For instance
$$
\arbreA \prec \arbreA = \ \vert \vee (\ \vert * \arbreA ) = \  \vert \vee \arbreA = \arbreC \ ,
$$
$$
\arbreA \succ \arbreA = ( \arbreA * \vert\  ) \vee \vert \ =  \arbreA \vee \vert\  = \arbreB \ ,
$$
or, equivalently, $[1] \prec [1] = [21] , \quad [1]\succ [1]=[12]$.

These operations are extended  to $K[Y_{\infty}]$ by linearity. Observe that $[0]*[0] = [0]$, but $[0]\prec
[0]$ and  $[0]\succ[0]$ are not defined.
\M

\N {\bf 5.6. Lemma.} {\it The vector space $\oplus _{n\geq 1} {K[Y_n]}$ equipped with the two
operations $\prec$ and $\succ$ described above is a dendriform algebra, which is generated by $[1] =
\arbreA$.}
 \S
\N {\it Proof.}  We prove this assertion by induction on the (total) degree of the
trees.

Let $y=y_1\vee y_2,\ y\rq=y\rq _1\vee y\rq _2 $ and $y''=y''_1\vee 
y''_2$ be planar binary trees. The following equalities follow by induction from the definitions
of the operations and the associativity of $*$: $$\matrix{
(i)\hfill &(y\prec y\rq )\prec y'' = (y_1\vee (y_2*y\rq ))\prec y'' \hfill\cr
 &=y_1\vee (y_2*y\rq *y'')=y\prec (y\rq \prec y'')\ +\ y\prec (y\rq \succ y'').\hfill\cr
&\cr
(ii)\hfill  & y\succ (y\rq \prec y'') = y\succ (y\rq _1\vee 
(y\rq _2*y''))=(y*y\rq _1)\vee (y\rq _2*y'')\hfill \cr
 &=((y*y\rq _1)\vee y\rq _2) \prec y'' = (y\succ y\rq )\prec y''.\hfill\cr
&\cr
(iii)\hfill & y\succ (y\rq \succ y'') =y\succ ((y\rq *y''_1) \vee y''_2)\hfill\cr
&=(y*y\rq *y''_1)\vee y''_2 =(y\prec y\rq)\succ y'' + (y\succ y\rq)\succ y''.\hfill\cr
}$$
Let us show that $\overline {K[Y_{\infty}]}$ is generated by $[1]$ under the operations $\prec$ and $\succ$. 
Let  $y=y_1\vee y_2$ be a tree. From the definitions of the operations we have 
$$\eqalign{
y_1 \vee y_2 & := [1] \qquad \qquad\hbox { if } y_1=[0]= y_2, \cr
 & :=    [1] \prec   y_2 \qquad \hbox { if } y_1=[0]\not= y_2, \cr
 & := y_1\succ [1]  \qquad \hbox { if } y_1\not= [0]= y_2, \cr
 & := y_1\succ [1] \prec y_2 \quad \hbox { if } y_1\not=[0]\not= y_2. \cr
}$$
Therefore, by induction, it is clear that $\overline {K[Y_{\infty}]}$ is generated by $[1]$. \hfill$\square$
\M

\N {\bf 5.7. Proposition.} {\it  The unique dendriform algebra map $Dend(K) \to \oplus _{n\geq 1}
 {K[Y_n]}$ which sends
 the generator $x$ of $Dend(K)$ to $[1]$ is an isomorphism.}
\S

\N {\it Proof.} Let us show that the dendriform algebra $( \overline {K[Y_{\infty}]} ,\prec ,\succ)$ defined in 5.5
satisfies the universal condition to be the free dendriform algebra on one generator.

Let $D$ be a dendriform algebra and $a$ an element in $D$. Define a linear map  $ \aa : \overline {K[Y_{\infty}]} \to D $ by its value on the
trees $y = y_1 \vee y_2$ as follows :

$$\eqalign{
\aa(y_1 \vee y_2) & := a \qquad \qquad\hbox { if } y_1=[0]= y_2, \cr
 & :=    a \prec \aa (y_2) \qquad \hbox { if } y_1=[0]\not= y_2, \cr
 & := \aa (y_1) \succ a  \qquad \hbox { if } y_1\not= [0]= y_2, \cr
 & := \aa (y_1) \succ a \prec \aa (y_2) \quad \hbox { if } y_1\not=[0]\not= y_2. \cr
}$$
We claim that $\aa$ is a morphism of dendriform algebras. The proof is by induction on the degree of the tree. Indeed, on
one hand
$$\eqalign{
\aa (y\prec z) &= \aa (y_1 \vee (y_2 *z)) \cr
&= \aa (y_1)\succ a \prec \aa (y_2 *z) \cr
&= \aa (y_1)\succ a \prec (\aa (y_2) *\aa(z)). \cr
}$$
On the other hand,
$$\eqalign{
\aa (y)\prec \aa (z) &= (\aa (y_1) \succ a \prec  \aa(y_2)) \prec \aa(z) \cr
&= ((\aa (y_1) \succ a) \prec  \aa(y_2)) \prec \aa(z) \cr
&= (\aa (y_1) \succ a) \prec  (\aa(y_2) * \aa(z)) \cr
&= \aa (y_1) \succ a \prec (\aa(y_2) * \aa(z)). \cr
}$$
Here we supposed that  $y_1\not=[0]\not= y_2$, but the proof is similar for the other cases.

Since by lemma 5.6 $\overline {K[Y_{\infty}]}$ is generated by $[1]$, the morphism $\aa$ such that $\aa([1]) = a $ is unique.

It follows that  $( \overline {K[Y_{\infty}]} ,\prec ,\succ)$ is the free dendriform algebra on one generator.
\hfill $\square$

\M

\N {\bf 5.8. Theorem (Free dendriform algebra).} {\it  The unique dendriform algebra map 
$$Dend(V) \to 
 \oplus _{n\geq 1} {K[Y_n]}\otimes V^{\otimes n}$$
 which sends
 the generator $v\in V$   to $[1]\t v$ is an isomorphism.}
\S
\N {\it Proof.} Define the dendriform algebra structure on  
$\oplus _{n\geq 1} {K[Y_n]}\otimes V^{\otimes n}$ by 
$$\eqalign{
y\otimes \oo \prec y'\otimes \oo ' &:= (y\prec y')\otimes \oo \oo ',\cr
y\otimes \oo \succ y'\otimes \oo ' &:= (y\succ y')\otimes \oo \oo '.\cr}
$$
Since in the
relations defining a dendriform algebra the variables stay in the same order, the free
dendriform algebra over $V$ is completely determined by the free dendriform algebra on one
 generator:
$$
Dend(V) = \bigoplus _{n\geq 1}Dend(K)_n\t V^{\t n},
$$
where $Dend(K)_n$ is the subspace of $Dend(K)$ generated by all the possible products of $n$ copies
of the generator. Hence, by proposition 5.7, one gets $Dend(V) \cong \oplus_n K[Y_n]\t V^{\t
n}$.
\hfill $\square$
\M

\N {\bf 5.9. Remark.} The inverse isomorphism is obtained as follows. From $y\in Y_n$ we construct 
a monomial in the variables $x_1, \cdots, x_n$ by first putting the variable $x_i$ in
between the leaves $i-1$ and $i$. Then, for each vertex of depth one, we replace the
local patterns with two vertices by local patterns with one vertex:
$$
\arbreBxy \mapsto \arbreAxdy \qquad \qquad \qquad \arbreCxy \mapsto \arbreAxgy
$$
Continue the process until one reaches $\arbreAz$. The element $z\in Dend(V)$ is
the image of $(y;x_1\t \cdots \t x_n)\in K[Y_n]\t V^{\t n}$.  Observe that sometimes one needs to choose an
order to perform the process. But, thanks to the second axiom of dendriform algebras, the result does not
depend on this choice:
$$
\arbretroisxyz \mapsto (x\succ y)\prec z =  x\succ (y\prec z)\quad .
$$
Observe that the right $S_n$-module $Dend(n)$ which is such that $Dend(V) = \oplus _n Dend(n)\t _{S_n} V^{\t n}$,
is the regular representation $K[Y_n]\t K[S_n]$.
\M

\N {\bf 5.10.  Nested sub-trees and quotients.} Given a planar binary tree with $n+1$ 
leaves and a consecutive sequence of $k+1$ leaves $\{i, \ldots , i+k\}$, the sub-tree of
$y$ which contains these leaves is isomorphic to a unique planar binary tree
with $k+1$ leaves, which we denote by $y'$. We say that the sub-tree
$y'$ is  {\it nested} in $y$ at $i$. By definition the quotient $y'' = y/y'$ is the planar binary
tree with
$n-k+3$ leaves obtained from $y$ by removing the leaves $\{i+1, \ldots, i+k-1\}$.
\S

Example with $n=6, i=1, k=4$ and $y=[131612]$ :

$$\SUBTREEun \longrightarrow \SUBTREEdeux \longrightarrow \SUBTREEtrois
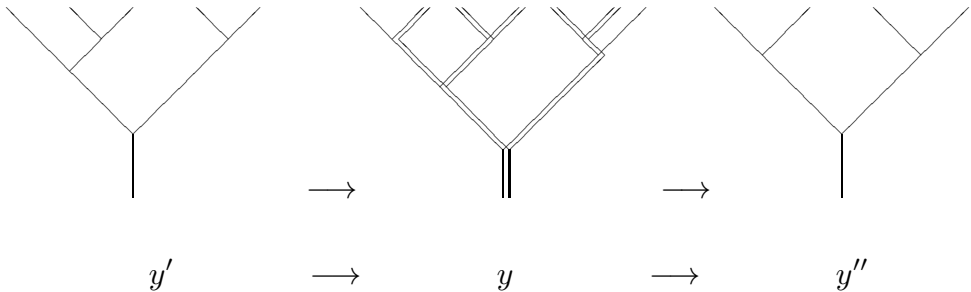
$$
$$y' \qquad \qquad \longrightarrow \qquad \qquad y \qquad \qquad \longrightarrow \qquad \qquad y''$$
\M

Observe that $y$ is a sub-tree of itself with quotient $[1] = \arbreA $, and that there are
$n$ nested sub-trees of $y$ of the form $[1]$ (with quotient $y$).

In terms of the permutation-like notation the names of $y'$ and of $y''$ are obtained as
follows. Let $y=[a_1 \ldots a_n]$. Start with the sequence of integers $[a_{i+1} \ldots  a_{i+k}]$ and
make it into a name of tree by first replacing the largest integer by $k$. Then proceed the same way
with the two remaining intervals, and so forth (like in A3). Ultimately one gets the name of $y'$.
In our example we get $[2141]$. The quotient tree $y''$ of $y$ by $y'$ is  obtained by
making the sequence of numbers $[a_1 \ldots a_i\  a_j\ a_{i+k+1} \ldots  a_n]$ into a
name of tree as before ($a_j$ is the largest integer in $[a_{i+1} \ldots  a_{i+k}]$). In our
example we get $[131]$. 
\M

\N {\bf 5.11. Proposition.} {\it Under the isomorphism of Proposition 5.7 the
composition operation in $Dend(K)$ induces the following composition operation in 
$\oplus _{n\geq 1} {K[Y_n]}$:
$$
y''\circ_i y' = \sum y,
$$
where the sum is extended over all the trees $y$ which contain $y'$ as a nested
sub-tree at $i$ and for which $y/y'= y'' $.}
\S

\N {\it Proof.} Since the operad is quadratic it suffices to check this assertion in low
dimension. The isomorphism gives  $[21]=x\prec x$ and  $[12] = x\succ x$. The eight
distinct cases of composition are: 
$$\eqalign{
[21]\circ_1[21] &= (x\prec x)\prec x = x\prec (x\prec x) + x\prec (x\succ x)=
      [213]+[312]\cr
[21]\circ_1[12] &= (x\succ x)\prec x =[131]\cr
[21]\circ_2[21] &= x\prec (x\prec x) =[321]\cr
[21]\circ_2[12] &= x\prec (x\succ x) =[312]\cr
[12]\circ_1[21] &= (x\prec x)\succ x =[213]\cr
[12]\circ_1[12] &= (x\succ x)\succ x =[123]\cr
[12]\circ_2[21] &= x\succ (x\prec x) =[131]\cr
[12]\circ_2[12] &= x\succ (x\succ x) = x\prec (x\succ x) + x\succ (x\succ x)=
       [213]+[312]\cr
}$$
In each case we verify that the trees of the right hand side are precisely such that
$y'$ is nested at $i$ with quotient $y''$. For instance both [213] and [312] have [21]
nested at 1.
\hfill $\square$.
\M

\N {\bf 5.12. Associative algebra structure on $K[Y_{\infty}]$ and $K[S_{\infty}]$.} By
lemma 5.2 the vector space $\overline {K[Y_{\infty}]}$ is an associative algebra for the
product  $$x*y = x\prec y + x \succ y,$$
 hence $K[Y_{\infty}]$ is a graded associative and unital algebra
whose structure is completely determined by the following two conditions: $\vert$ is a unit, 
the recursive formula
$$y*z = y_1 \vee (y_2 *z) + (y*z_1)\vee z_2,
$$
holds.

Observe that this algebra has an obvious
involution: $[i_1, \ldots , i_n] \mapsto [i_n, \ldots , i_1]$ on trees.  In theorem 3.8 of [LR]
we prove that it is isomorphic to the tensor algebra over
$K[Y_{0,\infty}]$, where 
$Y_{0,n-1}:= \lbrace\ y = |\vee y\rq \in Y_n,\ {\rm for}\ y\rq \in 
Y_{n-1}\rbrace $. It is also isomorphic to $T(K[Y_{\infty, 0}])$, where 
$Y_{n-1,0}:= \lbrace\ y = y\rq \vee \vert\rbrace $.

The map $\psi_n: S_n \to Y_n$ (cf. Appendix A) induces a linear map  $\psi: K[S_{\infty}] \to
K[Y_{\infty}]$. There is an associative and unital algebra structure on $ K[S_{\infty}]$ given
by 
$$
x*y:= sh_{n,m}\cdot (x\times y) \in K[S_{n+m}], $$
where  $x\in S_n$, $y\in
S_m$ and  $sh_{n,m} \in  K[S_{n+m}]$ is  the sum of all the $(n,m)$-shuffles. It is proved in
[LR] that $\psi$ is an associative algebra homomorphism.
\M

\N {\bf 5.13. Hopf structure on $Dend(V)$.} It is well-known that the free
associative algebra $T(V)$ is in fact a cocommutative Hopf algebra, where the coproduct
is given by the shuffle.  Similarly there exists a structure of Hopf algebra on the associative unital
algebra $K\oplus Dend(V)  = \oplus_{n\ge 0} K[Y_n]\t V^{\t n}$. The coproduct is completely determined
by the shuffles and the coproduct on $ \oplus_{n\ge 0} K[Y_n]$ was constructed in [LR].
Observe that this coproduct is not cocommutative. It is related to the ``brace algebras" and has been studied in
details in [R].
\vfill \eject



\noindent {\bf 6. (CO)HOMOLOGY OF DENDRIFORM  ALGEBRAS}
\medskip

In this section we show that there exists a chain complex (of Hochschild type), for any
dendriform  algebra. It enables us to construct a  homology and a cohomology theory for
dendriform  algebras. It will be proved in section 8 that
these theories are the ones predicted by the operad theory in characteristic zero. 
\medskip

\N {\bf 6.1. The chain complex of a dendriform  algebra.} Let $E$ be a dendriform  algebra and let $C_n
$ be the set $\{1,\cdots ,n\}$. We define the module of $n$-chains of $E$
as
$$
C^{Dend}_n(E):= K[C_n]\otimes E^{\otimes n},
$$
and the differential $d =- \sum_{i=1}^{n-1} (-1)^i d_i: C^{Dend}_n(E) \to C^{Dend}_{n-1}(E)$ as
follows. First, we define the face operators $d_i, 1\leq i\leq n-1$, on $r\in C_n$ by
$$
d_i(r) = \cases{
r-1  & if $i\leq r-1,$\cr
r   & if $ i\geq r.$\cr} 
$$
These maps are extended linearly to maps
$$
d_i: K[C_n] \to K[C_{n-1}],\ 1\leq i\leq n-1.
$$
Second, we define the symbol $\circ _i^r$ as follows:
$$
\circ _i^r = \cases{
*   & if $ i<r-1,$\cr
\succ    & if $i=r-1,$\cr
\prec    & if $ i=r,$\cr
*    & if $ i>r.$\cr}
$$
Recall that $x*y= x\prec y + x\succ y$.

Finally the map $d_i: C^{Dend}_n(E) \to C^{Dend}_{n-1}(E)$ is given by
$$
d_i(r; x_1\otimes \cdots x_n):= (d_i(r);
 x_1\otimes \cdots \otimes x_{i-1}\otimes x_{i}\circ _i^r x_{i+1}\otimes \cdots \otimes x_n)
$$
for $1\leq i\leq n-1$.
\M

\N {\bf 6.2. Lemma.} {\it The maps $d_i: C^{Dend}_n(E) \to C^{Dend}_{n-1}(E)$ satisfy the simplicial
relations
 $d_id_j = d_{j-1}d_i$, for $i<j$, and so $(C^{Dend}_*(E),d)$ is a chain complex.}
\S
\N {\it Proof.} Let us first prove the lowest dimensional case, that is $d_1d_1 = d_1d_2$
on $C^{Dend}_3(E)= 3\ E^{\otimes 3}$. 

On the first component $(r=1)$ one gets:
$$\eqalign{
d_1d_1(1; a\otimes b \otimes c) &= d_1(1; (a\prec b) \otimes c) = (1; (a\prec b) \prec
c),\cr
d_1d_2(1; a\otimes b \otimes c) &= d_1(1; a\otimes (b*c)) = (1; a\prec (b*c)),\cr}
$$
hence $d_1d_1 = d_1d_2$ by axiom (i).

On the second component $(r=2)$ one gets:
$$\eqalign{
d_1d_1(2; a\otimes b \otimes c) &= d_1(1; (a\succ b) \otimes c) = (1; (a\succ b) \prec
c),\cr
d_1d_2(2; a\otimes b \otimes c) &= d_1(2; a\otimes (b\prec c)) = (1; a\succ (b\prec
c)),\cr} $$
hence $d_1d_1 = d_1d_2$ by axiom (ii).

On the third component $(r=3)$ one gets:
$$\eqalign{
d_1d_1(3; a\otimes b \otimes c) &= d_1(2; (a* b) \otimes c) = (1; (a*b) \succ
c),\cr
d_1d_2(3; a\otimes b \otimes c) &= d_1(2; a\otimes (b\succ c)) = (1; a\succ
(b\succ c)),\cr} $$
hence $d_1d_1 = d_1d_2$ by axiom (iii).

Higher up, the verification of $d_id_j = d_{j-1}d_i$ splits up into two different cases. First,
if $j=i+1$, then it is the same kind of computation as above, so it is a consequence of the
axioms of a dendriform  algebra. Second, if $j>i+1$, then both operations agree on $C_n$ (direct
checking) and the image of $(a_1\otimes \cdots \otimes a_n)$ is, in both cases, 
$$
(a_1\otimes \cdots \otimes a_i \circ _i^r a_{i+1}\otimes \cdots \otimes a_j \circ _j^r
a_{j+1}\otimes
\cdots \otimes a_n).
$$
\hfill $\square$
\M

\N {\bf 6.3. The chain bicomplex of a dendriform  algebra.}
Observe that $C^{Dend}_*(E)$ is in fact the total complex of a bicomplex. Indeed, let
$$
C^{Dend}_{p,q}(E) = \{p\} E^{\t p+q}\quad {\rm for }\quad p\geq 1, q\geq 0.
$$
The map 
$$
d^h: C^{Dend}_{p,q}(E) \to  C^{Dend}_{p-1,q}(E)\quad {\rm  by}\quad  d^h:= \sum_{i=1}^{p-1} (-1)^id_i ,
$$
is well-defined because $d_i(p) = p-1$, when $i\leq p-1$, and
$$
d^v: C^{Dend}_{p,q}(E) \to  C^{Dend}_{p,q-1}(E)\quad {\rm  by}\quad  d^v:= \sum_{i=p}^{p+q} (-1)^id_i .
$$
is well-defined because $d_i(p) = p$, when $i\geq p$.
In other words, the $p$-th component of $C^{Dend}_n(E)$ is put in bidegree $(p,n-p)$.
\M

\N {\bf 6.4. (Co)homology of dendriform  algebras.} By definition the homology (with trivial
coefficients) of a dendriform  algebra $E$ is
$$
H^{Dend}_*(E):= H_*( C^{Dend}_*(E), d)\ ,
$$
and the cohomology of a dendriform  algebra $E$ (with trivial coefficients) is
$$
H_{Dend}^*(E):= H^*(\Hom (C^{Dend}_*(E), K)).
$$
Let us use freely the interpretation of the preceding results in terms of operads as devised
in the next section.
From Koszul duality and Appendix B5d, the graded module $H_{Dend}^*(E)$ is naturally
equipped with a structure of graded dialgebra (and hence a structure of graded Leibniz
algebra).
\M
\N {\bf 6.5. Theorem.} {\it The dendriform  algebra homology of a free dendriform algebra is
trivial. More precisely
$$\eqalign{
H^{Dend}_n(Dend(V)) &= 0 \ for \ n>0,\cr
H^{Dend}_1(Dend(V)) &=V .\cr}
$$
}
\S
\N {\it Proof.} Since we know already that the analogous theorem for associative
dialgebras is true (Theorem 3.8), it is a consequence (by the operad theory, cf.
Appendix B), of the operad duality between $Dias$ and $Dend$ proved in Proposition 8.3.  
\hfill $\square$

\vfill \eject



\noindent {\bf 7. ZINBIEL ALGEBRAS, DENDRIFORM ALGEBRAS AND HOMOLOGY}
\medskip

In this section we introduce Zinbiel (i.e. dual-Leibniz) algebras  and we compare them
with dendriform  algebras. In particular we compute the natural map from a free dendriform  algebra to the
free Zinbiel algebra, considered as a dendriform  algebra. Finally we compare the homology
theories.
 \medskip

\noindent {\bf 7.1. Zinbiel algebras} [L3].
A {\it  Zinbiel algebra}  $R$  
\footnote {(*)}{Terminology proposed by J.-M. Lemaire}
(also called {\it dual-Leibniz algebra}) is a module over $K$ equipped with a binary operation
$(x,y)\mapsto x\cdot y$, which satisfies the identity
$$
(x\cdot y)\cdot z=x\cdot (y\cdot z)+x\cdot (z\cdot y)\ ,\ \ for\ all\  x,y,z\in R.\leqno (7.1.1)
$$
The category of Zinbiel algebras is denoted ${\bf ZinbÊ}$. The free Zinbiel algebra over the
vector space $V$ is ${\overline T}(V) = \oplus_{n\geq 1}V^{\t n}$ equipped with the following
product
$$
(x_0 \ldots x_p) \cdot (x_{p+1}\ldots x_{p+q}) = x_0 sh_{p,q} (x_1 \ldots x_{p+q})
$$
where $sh_{p,q}$ is the sum over all $(p,q)$-shuffles. We denote it by $Zinb(V)$. Observe
that $Zinb(V)= \bigoplus_n Zinb(n)\t_{S_n} V^{\t n}$ for $Zinb(n)= K[S_n]$. It is immediate
to check that the symmetrized product 
$$xy:=x\cdot y+y\cdot x\leqno (7.1.2)$$
 is associative
 (cf. [R], [L3]), so,  under the symmetrized product, $R$ becomes an
associative and commutative algebra. This construction gives a functor
$$
{\bf ZinbÊ}\buildrel + \over \longrightarrow {\bf Com}.
$$
Let us construct functors to and from the category of 
dendriform  algebras ${\bf Dend}$.
\M

\N {\bf 7.2. Lemma.} {\it  Let $R$ be a Zinbiel algebra and put
$$ x\prec y:=x\cdot y\ ,\ x\succ y:=y\cdot x\ ,\ \forall\  x,y\in R.
$$ 
Then $(R,\prec,\succ)$ is a dendriform  algebra denoted $R_{Dend}$. Conversely, a commutative
dendriform  algebra (i.e. a dendriform  algebra for which $x\succ y = y \prec x$) is a Zinbiel algebra. }
\S
\N {\it Proof.} Indeed, relation (i) is exactly relation
(3.3.1) and so is relation (iii). Relation (ii) also follows from (3.3.1) since
$$ x\succ (y\prec z)=(yz)x,\quad (x\succ y)\prec z=(yx)z,
$$ 
and the relation $(y\cdot z)\cdot x=(y\cdot x)\cdot z$ follows from (3.3.1) for $y,z,x$ and
$y,x,z$.  \hfill $\square$
\M

So we have constructed a functor
$$ 
{\bf ZinbÊ}\to {\bf Dend} .
$$

From Lemma 5.2 we get a functor
$$ 
{\bf Dend} \buildrel + \over \longrightarrow {\bf As}.
$$
Summarizing we get the following
\M

\noindent {\bf 7.4. Proposition.}
{\it The following  diagram of functors between categories of algebras is commutative
$$
\matrix{
{\bf ZinbÊ}&\hookrightarrow &{\bf Dend}\cr
&&\cr
\downarrow\scriptstyle{+} &&\downarrow\scriptstyle{+} \cr
&&\cr
{\bf Com}&\hookrightarrow &{\bf As}\cr
}
$$}
\S
\N{\it Proof}. 
The commutativity of the diagram is immediate since for a Zinbiel
algebra the associated associative products are equal :
$$
x*y=x\prec y+x\succ y=x\cdot y+y\cdot x = xy \quad .
$$
\hfill $\square$
\medskip

\noindent {\bf 7.5. Theorem.} {\it For any vector space $V$, the natural map of dendriform  algebras
$Dend(V) \to Zinb(V)_{Dend}$ is induced, in degree $n$, by the map
$$
K[Y_n]\t K[S_n] \to K[S_n] , \quad y\t \oo \mapsto \sum_{\{\ss\in S_n\  \mid\  \psi '(\ss)
= y\}} \ss \oo\ ,
$$
where $\psi ': S_n \proj Y_n$ is the surjective map described in Appendix A6.}
\B

\noindent {\it Proof.} First, observe that the map $V\to Zinb(V)$ gives a map (the same on
the  underlying vector spaces) $V\to Zinb(V)_{Dend}$. Since this latter object is a
dendriform  algebra, the map factors through the free dendriform  algebra on $V$, whence a natural
dialgebra map $\phi : Dend(V) \to  Zinb(V)_{Dend}$.

The proof will be done by induction on
$n$.

Restricting $\phi$ to   the degree $n$ part gives a commutative square (cf. 5.8 and 5.9):
$$
\matrix{
K[Y_n]\otimes V^{\otimes n} &\buildrel {\phi_n}\over \longrightarrow &V^{\otimes n} \cr
{}\cr
\bigcap &&\bigcap \cr
{}\cr
Dend(V) &\buildrel {\phi}\over \longrightarrow &Zinb(V)_{Dend}\ .\cr}
$$

For $n=1$, $\phi_1$ is clearly the identity of $V$.
\S
For $n=2$, $\phi_2 : K[Y_2]\t V^{\t 2} \to V^{\t 2}$ is given by
$$\displaylines{
\arbreC \t (x,y) = x\prec y \mapsto x\cdot y = [12].(x,y), \cr
\arbreB \t (x,y) = x\succ y \mapsto y\cdot x = [21]. (x,y), \cr
}$$
so $\arbreC \mapsto [12]$ and $\arbreB \mapsto [21]$, which is the map $\psi '$ of Appendix A6.
\S

For n=3, one has 
$$\eqalign{
\arbrecinq\t (x,y,z) = x\prec (y\prec z) &\mapsto x\cdot (y\cdot z) = [123].
(x,y,z),\cr
\arbrequatre\t (x,y,z) = x\prec (y\succ z) &\mapsto x\cdot (z\cdot y) = [132].
(x,y,z),\cr
\arbretrois\t (x,y,z) = (x\succ y)\prec z &\mapsto (y\cdot x)\cdot z = 
y\cdot (x\cdot z+z\cdot x)\cr
&\qquad =  ([213]+ [312]).
(x,y,z),\cr
\arbredeux\t (x,y,z) = (z\prec x)\succ y &\mapsto x\cdot (y\cdot z) = [231].
(x,y,z),\cr
\arbreun\t (x,y,z) = (z\succ y)\succ x &\mapsto x\cdot (y\cdot z) = [321].
(x,y,z),\cr
}$$
so the map is precisely $\psi '$.
\S
Let us now prove it for any $n$. We suppose that the theorem has been proved for any $p<n$. We
are going to prove it for $n$. We use the {\it grafting} operation on trees (cf. Appendix A).

The element $(y_1\vee y_2; x_1 \ldots x_{p+q+1})$ can be written as  
$$
(y_1; x_1 \ldots
x_p)\succ x_{p+1}\prec (y_2; x_{p+2} \ldots x_{p+q+1})
$$
 in $Dend(V)$ when $p\geq 1, q\geq 1$. Here $x_{p+1}$ stands for $([1];x_{p+1})$. We
do not need to put any parenthesis because of the second relation of dendriform algebras. The
image of this element under $\phi$ is
$$
\matrix{
 (x_{p+1}\cdot \phi_q(y_2; x_{p+2} \ldots  x_{p+q+1}))\cdot 
\phi_p(y_1; x_1 \ldots  x_{p+q+1})\qquad\hfill \cr
\hfill= x_{p+1}\cdot (\phi_q(y_2; x_{p+2} \ldots  x_{p+q+1})\cdot 
\phi_p(y_1; x_1 \ldots  x_{p+q+1}) \cr
\hfill + \phi_p(y_1; x_1 \ldots 
x_{p+q+1})\cdot \phi_q(y_2; x_{p+2} \ldots  x_{p+q+1}))
.\cr}
$$
By induction we know that $\phi_p (y_1;x_1 \ldots  x_p) = \sum_{\sigma_1 \in \psi '^{-1}
(y_1)} \sigma_1 (x_1 \ldots  x_p)$ and that $\phi_q (y_2;x_{p+2} \ldots  x_{p+q+1}) =
\sum_{\sigma_2 \in \psi '^{-1} (y_2)} \sigma_2 (x_{p+2} \ldots  x_{p+q+1})$. For  $z_i\in V$
one has, in the free Zinbiel algebra, the equality 
$$z_0\cdot(z_1\ldots z_p
\cdot z_{p+1}\ldots  z_{p+q} + z_{p+1}\ldots 
z_{p+q}\cdot z_1\ldots z_p) = \sum_{\tau} \tau(z_0 \ldots  z_{p+q}),$$
 where the sum is extended
over all the $(p,q)$-shuffles $\tau$ (acting on the set $\{1, \ldots , p+q\}$. Hence we get 
$$
\phi_n({y_1\vee y_2};x_1\ldots x_n) = \sum_{\pi} \pi (x_1\ldots x_n),
$$
where the permutation $\pi$ is of the following form :
$$\matrix{
\pi (i) = \tau \sigma_1(i) {\rm \ for\ } 1\leq i\leq p,\hfill\cr
\pi (p+1) = 1,\hfill\cr
\pi (i) = \tau (p+ \sigma_2(i))  {\rm \ for\ } p+2\leq i\leq p+q+1, \hfill\cr
}
$$
where $\sigma_1 \in \psi '^{-1}
(y_1),\ \sigma_2 \in \psi '^{-1} (y_2)$
 and $\tau$ is a $(p,q)$-shuffle.

On the other hand one easily checks that all the permutations $\sigma$ which belong to
$\psi '^{-1}(y_1\vee y_2)$ are precisely obtained by choosing such a $\sigma_1$ and
such a $\sigma_2$, and then shuffle the associated levels.
 So we have proved the formula for $n$.

In the preceding proof we assumed that $p\geq 1, q\geq 1$. We let to the reader the task
of modifying the proof when either $p=0$ or $q=0$.
\hfill $\square$
\M

\M
\N {\bf 7.6. Comparison of homology theories.} As mentioned in Appendix B, for any algebra
over a Koszul operad, there is a small chain complex, modelled on the dual operad, whose homology
is the homology of the algebra. For Zinbiel algebras it takes the following form (cf. [Li2]):
$$\cdots\longrightarrow R  ^{\otimes n}\buildrel {d}\over
\longrightarrow R  ^{\otimes
n-1}\buildrel {d}\over
\longrightarrow\cdots\buildrel {d}\over
\longrightarrow R  ^{\otimes 2}\buildrel {\cdot}\over
\longrightarrow R  \leqno C_*^{Zinb}(R  ):$$
where
$$
d(x_1,\ldots ,x_n)=(x_1\cdot x_2,x_3,\ldots ,x_n) + \sum_{i=2}^{n-1}(-1)^{i-1}(x_1,\ldots
,x_i * x_{i+1}, \ldots ,x_n).
$$
The homology groups of this complex are denoted $H^{Zinb}_n(R)$, for $n\geq 1$. By Appendix
B5e, there is a natural map of complexes
$$
C^{Dend}_* (R_{Dend}) \longrightarrow C_*^{Zinb}(R)
$$
inducing
$$
H^{Dend}_* (R_{Dend}) \longrightarrow H_*^{Zinb}(R).
$$
Let us describe explicitly the chain complex map.
\M
\N {\bf 7.7. Proposition.} {\it Let  $\theta_n^r \in K[S_n]$ be the elements defined
recursively by the formulas :
$$\displaylines{
  \theta_1^1(1) = (1), \qquad \theta_n^r(1, \ldots, n) =(r, \tilde
{\theta}_n^r(1,\ldots, \widehat r \ldots, n), \cr
\tilde {\theta}_n^r(1,\ldots , n-1) = (r-1, \tilde {\theta}_{n-1}^{r-1}(1,\ldots, \widehat {r-1}
\ldots, n) +(r, \tilde {\theta}_{n-1}^r(1,\ldots, \widehat {r}
\ldots, n)
.\cr
}$$
In particular, $\theta_n^1(1, \ldots, n)= (1, \ldots, n)$ and $\theta_n^n(1, \ldots, n) = 
 (n, \ldots, 1)$.
The map $\Theta_n : K[C_n]\t R^{\t n} \to R^{\t n}$ defined by
$$
\Theta_n ([r]\t (x_1, \cdots x_n) ) := \theta_n^r (x_1, \cdots x_n) 
$$
 is the chain complex map 
$\Theta_* : C^{Dend}_* (R_{Dend}) \longrightarrow C_*^{Zinb}(R) $ induced by the operad
morphism $Dend \to Zinb$.}
\M
\N {\it Proof.} This is a consequence of the explicit description of the morphism of
Leibniz algebras $Leib(V) \to Dias(V)_{Leib}$, cf. Appendix B5e.

Observe that $\theta_n^r$ is the sum of the signed action of $({n \atop r})$ permutations. Since
$\sum ({n \atop r}) = 2^n$, we recover (up to permutation) the $2^n$ monomials of $[x_1, [x_2,
[\ldots, x_n]]]$ (cf. 4.10). \hfill $\square$

 \vfill \eject



\noindent {\bf 8. KOSZUL DUALITY FOR THE DIALGEBRA OPERAD}
\medskip

In this section we show that the operad associated to dendriform  algebras is dual, in the
operadic sense, to the operad associated to dialgebras.
Moreover we use the results of section 4 to show that the operad of associative dialgebras is a
Koszul operad (and so is the operad of dendriform  algebras). The reader not familiar with the
notions of operad and Koszul duality may have a look at Appendix B, from which we take the
notation.
\M

\N {\bf 8.1. The associative  dialgebra operad.} A dialgebra is determined by two operations (left
and right product) on two variables and by relations which make use of the composition of two
such operations. Hence the operad
$Dias$  associated to the notion of dialgebra is a  binary (operations on two variables) quadratic
(relations involving two operations)  operad. Moreover, there is no symmetry property for these
operations, and, in the relations, the variables stay in the same order. Hence the operad is a {\it
non-$\SS$-operad}, that is, as a representation of
$S_n$,  the space
$Dias(n)$ is a sum of copies of the regular representation.
It was proved in 2.5 that the free dialgebra over the vector space $V$ is $Dias(V)=
T(V)\otimes V\otimes T(V)$. The degree $n$ part of it is $\bigoplus_{i+1+j=n} V^{\otimes i}\otimes
V\otimes V^{\otimes j}$. Hence the operad
$Dias$ is such that 
$$Dias(n) = n K[S_n]\qquad \hbox {($n$ copies of the regular representation)}.$$
In particular $Dias(1) = K$ (the only unary operation on a dialgebra is the identity), $E:=Dias(2)=
E'\t K[S_2]$, where $E'$ is 2-dimensional  generated by $\g$ and $\d$. The space $\Ind _{S_2}^{S_3}(E\t E)$
is the sum of 8 copies of the regular representation of $S_3$. Each copy corresponds to a choice of
parenthesizing:  $(-\circ_1 ( - \circ_2 -))$ or $((-\circ_1 -)\circ_2 -)$, and a choice for
the two operations $\circ_1$ and $\circ_2$. The space of relations $R\subset \Ind _{S_2}^{S_3}(E\t E)$
is of the form $R'\t K[S_3]$, where $R'$ is the subspace of $E'^{\t 2}\oplus E'^{\t 2}$ determined by
the relations 1 to 5 of a
 dialgebra (cf. 2.1).
\M

\N {\bf 8.2. The dendriform  algebra operad.} Analogously the operad $Dend$ associated to the notion of
dendriform  algebra is binary and quadratic, and is a {\it non-$\SS$-operad} since there is no symmetry
for the operations and, in the relations, the variables stay in the same order. It was proved in 5.7
that the free dendriform  algebra over the vector space
$V$ is $Dend(V)=
\bigoplus_{n\geq 1} K[Y_n]\otimes V^{\otimes n}$. Hence the operad
$Dend$ is such that 
$$Dend(n) = K[Y_n]\t K[S_n].$$
In particular $Dend(1) = K, Dend(2)=
F'\t K[S_2]$, where $F'$ is 2-dimensional generated by $\prec$ and $\succ$. The space 
$$\Ind
_{S_2}^{S_3}(Dend(2)\t Dend(2))\cong (F'^{\t 2} \oplus F'^{\t 2})\t K[S_3]$$
is the sum of 8 copies of the regular representation of
$S_3$. The space of relations is of the form
$S'\t K[S_3]$, where $S'\subset (F'^{\t 2} \oplus F'^{\t 2})$ is the subspace determined by the 3
relations of a dendriform  algebra (cf. 5.1).
\medskip

\N {\bf 8.3. Proposition.} {\it The operad $Dend$ of dendriform  algebras is dual, in the operad
sense, to the operad $Dias$ of dialgebras : $Dias^! = Dend$.}
\S
\N {\it Proof.} Let us identify $F'=E'^{\vee}$ with $E'$ by identifying the basis $(\prec,
\succ)$ with the basis $(\g, \d)$.  Since $R=R'\t K[S_3]$, the space of relations for the dual
operad is of the form
$R'^{ann}\t K[S_3]$, where, according to Proposition B3, $R'^{ann}$ is the annihilator of
$R'$. Recall from Proposition B3 that the scalar product on $E'^{\t 2}\oplus E'^{\t 2}$ is given by
the matrix
$\left [
\matrix{\Id & 0\cr 0&-\Id\cr}\right ]$.

With obvious notation, the subspace $R'$ is determined by the relations
$$\left\{ 
\matrix {
(\g,\g)_2 - (\g, \g)_1 = 0,\cr
(\g,\g)_1 - (\g, \d)_2 = 0,\cr
(\d,\g)_1 - (\d, \g)_2 = 0,\cr
(\g,\d)_1 - (\d, \d)_2 = 0,\cr
(\d,\d)_2 - (\d, \d)_1 = 0.\cr}\right.
$$
It is immediate to verify that its annihilator $R'^{\perp}$ with respect to the given scalar product
is the subspace determined by the relations
$$\left\{ 
\matrix {
(\g,\g)_1 - (\g, \g)_2 - (\g, \d)_2 = 0,\cr
(\d,\g)_1 - (\d, \g)_2 = 0,\cr
(\d,\d)_2 - (\g, \d)_1 - (\d, \d)_1 = 0.\cr}\right.
$$
This is precisely the relations of dendriform  algebras (once we have changed $\g$ into $\prec$ and $\d$
into $\succ$). \hfill $\square$
\M

\N {\bf 8.4. Proposition.} {\it The chain complex $CY_*(D)$ associated to a dialgebra $D$ is the
chain complex of $D$ in the operad sense (cf. Appendix B4). Hence $HY$ is the (co)homology theory for dialgebras
predicted by the operad theory.}
\S

\N {\it Proof.}  From the theory of operads recalled in Appendix B we have
$$
C_n^{Dias}(D) = Dend(n)^{\vee} \otimes_{S_n} D^{\otimes n}.
$$ 
By Proposition 5.7 we get
$$
C_n^{Dias}(D) = K[Y_n]\otimes K[S_n]\otimes _{S_n} D^{\otimes n} =  K[Y_n]\otimes D^{\otimes n} =
CY_n(D).
$$
So, it suffices to prove that the boundary operator $d$ of $CY_*(D)$
agrees with the dialgebra structure of $D$ on $CY_2(D)$ and that it is a coalgebra
 derivation.

The boundary map on $CY_2(D) = D^{\t 2}\oplus D^{\t 2}$ is given by $x\t y \mapsto x\g y$ on the first
 component and by  $x\t y \mapsto x\d y$ on the second one. So the first condition is
fulfilled. Checking the coderivation property is analogous to the associative case (cf. B4).
\hfill $\square$
\M

\N {\bf 8.5. Theorem.} {\it The operad $Dias$ of dialgebras is a Koszul operad, and so is the operad
 $ Dend$ of dendriform algebras.}
\S
\N {\it Proof.} From the definition of a Koszul operad given in B4, this theorem follows from
the vanishing of $H_*^{Dias} = HY_*$ of a free dialgebra, as proved in Theorem 3.8. Since $\P$ being
Koszul implies $\P^!$ is Koszul (cf. Appendix B5b), the operad $Dend$ is Koszul.
\hfill
$\square$
\M

Observe that this theorem, together with the general property of Koszul operads recalled in B4a,
implies proposition 5.2.
\medskip

\N {\bf 8.6. Poincar\'e series.} From the description of the operad $Dias$ it follows that its
Poincar\'e series  is
$$
g_{Dias}(x) = \sum_{n\geq 1} (-1)^n n n! {x^n \over n!} =  \sum_{n\geq 1} (-1)^n n x^n = {-x \over
(1+x)^2}.
$$
On the other hand, the Poincar\'e series of $ Dend$ is
$$
g_{Dend}(x) = \sum_{n\geq 1} (-1)^n c_n n! {x^n \over n!} =  \sum_{n\geq 1} (-1)^n c_n x^n =
{-1-2x+\sqrt{1+4x}\over 2x}.
$$
As expected (cf. Appendix B5c) we verify that $g_{Dend}(g_{Dias}(x)) = x$.
\M

\N {\bf 8.7. Operad morphisms.} Any morphism of operads induces a functor between the associated
categories of algebras. Taking the dual gives also a morphism, but in the other direction.  For
instance the dual of the inclusion of the category of commutative algebras into the category of
associative algebras is the ``-" functor which transforms an associative
algebra into a Lie algebra. All the functors between categories of algebras that we met in the
previous sections come from morphisms of operads. They assemble into the following commutative
diagram of functors (cf. Proposition 4.4 and Proposition 7.4):
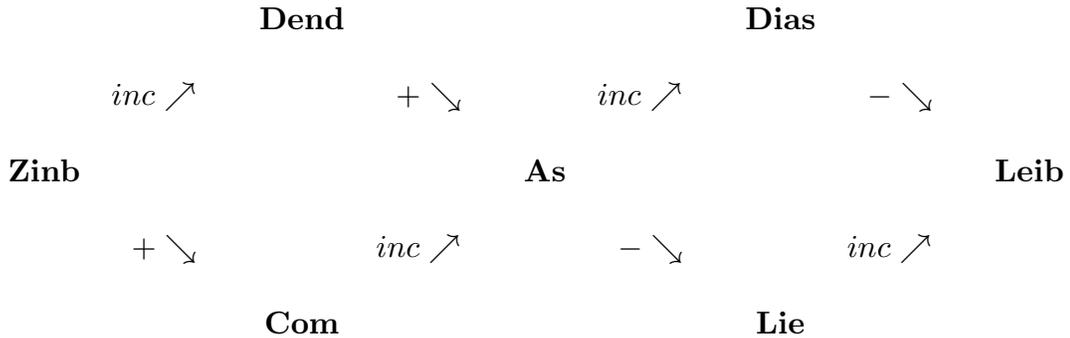
\B

\def\mydiagram{
\matrix{
&&{\bf Dend}&&&&{\bf Dias}&\cr
&&&&&&&\cr
&inc \nearrow\quad&&+ \searrow \ &&inc \nearrow\quad&&-\searrow \ \cr
&&&&&&&\cr
{\bf Zinb}&&&&{\bf As}&&&&{\bf Leib}\cr
&&&&&&&\cr
&+ \searrow \ &&inc \nearrow\quad&&- \searrow \ &&inc \nearrow\quad\cr
&&&&&&&\cr
&&{\bf Com}&&&&{\bf Lie}&\cr
}}
\B
 $\mydiagram$
\B

The symmetry (around a vertical axis passing through ${\bf As^! =As}$) reflects the Koszul duality of
quadratic operads : $Lie^! = Com, Dias^! = Dend, Leib^! = Zinb$.

\vfill \eject



\noindent {\bf 9. STRONG HOMOTOPY ASSOCIATIVE DIALGEBRAS}
\medskip
Strong homotopy $\P$-algebras are governed by the Koszul dual operad $\P ^!$ (cf. Appendix B6). Since
 for $Dias$ we know of an explicit description of its dual $Dias ^!$, we are able to describe
explicitly the notion of {\it strong homotopy dialgebra}. In order to do it we use the
notion of nested sub-tree of a planar binary tree.
 \M

\N {\bf 9.1. Nested sub-trees.} Recall that in section 5 we defined the notion of nested
sub-tree $y'$ of a
 tree $y$ with quotient $y''=y/y'$ and we described composition in the free dendriform  algebra in terms of
nested sub-trees, cf. 5.10 and Proposition 5.11.
\M

\N {\bf 9.2. Theorem.} {\it A strong homotopy dialgebra is a graded vector space $A = \oplus
_{i\in \ZZ}A_i$ equipped with operations 
$$ m_y : A^{\t n} \to A , \ for\  any \ y\in Y_n ,\  n\geq 1,$$
which are homogeneous of degree $n-2$ and which satisfy the following relations for any $y$:
$$\sum _{y'\subset y,\ y''=y/y'} \pm m_{y''}(a_1, \ldots ,a_i, m_{y'}(a_{i+1}, \ldots,
a_{i+k}), \ldots, a_n) =0,\leqno (*)_y
$$
where the sign $\pm$ is $+$ or $-$ according to the parity of $(k+1)(i+1) + k(n+ \sum_{j=1}^{k} \vert
a_j\vert)$.

 In this relation the tree $y$ is fixed and the sum runs over all the nested sub-trees $y'$
of
$y$ with $y''=y/y'$ (as described in 5.10).}
\S
Compare with the definition of a $A_{\infty}$-algebra  [St, p. 294].
\S

\N {\it Proof.} Since the operad $Dias$ is Koszul (cf. Theorem 8.5), we may apply Theorem
B.7. By Proposition 5.8 the dual operad $Dend$ is generated in dimension $n$ (as a free
$S_n$-module) by the set of $n$-trees $Y_n$. Hence the operations on a  ${\cal B} (\P ^!)^*$-algebra
$A$  are generated by operations $m_y$, for $y\in Y_n$,
$$ m_y : A^{\t \vert y\vert } \to A\quad  \hbox  { of degree }\vert y\vert - 2.$$
The relations satisfied by these operations are obtained as follows. First one extends them
in order to get a coderivation
$$
m : \oplus_n K[Y_n]\t A^{\t n} \longrightarrow  \oplus_n K[Y_n]\t A^{\t n},
$$
and then one writes $m \circ m = 0$. The component in $K[Y_1]\t A = A$ of the image under $m$
of an element $K\{y\}\t A^{\t n} \subset  K[Y_n]\t A^{\t n}$ is given by $m_y$ (up to sign). More
precisely we put
$$
m_y(a_1, \cdots , a_n) := (-1)^{(k-1)\vert a_1\vert + (k-2)\vert a_2\vert + \cdots a_{k-1}} m(y;
a_1, \cdots , a_n)_1.
$$
 In order to obtain the component in $K\{y'\}\t
A^{\t k}$ we need to look at the composition in the cofree co-dendriform  algebra, or dually, in the free
dendriform algebra:
$$
Dias ^!(k)\t Dias ^!(i_1)\t \cdots \t Dias ^!(i_k) \longrightarrow Dias ^!(i_1+\cdots +i_k)$$
or equivalently,
$$
K[Y_k]\t K[Y_{i_1}]\t \cdots \t K[Y_{i_k}]\longrightarrow K[Y_{i_1+\cdots +i_k}].
$$
It is sufficient to write the relation $m \circ m = 0$ for the component in $K[Y_1]\t A = A$ of
the image,
since the vanishing of the other components is a consequence of that one.  Hence it is sufficient to
compute the composition product for $i_u=1$ for all $i_u$ except one of them, let us say $i_j=m$. In
this case it is precisely the result of Proposition 5.11. \hfill $\square$
\M

\N {\bf 9.3. Strong homotopy associative dialgebra in low dimensions.} Let us write $m_y = m_{i
j
\cdots k}$  in place of $m_{[i j \cdots k]}$, when $y= [i j \cdots k]$.

For $n=1$ the operation $\dd :=m_1 : A\to A$
is of degree $-1$. Since the only nested sub-tree of $[1]$ is $[1]$ itself, the relation
$(*)_1$ is  $$
\dd \circ \dd = 0.
$$
Hence $(A, \dd)$ is a chain complex.
\M
For $n=2$, there are two maps $m_{12}$ and $m_{21} : A^{\t 2} \to A$, which are of degree 0.
The tree $[12]$ (resp. $[21]$) has three nested sub-trees, itself  with quotient $[1]$,
and two (different) copies of $[1]$, both with quotient $[12]$ (resp. $[21]$). Hence the
relation $(*)_{12}$ takes the form (where 1 stands for $\Id$):
$$
\dd\circ  m_{12} - m_{12} \circ (\dd\t 1)  - m_{12} \circ (1\t \dd) = 0,
$$
and similarly for $[21]$. In other words $\dd$ is a derivation for $m_{12}$ and for $m_{21}$.
\M
For $n=3$, there are five maps $A^{\t 3} \to A$, denoted $m_{123}, m_{213}, m_{131},
 m_{312}, m_{321}$, corresponding to the five trees of $Y_3$. The five relations $(*)_y$ for $y$ a
tree of degree 3 are: 
$$\displaylines{
\dd \circ  m_{123} +  m_{123}\circ (\dd\t 1 \t 1 + 1\t \dd \t 1 + 1\t 1\t \dd ) = \hfill \cr
\hfill  m_{12} \circ ( m_{12}\t 1) -  m_{12} \circ (1\t m_{12})  , \cr
\dd \circ  m_{213} +  m_{213}\circ (\dd\t 1 \t 1 + 1\t \dd \t 1 + 1\t 1\t \dd ) = \hfill \cr
\hfill  m_{12} \circ ( m_{21}\t 1) -  m_{12} \circ (1\t m_{12})  , \cr
\dd \circ  m_{131} +  m_{131}\circ (\dd\t 1 \t 1 + 1\t \dd \t 1 + 1\t 1\t \dd ) = \hfill \cr
\hfill  m_{21} \circ ( m_{12}\t 1) -  m_{12} \circ (1\t m_{21})  , \cr
\dd \circ  m_{312} +  m_{312}\circ (\dd\t 1 \t 1 + 1\t \dd \t 1 + 1\t 1\t \dd ) = \hfill \cr
\hfill  m_{21} \circ ( m_{21}\t 1) -  m_{21} \circ (1\t m_{12})  , \cr
\dd \circ  m_{321} +  m_{321}\circ (\dd\t 1 \t 1 + 1\t \dd \t 1 + 1\t 1\t \dd ) = \hfill \cr
\hfill  m_{21} \circ ( m_{21}\t 1) -  m_{21} \circ (1\t m_{21})  . \cr
}$$

One observes that, as expected, if all the maps   $m_{ijk}$ are trivial (and also higher up), then
$m_{12}=\d\, ,\,  m_{21}=\g\, $, and the algebra is a graded dialgebra.

\vfill \eject



\noindent {\bf Appendix A. PLANAR BINARY TREES AND PERMUTATIONS.}
\BB

\N {\bf A.1. Planar binary trees. } A planar tree is {\it binary} if any vertex is trivalent.
We denote by
$Y_n$ the set of planar binary trees with $n$ vertices, that is with $n+1$ leaves (and one
root). Since we only use planar binary trees in this section we abbreviate it into tree (or
$n$-tree). The integer $n$ is called the degree of the tree. For any $y\in Y_n$ we label the $n+1$
leaves by
$\{0,1,\ldots ,n\}$ from left to right. We label the vertices by $\{1,\ldots ,n\}$ so that the $i$-th
vertex is in between the leaves
$i-1$ and $i$. In low dimension these sets are:
$$
Y_0 = \{ |\},\  Y_1= \{\  \arbreA \  \},\  Y_2=  \{\  
\arbreB ,\arbreC \  \},\   Y_3= \{\  \arbreun ,\arbredeux
,\arbretrois ,\arbrequatre ,\arbrecinq \  \}.
$$
The number of elements in $Y_n$ is $c_n={(2n)!\over n!(n+1)!}$,
 so ${\bf c}=(1,2,5,14,42,132, \dots )$ is the sequence of the Catalan numbers.
\M
We first introduce a notation for these trees as follows. The only element $\vert$ of $Y_0$
is denoted by $[0]$. The only element of $Y_1$ is denoted by $[1]$. 

The {\it grafting} of a $p$-tree $y_1$ and a $q$-tree $y_2$ is a
$(p+q+1)$-tree denoted by $y_1\vee y_2$ obtained by joining the roots of $y_1$ and $y_2$ and
 creating a new root from that vertex.
$$
y_1 \vee y_2 = \quad \arbreAbis
$$
Its name is written $[\ y_1\ p+q+1 \ y_2\ ] $ with the convention that all 0's are deleted
(except for the element in $Y_0$). For instance one has: $[0]\vee [0] = [1],\ [1]\vee [0] =
[12],\ [0]\vee [1] = [21],\ [1]\vee [1] = [131]$,
 and so on. So the names of the trees pictured above are, from left to right: 
$$[0]\ , [1]\ , [12]\ , [21]\ , [123]\ , [213]\ , [131]\ , [312]\ , [321]\ .$$
 This labelling has several
advantages. For instance if we draw the tree metrically, with the leaves regulary spaced,
and the lines at 45 degree angle, then the integers in the sequence are precisely the depth of
the successive vertices. 

Example 
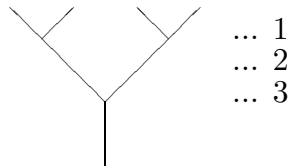
\M
$$\arbretroisbis
$$
\centerline {$[131]$}
\M

The orientation of the leaves can be read from the name as follows. Let $y = [a_1\ ...\
a_n]$. The $i$-th leaf is oriented SW-NE (resp. SE-NW) when $a_i<a_{i+1}$ (resp.
$a_i>a_{i+1})$.

The following is an inductive criterion to check whether a sequence of integers is the name
of a tree. 
\M

\noindent {\bf A.2. Proposition.} {\it A sequence of positive integers $[a_1\ ...\ a_n]$ is
the name of a tree if and only if it satisfies the following
conditions:

-- there is a unique
integer $p$ such that $a_p=n$,

-- the two sequences $a_1 ... a_{p-1}$ and $a_{p+1} ... a_n$ are either empty or name of
trees.}
\M

{\it Proof}. It suffices to remark that any tree $y$ is of the form $y_1\vee y_2$ for two
uniquely determined trees $y_1$ and $y_2$, whose degree is strictly smaller than the
degree of $y$. 
\hfill $\square$
\B
\noindent {\bf A.3.  From permutations to trees.} 
There is defined a surjective map 
$$\psi: S_n \proj Y_n$$
 as follows. The image of
 $\{1, ..., n\}$ under the permutation $\sigma$ is a sequence of positive integers 
$[\sigma _1\ ...\ \sigma _n]$. We 
convert it into the name of a tree by the following inductive rule. Replace the largest
integer in the interval $\sigma _1\ ...\ \sigma _n$, say $\sigma_p$,  by the length of the
interval (which is $n$ here). Then repeat the modification for the intervals $\sigma _1\ ...\
\sigma _{p-1}$ and $\sigma _{p+1}\ ...\ \sigma _n$, and so on, until each integer has been
modified. This gives the name of a tree (hence a tree), since it obviously satisfies the
above criterion. Let us perform this construction on an example, where the successive
modifications are underlined:
$$ \ss =[341652] \mapsto [341\underline
652]\mapsto [3\underline 31\underline 6\underline 22]\mapsto [\underline {131621}]= \psi
(\ss).$$

\N {\bf A.4.  Planar binary increasing trees.}
 Let us introduce a variation of trees: the {\it planar binary trees with levels} also
called {\it  increasing trees} in the literature. A tree with levels is an $n$-tree
together with a given level for each vertex. This level takes value in $\{1, ..., n\}$,
and we suppose that each vertex has a different level, and that the levels are increasing, that is
they respect the partial order structure of the tree (the level is the depth of the vertex).

Example: the following are two distinct increasing trees
$$
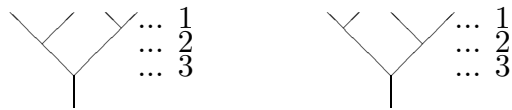
\arbretroisun \qquad \qquad\qquad \arbretroisdeux
$$

We denote by $\widetilde Y_n $ the set of increasing $n$-trees.

The following is a well-known result which was brought to my attention by Phil Hanlon [Ha].
\M

\noindent {\bf A.5. Proposition.} {\it The map which assigns a level to each vertex
determines a permutation. This gives a bijection $\widetilde Y_n \cong S_n$ between
increasing trees and permutations.}
\S

{\it Proof.} Label the vertices by their right leaf (i.e. vertex $i$ is in between the
leaf $i-1$ and the leaf $i$). Since each vertex has a different level it is clear that we get
a bijection. So any increasing tree gives rise to a permutation.

On the other hand a permutation gives rise to a tree under $\psi$, cf. A3. Labelling the level of the
$i$-th vertex by $\ss (i)$  gives an increasing tree.

It is immediately seen that the two constructions are inverse to each other.
$\square$
\B

Hence the map $\psi: S_n \r Y_n$ is the composite of the
bijection from the set of permutations to the set of increasing trees  with the
forgetful map (forgetting about the levels).
\S

In low dimension $\psi$ is given by:

$$\matrix{
\sigma &\ &[12]&[21]&\  &[123]&[213]&[132]&[231]&[312]&[321]&\hfill\cr
\psi(\sigma )&\ &[12]&[21]&\  &[123]&[213]&[131]&[131]&[312]&[321]&\hfill\cr
\cr
{\rm tree}&&\arbreB &\arbreC &&\arbreun &\arbredeux &\arbretrois &\arbretrois &\arbrequatre
&\arbrecinq \cr }
$$
\M

\M
\noindent {\bf A.6. The other choice $\psi '$ to code trees.} It is sometimes useful to code 
the vertices of a planar binary trees by {\it height} rather than by depth. The resulting
map $\psi ':S_n \to Y_n$ is related by the formula
$$
\psi ' (\ss) = \psi( \oo \ss \oo ^{-1}), \quad for \quad \oo=[n \cdots 2\ 1].
$$
For instance
$$
\psi ' ([12]) = \arbreC \qquad and \qquad \psi ' ([21]) = \arbreB\ .
$$
\M

\noindent {\bf A.7. Geometric interpretation of $\psi$}.
The map $\psi$ can be considered as the
restriction to the vertices of a cellular map between polytopes. Indeed, starting from the
set $<n> =
\{1, ... , n\}$ let us define a poset as follows. An element of the poset is an ordered
partition of
$<n>$. The element $X$ is less than the element $Y$ if $Y$ can be obtained from $X$ by
reuniting successive subsets. For instance $\{(35)(2)(14)(7)(6)\} < \{(235)(14)(67)\}$.
Note that the minimal elements are precisely the permutations. It can be shown that if we
exclude the trivial partition $\{(12 ... n)\}$ from the poset, then the geometric
realization is homeomorphic to the sphere $S^{n-2}$ (it is called the {\it permutohedron}).

Similarly, starting with planar trees with $n$ interior vertices one defines a poset as
follows: a tree $x$ is less than a tree
$y$ if $y$ can be obtained from $x$ by scratching some interior edges. Note
 that the minimal elements are precisely the planar binary trees.  It can be shown that if we
exclude the trivial tree with only one vertex from the poset, then the 
geometric realization is a polytope homeomorphic to the sphere $S^{n-2}$ (it is called the
Stasheff polytope, or the {\it associahedron}). The map
$\psi
$ described above can be obviously extended to a map of posets, and its geometric realization
is a homotopy equivalence of spaces, both homeomorphic to the sphere $S^{n-2}$. One can even compare
the associahedron with the hypercubes by looking at the orientation of the leaves, cf. [LR]. 

\vfill\eject


\def\P{{\cal P}}

\noindent {\bf Appendix B. ALGEBRAIC OPERADS}
\B

In this short appendix we briefly survey Koszul duality for algebraic operads as studied
by Ginzburg and Kapranov [GK] and Getzler-Jones [GJ].

The ground field $K$ is supposed to be of characteristic zero. The category of finite
dimensional vector spaces over $K$ is denoted by ${\bf Vect}$. When $V$ is graded its
suspension   $sV$ is such that $(sV)_n = V_{n-1}$.
\M

\N {\bf B.1. Algebraic operads}. For a given ``type of algebra" $\P$ (for instance
associative algebras, or Lie algebras, etc ...), let 
$\P (V)$ be the
 free algebra over the vector space $V$. One can view $\P$ as a functor from the
category {\bf Vect} to itself, which preserves filtered limits. The map $V\to \P (V) $
gives a natural transformation $\Id \to \P$.  The functor 
$\P$ is {\it analytical}. In characteristic zero it is equivalent (by Schur's lemma) to:
$\P$ is of the form
$$
\P(V) = \bigoplus _{n\geq  0}\P(n)\t _{S_n}V^{\t n},\leqno (B.1.1)$$
for some right
$S_n$-module  $\P(n)$.   

From the universal property of
the free algebra applied to $\Id: \P(V) \to \P (V)$ one gets a natural map $\P(\P(V)) \to
\P(V)$, that is a transformation of functors $\cc: \P\circ \P \to \P$.
By universality property of the free algebra functor, one sees that the operation $\cc$ is
associative and has a unit. 
\M
In order to make precise the notion of ``type of algebras", one axiomatizes the above
properties and puts up the following definition.
\M
By definition an {\it algebraic operad} over a characteristic zero field is an
analytical  functor
$\P: {\bf Vect} \to {\bf Vect}$, such that $\P(0) = 0$,   equipped with a natural
transformation of functor
$\cc :\P\circ \P \to \P$ which is associative and has a unit $1: \Id \to \P$. In other
words
$(\P,\cc, 1)$ is a {\it monad} in the tensor category $({\bf Funct}, \circ)$, where {\bf
Funct} is the category of analytical functors from {\bf Vect} to itself,
and $\circ$ stands for composition of functors.

By definition a {\it
$\P $-algebra} (that is an algebra over the operad $\P$) is a vector space $A$ equipped with
a map
$\cc_A:
\P (A)
\to A$ compatible with the composition $\cc$ in the following sense: the diagram
$$\matrix{
\P(\P(A)) & \buildrel {\P(\cc_A)} \over \longrightarrow & \P(A)\cr
&&\cr
\downarrow \cc(A)&&\downarrow \cc_A\cr
&&\cr
\P(A)&  \buildrel { \cc_A} \over \longrightarrow &A\cr}
$$
is commutative.

By writing $\P (V)$ and $\P\circ \P (V)$ in terms of the vector spaves $\P (n)$'s we see
that the operation
$\cc$ is determined by linear maps
$$
\P(n)\t \P(i_1)\t \cdots \t \P(i_n) \longrightarrow \P(i_1+\cdots +i_n)
$$
which satisfy some axioms deduced from the $S_n$-module structure of $\P(n)$ and from the
associativity of $\cc$ (cf. for instance [M]). 

From this point of view, a $\P$-algebra is determined by linear maps 
$$
\P(n)\t_{S_n}A^{\t n} \to A
$$
satisfying compatibility properties (cf. [M]). Observe that the space $\P (n)$ is
the space of all  operations that can be performed on $n$ variables.

The family of $S_n$-modules $\{\P(n)\}_{n\geq 1}$ is called an ${\bf S}$-module. There is an
obvious forgetful functor from operads to ${\bf S}$-modules. The left adjoint functor exists
and gives rise to the {\it free operad} over an ${\bf S}$-module (cf. for instance [BJT]).
\M 

We are only dealing here with {\it binary} operads, that is operads generated by operations
on two variables. More explicitly, let $E$ be an $S_2$-module (module of generating
operations) and let ${\cal T}(E)$ be the free operad on the ${\bf S}$-module
$$
(0,E,0,\ldots\quad ).
$$
The first components of ${\cal T}(E)$ are 
$$\eqalign{
{\cal T}(E)(1) &= K ,\cr
{\cal T}(E)(2) &= E ,\cr
{\cal T}(E)(3) &= \hbox{Ind}_{S_2}^{S_3} (E\t E ) = 3 E\t E,\cr
}$$
where the action of $S_2$ on  $E\t E$ is on the second factor only. Explicitly, ${\cal
T}(E)(3)$ is the space of all the operations on 3 variables that can be performed out of
the operations on 2 variables in $E$.

A binary operad is {\it quadratic} if it is the quotient of ${\cal T}(E)$ (for some
$S_2$-module $E$) by the ideal generated by some $S_3$-submodule $R$ of ${\cal T}(E)(3)$.
The associated operad is denoted $\P (E,R)$. 

For instance, for $\P = As$, one has
$E= K[S_2]$, the regular representation, and so  $\hbox{Ind}_{S_2}^{S_3} (E\t E ) =
K[S_3]\oplus K[S_3]$. Labelling the generators of the first summand by
$x_i(x_jx_k)$ and the generators of the second summand by $(x_ix_j)x_k$, the space $R$ is
the subspace generated by all the elements  $x_i(x_jx_k) - (x_ix_j)x_k$.

  All the operads
$ Com, Lie, As, Pois, Leib,  Zinb, Dias, Dend$ are binary
and quadratic. For all these cases the operadic notion of algebra coincides with the one we
started with.
\M

\N {\bf B.2. Dual operad.} For any right $S_n$-module $V$ we denote by $V^{\vee}$ the right
$S_n$-module $V^*\t (sgn)$, where $(sgn)$ is the one-dimensional signature representation.
Explicitly, if $^{\ss}f$ is given by  $^{\ss}f(x) = f(x^{\ss})$, then $f^{\ss} =
sgn(\ss) (^{\ss^{-1}}\! f)$.  The pairing between $V^{\vee}$ and $V$ given by:
$$<-,->: V^{\vee} \t V \longrightarrow  K ,\quad <f,x> = f(x),$$
is sign-invariant, that is $<f^{\ss}, x^{\ss}> = sgn (\ss) <f,x>$.

 To any quadratic binary operad $\P = \P (E,R)$ is associated
its {\it dual operad} $\P ^!$ defined as follows. Since $E$ is an $S_2$-module, so
is $E^{\vee}$. There is a canonical isomorphism
$${\cal T} (E^{\vee})(3) = {\rm Ind}\ _{S_2}^{S_3} (E^{\vee}\t E^{\vee})\cong 
 ({\rm Ind}\ _{S_2}^{S_3}(E\t E))^{\vee}= {\cal T}(E)(3)^{\vee} ,
$$
and one defines the orthogonal space of $R$ as
$$
R^{\bot}:= \Ker\  ({\cal T}(E^{\vee})(3) \to R^{\vee}).
$$
By definition one puts 
$$
\P^!:= \P(E^{\vee}, R^{\bot}).
$$

It is shown in [GK] that $As^! = As,\  Com^! = Lie,\  Lie ^! = Com$. The notation
$Zinb$ indicates that the  operad of Zinbiel algebras is precisely the dual, in
the above sense, of the operad of Leibniz algebras (cf. [L3]). It is proved in section 6
that the operad associated to dendriform   algebras, denoted $Dend$ is the dual of the operad
associated to dialgebras, as expected. It is a consequence of the following proposition,
which makes explicit the dual operad when
the operations bear no symmetry. This hypothesis says that the space of
operations is of the form $E=E'\t K[S_2]$ for some vector space $E'$. Then, the space of
operations on 3 variables is 
${\cal T}(E)(3) = (E'^{\t 2}\oplus E'^{\t 2}) \t K[S_3]$. For any $\xi \in E'$ we denote by
$(\xi)_1$ (resp. $(\xi)_2$ the element corresponding to an operation of the type $(-(--))$
(resp. $((--)-)\ )$. The first (resp. the second) component of $E'^{\t 2}\oplus E'^{\t 2}$
is made of the elements $(\xi)_1$ (resp. $(\xi)_2)$.
\M

\N {\bf B.3. Proposition.} {\it Let $\P$ be an operad whose generating operations
have no symmetry, in other words $\P = \P(E'\t K[S_2], R)$ for some vector space $E'$, and
some $S_3$-sub-module $R$ of $(E'^{\t 2}\oplus E'^{\t 2}) \t K[S_3]$.

The dual operad of $\P$ is 
$$
\P ^! = \P(E'\t K[S_2],R^{ann}),
$$
where $R^{ann}$ is the annihilator of $R$ for the scalar
product on $(E'^{\t 2}\oplus E'^{\t 2}) \t K[S_3]$ given by
$$\displaylines{
< \xi_1\t \ss , \xi_1\t \ss > = {\rm sgn}(\ss ),\cr
< \xi_2\t \ss , \xi_2\t \ss > = - {\rm sgn}(\ss ),\cr
}$$
all other scalar products are $0$, where $\xi_1$ (resp $\xi_2$) is a basis vector of
the first (resp. second) summand of $E'^{\t 2}\oplus E'^{\t 2}$ and $\ss \in S_3$.}

\S
\N {\it Proof.} Let $\phi: E \to E^{\vee}$ be the isomorphism of $S_2$-modules deduced
from the preferred basis of $E$. It induces an isomorphism of $S_3$-modules
$$
{\cal T} (\phi):  {\cal T}(E^{\vee})(3) \cong {\cal T}(E)(3).
$$
Hence, the natural evaluation map  ${\cal T}(E)(3)^{\vee} \t {\cal T}(E)(3)\to K$ gives a
scalar product  ${\cal T}(E)(3) \t {\cal T}(E)(3)\to K$.

Suppose that $E=E'\t K[S_2]$. Then ${\cal T}(E)(3) = (E'^{\t 2}\oplus E'^{\t 2}) \t
K[S_3]$ as mentioned above. In order to prove the Proposition, it suffices to prove it when
$E'$ is one-dimensional, generated by the unique product $(x_1x_2)$. The basis of $E$ is
$\{(x_1x_2), (x_2x_1)\}$. The isomorphism $\phi$ is given by $\phi ((x_1x_2)) = (x_1x_2)^*,
\phi ((x_2x_1)) = -(x_2x_1)^* $. The map ${\cal T}(\phi) $ sends a basis element to itself
(up to sign), so the scalar product is diagonal. Since it is sign-invariant, it
suffices to compute $<(x_1(x_2x_3)), (x_1(x_2x_3))>$ and  $<((x_1x_2)x_3), ((x_1x_2)x_3)>$.
With the choice at hand we find $+1$ in the first case and $-1$ in the second. \hfill
$\square$
\M

It is clear from this Proposition that $As^! = As$. Indeed, $E'$ is one dimensional
generated by $(\cdot \  \cdot )$, and $E'^{\t 2}\oplus E'^{\t 2}$ is 2-dimensional generated
by
$(\cdot (\cdot \ \cdot ))$ and $((\cdot \ \cdot )\cdot )$. The space $R$ is of the form
$R'\t K[S_3]$ because, in the associative relation, the variables stay in order. Since
$R'$ is determined by the equation $(\cdot (\cdot \ \cdot )) - ((\cdot \ \cdot )\cdot ) =0$,
it is immediate to check that it is its own annihilator.
\M

\N {\bf B.4. Homology and Koszul duality.} Let $\P$ be a quadratic binary operad and $\P
^!$ its dual operad. Let $A$ be a $\P$-algebra. There is defined a chain-complex $C_*^{\P}
(A)$:
$$
 \cdots \to  \P ^!(n)^{\vee}\t _{S_n} A^{\t n} \buildrel d \over
{\longrightarrow } \P ^!(n-1)^{\vee}\t _{S_{n-1}} A^{\t n-1} \to \cdots \to   \P
^!(1)^{\vee}\t A
$$
where the differential $d$ agrees, in low dimension, with the $\P$-algebra structure of $A$ 
$$
\cc_A(2): \P(2)\t A^{\t 2}\to A.
$$
In fact $d$ is characterized by this condition plus the fact that on the cofree coalgebra
$\P^{!*}(sA)$ it is a graded coderivation. The associated  homology
groups are denoted by
$H_n^{\P}(A)$, for $n\geq 1$. Taking the linear dual of $C^{\P}_*(A)$ over $K$
gives a cochain complex, which permits us to define the cohomology groups $H^n_{\P}(A)$.

If, for any vector space $V$, the groups $H_n^{\P}(\P (V))$ are trivial for $n>1$, then
the operad $\P$ is called a {\it Koszul} operad.
\S

One can check that the chain-complex $C^{\P}_*$ is 

-- the Hochschild complex (of nonunital algebras) for $\P = As$,

-- the Harrison complex (of nonunital commutative algebras) for  $\P = Com$,

-- the Chevalley-Eilenberg complex for  $\P = Lie$,

-- the chain-complex constructed in [L1] for  $\P = Leib$,

-- the chain-complex constructed in [Li2] for  $\P = Zinb$.
\M
Let us give some details about the check in the case of nonunital associative
 algebras. Since $As^! = As$, one has 
$$
C_n^{As}(A)= As^!(n)^{\vee}\t_{S_n}A^{\t n} = K[S_n]\t_{S_n}A^{\t n}
=A^{\t n} .$$
Since the lowest differential coincides with the product on $A$, the map $d_2:A^{\t
2}\to A$ is given by $d_2(x,y)=xy$. The coalgebra structure of $C_*^{As}(A)$
is given by ``deconcatenation", that is
$$
\DD (a_1, \ldots , a_n) = \sum_{i=0}^{n} (a_1, \ldots , a_i) \t (a_{i+1}, \ldots
, a_n).
$$
Since $d$ is a graded coderivation, we have on $A^{\t 3}$
$$\eqalign{
\DD d_3(a_1, a_2, a_3) &= (d_2\t 1 + 1\t d_2)\DD (a_1, a_2, a_3)\cr
&= (d_2\t 1 + 1\t d_2)( (a_1,a_2)\t a_3 + a_1\t (a_2,a_3))\cr
&=  (a_1a_2)\t a_3 - a_1\t (a_2a_3).\cr}
$$
Hence we obtain $d_3(a_1, a_2, a_3)=(a_1a_2, a_3) - (a_1,a_2a_3)$.

More generally, the same kind of computation shows that
$$
d_n(a_1, \ldots , a_n) = \sum_{i=1}^{n-1}(-1)^{i-1}(a_1, \ldots , a_ia_{i+1},
\ldots , a_n).
$$
So $(C_*^{\P}(A), d)$ is precisely the Hochschild complex for non-unital algebras
(also called the $b'$-complex in the literature, cf. [L0] for instance).
\M

\N {\bf B.5. Properties of the operad dualization.}
Here are some properties of the dualization of operads.
\M

\N (a) {\it Lie algebra property}. Let $A$ be a $\P$-algebra and $B$ be a $\P ^!$-algebra. The
following bracket makes    $A\t B$ into a Lie algebra :
$$
[a\t b, a'\t b']:= \sum \big( \mu (a,a')\t \mu^{\vee}(b,b') - \mu (a',a)\t
\mu^{\vee}(b',b)\big),
$$ where the sum is over a basis $\{\mu\}$ of the binary operations of $\P$, 
$\{\mu^{\vee}\}$ being the dual basis in $\P^!$.
\M

\N (b) {\it  Koszulness.} $\P^{!!} = \P $. If $\P$ is Koszul, then so is $\P ^!$.
\M

\N (c) {\it Poincar\'e series.} Define the Poincar\'e series of $\P$ as
$$
g_{\P}(x):= \sum_{n\geq 1} (-1)^n \dim\ \P (n) {x^n \over n!}.
$$
If $\P$ is Koszul, then the following formula holds:
$$
g_{\P}(g_{\P ^!}(x)) = x.
$$
For instance $g_{As}(x) = {-x \over 1+x},\  g_{Comm}(x) = {\rm exp}(-x)-1,\ g_{Lie}(x) =
-{\rm log} (1+x)$.
\M

\N (d) {\it Multiplicative structure.} For any $\P$-algebra $A$ the homology groups
$H_*^{\P}(A)$ form a graded $\P^!$-coalgebra. Equivalently,  $H^*_{\P}(A)$ is a
graded $\P ^!$-algebra.
\M
\N (e) {\it Functoriality.} Let $\phi:\P \to \Q$ be a morphism of operads, inducing
 a functor
($\Q $-algebras) $\to$ ($\P $-algebras), $A\mapsto A_{\P}$ between the categories
 of algebras. For any $\Q$-algebra $A$ there is a well-defined graded morphism
$$
\phi_*: H_*^{\P}(A_{\P}) \to  H_*^{\Q}(A),
$$
which is induced, at the complex level, by the unique morphism of $\Q^!$-algebra
$$
\Q^!(V) \to \P^!(V)_{\Q^!},
$$
which commutes with the embeddings of the vector space $V$.
\M

\N (f) {\it Quillen homology.} If $\P$ is Koszul, then the homology theory $H_*^{\P}$ for
$\P$-algebras, as defined above, coincides with the Quillen homology with trivial
coefficients (cf. [Q], [Li3], [FM]). This is a consequence of the existence of a
quasi-isomorphism ${\cal B}(\P^!)^*\to \P$ (see below).
\M

\N {\bf B.6. Homotopy algebras over an operad.} A {\it quasi-free resolution}  ${\cal Q}\to
\P$ of the operad $\P$ is a differential graded operad ${\cal Q}$ which is a free operad over
some graded {\bf S}-module $V$, and such that, for all $n$, ${\cal Q}(n)$ is a chain-complex
whose homology is trivial, except $H_0$ which is equal to $\P (n)$. The isomorphism
$H_0({\cal Q}) \cong \P$ is supposed to be an isomorphism of operads.  Observe that one
does not require ${\cal Q}$ to be free in the category of differential graded operads. By
definition a  {\it homotopy $\P$-algebra} is a ${\cal Q}$-algebra for some quasi-free
resolution ${\cal Q}$ of $\P$.

A quasi-free resolution ${\cal Q}$ over $V$, with augmentation ideal $\bar {\cal Q}$  is
called {\it minimal} when the differential $d$ is {\it quadratic}, that is, satisfies $d(V)
\subset  {\bar {\cal Q}}\cdot {\bar {\cal Q}}$.
For a given $\P$ such a minimal model always exists in characteristic zero  and is unique
up to homotopy. 

By definition a {\it strong homotopy $\P$-algebra} is a ${\cal Q}$-algebra
for the minimal model ${\cal Q}$ of $\P$.

 If $\P$ is Koszul, then there is a way of constructing the minimal model as follows. Let 
$\P ^!$ be the Koszul dual of $\P$ and let ${\cal B}(\P^!)$ be the bar-construction over 
$\P^!$. It is the cofree cooperad ${\cal T}'(s\bar \P ^!)$   equipped with the unique 
coderivation $d$ which is 0 on $\P ^!$ and coincides with the composition on $\P ^! \circ
\P ^!$. It is immediate that $d^2=0$, and hence ${\cal B}(\P^!):= ({\cal T}'(s\bar \P ^!),d)
$ is a complex called the {\it bar-complex} of $\P ^!$. It is shown in [GK] (cf. also
[GJ]), that, if the operad $\P$ is Koszul, then the differential graded operad ${\cal
B}(\P^!)^*$ is a resolution of
$\P$. Since it is quasi-free and since the differential $m$ of ${\cal B}(\P^!)^*$ is
quadratic,
${\cal B}(\P^!)^*$ is the minimal model of $\P$.
\M

\N {\bf B.7. Theorem.} {\it For a Koszul operad $\P$,  strong homotopy $\P$-algebras are
equivalent to  ${\cal B}(\P^!)^*$-algebras.}
\M
\N {\bf B.8. Strong homotopy associative algebras.} The associative operad $\P= {\bf As}$ is
Koszul and self-dual. Since ${\bf As}(n)$ is one-dimensional as a free $S_n$-module, it
gives rise to a generating operation $m_n: A^{\t n} \to A$  of the strong  homotopy
associative algebra $A$. The condition $m\circ m=0$ gives, for each $n\geq 1$
$$
\sum _{k+l=n-1} m_l\circ m_k = 0.
$$
Hence a strong homotopy associative algebra is exactly what is called an
$A_{\infty}$-algebra in the literature (cf. [St]).

\vfill \eject



\centerline{\bf References}
\M
\noindent [AG] M. AYMON, and P.-P. GRIVEL, {\it Un th\'eor\`eme de Poincar\'e-Birkhoff-Witt pour
les alg\`ebres de  Leibniz}, Comm. Alg. (to appear).
\smallskip

\noindent [Ak] F. AKMAN, {\it On some generalizations of Batalin-Vilkovisky algebras}, J. Pure 
Applied Alg. 120 (1997), 105--141.
\smallskip

\noindent [BJT] H.J. BAUES, M. JIBLADZE, A. TONKS, {\it 
Cohomology of monoids in monoidal categories}. 
Operads: Proceedings of Renaissance Conferences (Hartford, CT/Luminy, 1995), 137--165, 
Contemp. Math., 202, 
Amer. Math. Soc., Providence, RI, 1997. 

\smallskip
\noindent [FM] T.F. FOX and M. MARKL, {\it  Distributive laws, bialgebras, and cohomology}.
Operads: Proceedings of Renaissance Conferences (Hartford, CT/Luminy, 1995), 167--205, Contemp.
Math., 202, Amer. Math. Soc., Providence, RI, 1997.
\smallskip

\noindent [F1] A. FRABETTI,   {\it Dialgebra homology of associative
algebras},
 C. R. Acad. Sci. Paris 325 (1997), 135-140. 
\smallskip 

\noindent [F2] A. FRABETTI, {\it Leibniz homology of dialgebras of matrices},
 J. Pure Appl. Alg. 129 (1998), 123-141.
\smallskip 

\noindent [F3] A. FRABETTI, {\it Simplicial properties of the set of planar binary trees.}
 J. Alg. Combinatorics, to appear.
\smallskip 

\noindent [F4] A. FRABETTI, {\it Dialgebra (co)homology with coefficients.}
\smallskip

\noindent [GJ] E. GETZLER and J.D.S. JONES, {\it Operads, homotopy algebra, and iterated integrals
for double loop spaces}, electronic prepublication \hfill
\break {\tt [ http://xxx.lanl.gov/
abs/hep-th/9403055 ]}.
\smallskip

\noindent [GK] V.  GINZBURG,  M.M. KAPRANOV, {\it Koszul duality for operads}, Duke Math. J. 
76 (1994), 203--272.  
\smallskip

\noindent [Go] F. GOICHOT, {\it Un th\'eor\`eme de Milnor-Moore pour les alg\`ebres de Leibniz},
\smallskip

\noindent [Ha] P. HANLON, {\it The fixed point partition lattices}, Pacific J. Math. 96, no.2,
(1981), 319--341.
\smallskip

\noindent [Hi] P. HIGGINS, Thesis, Oxford (UK).
\smallskip

\noindent [I] H.N. INASSARIDZE, {\it Homotopy of pseudosimplicial groups, nonabelian derived
functors, and algebraic
$K$-theory.} (Russian) Mat. Sb. (N.S.) 98(140) (1975), no. 3(11), 339--362, 495.
\smallskip

\noindent [KS] Y. KOSMANN-SCHWARZBACH, {\it  From Poisson algebras to Gerstenhaber
algebras}. Ann. Inst. Fourier
 46 (1996), no. 5, 1243--1274.
\smallskip

\noindent [Kub] F. KUBO, {\it
Finite-dimensional non-commutative Poisson algebras.} 
J. Pure Appl. Algebra 113 (1996), no. 3, 307--314. 
\smallskip

\noindent [Kur] R. KURDIANI, {\it Cohomology of Lie algebras in the tensor category of linear maps.}
Comm. Alg.27 (1999), 5033--5048.
\smallskip

\noindent [Li1] M. LIVERNET, {\it Homotopie rationnelle des alg\`ebres de Leibniz},
C. R. Acad. Sci. Paris  325 (1997), 819--823.
\smallskip

\noindent [Li2] M. LIVERNET, {\it Rational homotopy of Leibniz algebras}, Manuscripta
Mathematica 96 (1998), 295--315. 
\smallskip

\noindent [Li3] M. LIVERNET, {\it On a plus-construction for algebras over an operad},
$K$-theory Journal 18 (1999), 317--337.
\smallskip

\noindent [L0] J.-L.  LODAY,  ``Cyclic homology". Second edition. Grund. math.
 Wiss., {\bf 301}. Springer-Verlag, Berlin, 1998.
\smallskip

\noindent [L1] J.-L.  LODAY, {\it Une version non commutative des alg\`ebres de Lie : les
alg\`ebres de Leibniz}, Ens. Math. 39 (1993), 269--293.
\smallskip 

\noindent [L2] J.-L. LODAY, {\it Alg\`ebres ayant deux op\'erations associatives (dig\`ebres)},
  C. R. Acad. Sci. Paris 321 (1995), 141-146. 
\smallskip 

\noindent [L3] J.-L. LODAY, {\it Cup-product for Leibniz cohomology and dual Leibniz
algebras}, Math. Scand. 77 (1995), 189--196.
\smallskip

\noindent [L4] J.-L. LODAY,  {\it Overview on Leibniz algebras, dialgebras and their homology}.
 Fields Inst. Commun. 17 (1997), 91--102.
\smallskip

\noindent [LP1] J.-L. LODAY, and  T.  PIRASHVILI, {\it Universal enveloping algebras of
Leibniz algebras and (co)homology}, Math. Ann. 296 (1993), 139--158.
\smallskip

\noindent [LP2] J.-L. LODAY, and  T.  PIRASHVILI, {\it  The tensor category of linear maps 
and Leibniz algebras}. Georgian Math. J. 5 (1998), 263-276.
\smallskip

\noindent [LR] J.-L. LODAY, and M.O. RONCO, {\it  Hopf algebra of the planar binary trees}.
Adv. in Maths 139 (1998), 293--309.
\smallskip

\noindent [Lo] J. LODDER, {\it Leibniz homology and the Hilton-Milnor theorem}, Topology 36 (1997),
729--743. 
\smallskip

\noindent [M]  J. P. MAY, {\it Definitions: operads, algebras and modules}. Operads: Proceedings of Renaissance Conferences
(Hartford, CT/Luminy, 1995), 1--7, Contemp. Math., 202, Amer. Math. Soc., Providence, RI, 1997.
\smallskip

\noindent [Q] D. QUILLEN, {\it On the (co-) homology of commutative rings}. 1970 Applications of Categorical Algebra (Proc.
Sympos. Pure Math., Vol. XVII, New York, 1968) pp. 65--87.
\smallskip

\noindent [R] M.O. RONCO, {\it A Milnor-Moore theorem for dendriform Hopf algebras}. C. R. Acad.
Sci. Paris S\'er. I Math. (2000), to appear.
\smallskip

\noindent [S] J.D. STASHEFF, {\it Homotopy associativity of $H$-spaces}. I, II. Trans. Amer.
Math. Soc. 108 (1963), 275-292; ibid. 108 1963 293--312.


\BB

 Institut de Recherche Math\'ematique Avanc\'ee

    CNRS et Universit\'e Louis Pasteur

    7 rue Ren\'e-Descartes

    67084 Strasbourg Cedex, France

    E-mail : loday@math.u-strasbg.fr
\BB
\hfill March 19,  1999

\hfill updated February 6, 2001
\end